\documentclass[11pt]{amsart}
\usepackage{graphicx}
\usepackage[all]{xy}
\usepackage{amssymb, amsmath}
\usepackage{setspace}
\usepackage{eufrak}
\usepackage{amsfonts}
\usepackage{float}
\usepackage{amsmath}

\def\qed{\hfill $\Box$}
\def\proof{\noindent {\sl Proof} :\;  }

\def\dis{\displaystyle}

\newcommand{\A}{\mathcal{A}}
\newcommand{\K}{\mathcal{K}}

\newcommand{\Gcal}{\mathcal{G}}
\newcommand{\E}{\mathcal{E}}
\newcommand{\F}{\mathcal{F}}
\newcommand{\Ost}{\mathcal{O}}
\newcommand{\Proj}{\mathbb{P}}
\newcommand{\C}{\mathbb{C}}
\newcommand{\Z}{\mathbb{Z}}

\newcommand{\Q}{\mathbb{Q}}
\newcommand{\R}{\mathbb{R}}

\newcommand\jeden {1\hskip-3.5pt1}

\newcommand\codim{ \mbox{\rm codim}\, }
\newcommand\Dual{\mbox{\rm Dual}\, }

\newcommand\Image{\mbox{\rm Im} }
\def\Diff{\mbox{\rm Diff}}
\def\Hom{{\rm Hom}}
\def\qed{\hfill $\Box$}
\def\proof{\noindent {\sl Proof} :\;  }

\def\rd{\partial}

\def\imageconst{\alpha_{\mbox{\rm \tiny image}}}
\def\disconst{\alpha_{\mbox{\rm \tiny dis}}}
\def\SM{{\rm SM}}
\def\Ma{{\rm Ma}}

\def\bv{\mbox{\boldmath $v$}}

\def\b0{\mbox{\boldmath $0$}}

\newtheorem{thm}{\bf Theorem}[section]
\newtheorem{"thm"}[thm]{\bf `Theorem'}
\newtheorem{cor}[thm]{\bf Corollary} 
\newtheorem{lem}[thm]{\bf Lemma} 
\newtheorem{prop}[thm]{\bf Proposition} 
\newtheorem{definition}[thm]{\bf Definition} 
 
\newtheorem{rem}[thm]{\bf Remark} 
\newtheorem{exam}[thm]{\bf Example} 
\newtheorem{conj}[thm]{\bf Conjecture} 

\begin{document}

\title[ Higher Thom polynomials]
{Singularities of Maps  and Characteristic Classes } 
\author[T.~Ohmoto]{Toru Ohmoto}
\address[T.~Ohmoto]{Department of Mathematics, 
Faculty of Science,  Hokkaido University,
Sapporo 060-0810, Japan}
\email{ohmoto@math.sci.hokudai.ac.jp}
\dedicatory{Dedicated to Professor Shyuichi Izumiya on the occasion of his 60th birthday.}
%
%
\maketitle
%
%
%

\tableofcontents


\section{Introduction} 

In classical algebraic geometry, 
numerical characters of complex projective varieties 
were extensively studied by means of enumerating 
singular points of naturally associated algebraic maps, 
e.g., the degree of loci of ramification, polar, 
multiple points, inflections ... and so on. 
A modern unified approach to such enumerative problems  
is {\it the theory of  Thom polynomials}
 based on 
the classification of mono and multi-singularities of maps. 
In this lecture  we introduce a new branch of the theory, 
in which we replace  
counting singular points by computing (weighted) Euler characteristics. 
This theory leads to a number of generalizations of classical enumerative formulas, 
while we here focus on an application to the vanishing topology 
of $\A$-finite map-germs. 

A simple toy example is the  {\it  Riemann-Hurwitz formula}: 
Let $f: M \to N$ be a surjective holomorphic map between compact complex curves. 
To each point of $M$ the multiplicity $\mu=\mu(f)$ is assigned so that 
the germ of $f$ at the point is written as $z \mapsto z^{\mu+1}+\cdots$. 
The classical formula says  that 
the number of critical points taking account of multiplicities $\mu$  
measures the different between 
the topological Euler characteristics $\chi$ of $M$ and $N$, 
that is written in a slightly modern form as follows: 
\begin{eqnarray*}
{\displaystyle \int_M \mu(f) \; d\chi} &=&\deg f\cdot \chi(N) - \chi(M)\\ 
&=& c_1(TN) \frown f_*[M] - c_1(TM) \frown [M]\\  
&=&  c_1(f^*TN-TM) \frown [M]. 
\end{eqnarray*}
Here appear major characters playing in this mini-course: 
\begin{itemize}
\item 
$c_i$ stands for the {\it Chern class} of vector bundles and 
$[ -]$ is the fundamental cycle in classical intersection theory (Section \ref{Chern}); 
\item 
$\int_M$ is the {\it integral of constructible functions}, 
which will soon be replaced by 
the {\it Chern-Schwartz-MacPherson class} (CSM class)   (Section \ref{CSM}); 
 \item 
 $tp(A_1)=c_1(f^*TN-TM)$,  the simplest {\it Thom polynomial} 
for $A_1$-singularity of equidimensional maps (Section \ref{thom}). 
\end{itemize}
The emphasis is that 
{\it integrating} local invariants of singularities of maps 
provides global invariants associated to maps, 
and conversely, 
{\it localizing} global invariants to a critical point 
(via torus-action) 
 gives local invariants at that point. 
Our main goal is to present 
a certain framework for generalizing this picture, 
based on the well-established classification theory of map-germs 
(the Thom-Mather theory) and characteristic classes for singular varieties 
(Chern-Mather and Chern-Schwartz-MacPherson classes 
and (singular) Todd class etc). 
We also show the effectivity of our approach 
by giving a number of actual computations in concrete examples. 

We works in the complex analytic/algebraic context throughout, 
however, almost all parts can suitably be repeated  over  
algebraically closed field in characteristic zero. 

The organization of this note is as follows. 

We begin with basic materials: 
In \S 2, some required knowledge in classification theory of map-germs and 
 classical intersection theory are briefly summarized. 

A quick introduction to the CSM class is given in \S 3. 
In particular this section contains a digest from \cite{Ohmoto06} 
about equivariant (co)homology,  
the algebraic Borel construction and 
the theory of equivariant CSM class: 
Theorem \ref{segre2} is the foundation of this lecture.  

The main body is \S 4. 
Given a stable singularity type $\eta$ 
of holomorphic map-germs from $\C^m$ to $\C^n$, 
the {\it Thom polynomial} $tp(\eta)$ is by definition 
a universal polynomial in 
the quotient Chern classes $c_i(f)=c_i(f^*TN-TM)$ which 
expresses the fundamental class of the closure of 
$$\eta(f):= \{\; x \in M \; | \; 
\mbox{the germ $f$ at $x$ is of type $\eta$} \; \}$$ 
for any stable maps $f:M \to N$ (Theorem \ref{tp}): 
$$\Dual [\overline{\eta}(f)] = tp(\eta)(c(f)).$$
Obviously, in case that the codimension of $\eta$ is equal to $\dim M$, 
$tp(\eta)$ for $f$ counts the number of $\eta$-singular points. 
Such universal polynomial expression  
can be considered for not only the fundamental class but also 
other certain invariants of the prescribed singular locus of maps. 
We then focus on the topological Euler characteristics - 
the {\it higher Thom polynomial} $tp^{\SM}(\overline{\eta})$ 
is introduced so that 
it universally expresses the CSM class  
of the $\eta$-type singular point locus $\overline{\eta}(f)$ (Theorem \ref{tp^SM}):  
$$\Dual c^{\SM}(\overline{\eta}(f))=c(TM) \cdot tp^{\SM}(\overline{\eta}).$$
In particular, the degree of the right hand side computes    
the Euler characteristics $\chi(\overline{\eta}(f))$. 
Here $tp^{\SM}(\overline{\eta})$ is a power series in $c_i=c_i(f)$ 
whose leading term is just the Thom polynomial $tp(\eta)$.  
To determine those polynomials, 
there is an effective method, which is described for a typical example in \S 4.3. 
We also discuss (higher) Thom polynomials for multi-singularities. 

Indeed we give several universal formulas for (weighted) Euler characteristics of 
singular loci in the source and the target; 
for instance, we show in Proposition \ref{complexIM} that  
for a closed singular surface $X$ in a projective $3$-fold $N$ 
having ordinary singularities, 
i.e., crosscaps ($A_1$) and normal crossings (double and triple points), 
and for its normalization 
$f: M \to X \subset N$,  it holds that 
$$\chi(X)=\frac{1}{6}  {\dis\int_M} 
\left(
\begin{array}{l}
3c_1(TM) c_1+6c_2(TM)-3c_1(TM) s_0\\
-c_1^2-c_2-2c_1s_0+s_0^2+2s_1
\end{array}
 \right) $$ 
where 
$c_i=c_i(f^*TN-TM), \; s_0=f^*f_*(1), \; s_1=f^*f_*(c_1)$. 
This is part of 
our more general formulas (Theorems \ref{imageSSM}, \ref{DisSSM}). 

We remark that 
as particular examples, 
applying  these (higher) Thom polynomials of multi-singularities 
to certain maps in projective algebraic geometry, 
e.g., normalizations of projective surfaces with ordinary singularities,   
leads us to rediscover a number of classical formulas in 19 century due to 
Salmon, Caylay, Zeuthen, Enriques, Baker, .... 
and actually it gives suitable generalizations: 
In particular, the computations on higher Thom polynomials involve 
the `exclusion-inclusion principle' among multi-singularity loci, 
that is quite similar to some typical argument in the classical works of those pioneers.

\S 5 and \S 6 are devoted to our main application.  
The purpose is to present a new method for studying {\it the vanishing topology 
of finitely determined weighted homogeneous map-germs} 
by localizing  (higher) Thom polynomials via $\C^*$-action: 
We exhibit a bunch of  numerical computations of 
\begin{itemize}
\item the number of stable singularities appearing in generic perturbation 
($0$-stable invariants) 
\item image and discriminant Milnor numbers 
\end{itemize}
for such map-germs of {\it any corank} in low dimensions. 
Our method can provide general formulas in terms of weights and degrees. 
Those are really new: 
In fact there has not been known any effective method  
for computing such invariants of germs without corank condition.

In this lecture note, mainly we deal with maps $f: M \to N$ of  
non-negative {\it relative-codimension} $\kappa:=\dim N-\dim M \ge 0$, 
and the negative codimensional case will be considered in another paper.  

The author would like to appreciate the organizers of the 12th Workshop 
of Real and Complex Singularities in Sao Carlos (2012) 
for their good organizing and kind hospitality, 
and also thank Maxim Kazarian for his comments  (\S \ref{multising_SSM}). 
This is partly supported by the JSPS grant no.24340007.

\section{Preliminaries}
\subsection{Basics in $\A$ and $\K$-classifications of map-germs} \label{singularity}
We describe some basic notions in the Thom-Mather theory, 
see , e.g., \cite{Arnold, MondLect, Wall}. 

Let $\Ost_m$ be the local ring of holomorphic function germs $\C^m,0 \to \C$ 
with the maximal ideal $\mathfrak{m}_m=\{h \in \Ost_m, f(0)=0\}$. 
Put $\E(m,n)$ to be the $\Ost_m$-module of all homolorphic map-germs 
$C^m,0 \to \C^n$, and also put 
$$\E_0(m,n)=\{\; f: \C^m, 0 \to \C^n,0 \;\; \mbox{holomorphic}\; \}=\mathfrak{m}_m \E(m,n).$$

\

\noindent
{\bf Equivalence.} 
The group of biholomorphic germs $\C^m, 0 \to \C^m,0$ 
is denoted by $\Diff (\C^m,0)$ (abusing the notation $\Diff$). 
There are two different kinds of natural equivalence relations on map-germs: 
\begin{itemize}
\item {\bf $\A$-classification} (right-left equivalence) 
classifies map-germs up to isomorphisms of source and target. 
The {\it right-left group} $\A \; (=\A_{m,n})$ is 
the direct product $\Diff (\C^m,0)\times \Diff (\C^n,0)$, 
which acts on $\E_0(m,n)$ by 
$$(\sigma, \tau).f:=\tau \circ f \circ \sigma^{-1}.$$
\item {\bf $\K$-classification} (contact equivalence) 
classifies up to the isomorphisms of source 
the zero locus $f^{-1}(0)$ as a scheme, i.e., classifies the ideal 
$$f^*\mathfrak{m}_n :=\langle f_1, \cdots, f_n\rangle_{\Ost_m}\; \subset \Ost_m$$ 
generated by the component functions of $f$; 
In other words, the $\K$-equivalence measures the tangency of 
the {\it graphs} $y=f(x)$ and $y=0$  in $\C^m \times \C^n$. 
The {\it contact  group} $\K \; (=\K_{m,n})$ consists of pairs $(\sigma, \Phi)$ 
of $\sigma \in \Diff (\C^m,0)$ and $\Phi: \C^m,0 \to GL(n,\C)$, 
which acts on $\E_0(m,n)$ by 
$$((\sigma, \Phi).f)(x)=\Phi(x)f(\sigma(x)).$$
\item 
If $f\sim_\A g$, then $f \sim_\K g$, i.e.,  
$\A.f \subset \K.f.$
\end{itemize}

\begin{exam}{\rm 
$f=(x^3+yx, y)$ and $g=(x^3, y)$ in $\E_0(2,2)$ 
are $\K$-equivalent but not $\A$-equivalent, so $\A.f \not= \K.f$.    
The $\A$-class of $f=(x^3+yx, y)$ is called  
an ordinary {\it cusp} or stable $A_2$-singularity. 
The discriminant (=singular value curves on the plane) is depicted in Fig. \ref{cusp}.  
\begin{figure}
\includegraphics[clip, width=2.5cm]{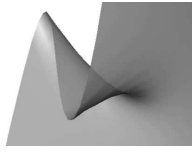}
\caption{\small Cusp ($A_2$) arises in a generic projection of a surface to the plane}
\label{cusp}
\end{figure} 
}
\end{exam}

\

\noindent
{\bf  Tangent spaces.} 
Let $f \in \E_0(m,n)$. 
An {\it infinitesimal deformation} of $f$ is 
a vector field-germ along $f$
$$v: \C^m, 0 \to T\C^n, \quad p \mapsto v(p) \in T_{f(p)}\C^n.$$
The space of infinitesimal deformations 
 is regarded as the `tangent space' of $\E(m,n)$ at $f$, and 
is denoted by 
$$\textstyle \theta(f) = \oplus_{i=1}^n \Ost_m \frac{\rd}{\rd y_i}.$$
Note that $\theta(f)$ admits two different module structures 
via multiplications with source functions in $\Ost_m$ and 
target functions in $\Ost_n$ through $f^*$. 
The subspace of infinitesimal deformations vanishing at the origin 
is just $\mathfrak{m}_m\theta(f)$,  regarded 
as the tangent space of $\E_0(m,n)$ at $f$. 

For the identity map $id_m$ of $\C^m$, 
$$\textstyle 
\theta_m:=\theta(id_m)= \oplus_{i=1}^m \Ost_m \frac{\rd}{\rd x_j}$$ 
is the space of germs of vector fields on $\C^m$ at the origin, 
in other words, the space of 
infinitesimal deformations of coordinate changes of $\C^m$ 
not necessarily preserving the origin.  
Instead, $\mathfrak{m}_m\theta_m$ is 
the space of infinitesimal deformations of coordinate changes 
preserving the origin. 
We set an $\Ost_m$-module homomorphism 
$tf:\theta_m \to \theta (f)$ and 
an $\Ost_n$-module homomorphism 
$\omega f:\theta_n \to \theta(f)$ by 
\begin{eqnarray*}
tf &: &\textstyle 
v=\sum v_j(x) \frac{\rd}{\rd x_j}\;\;  \longmapsto \;\; 
df(v)=\sum \frac{\rd f_i}{\rd x_j}(x)v_j(x) \frac{\rd}{\rd y_i}, \\
\omega f&: &\textstyle
 w=\sum w_i(y) \frac{\rd}{\rd y_i} \;\; \longmapsto \;\; 
 w\circ f=\sum w_i(f(x))\frac{\rd}{\rd y_i}.
\end{eqnarray*}
Then the {\it  tangent spaces} of $\A$ and $\K$-orbits of $f$ 
in $\mathfrak{m}_m\theta(f)$ 
and the {\it extended} tangent spaces in $\theta(f)$ are defined as follows: 
\begin{eqnarray*}
T\A.f&:=&tf(\mathfrak{m}_m\theta_m)+\omega f(\mathfrak{m}_n\theta_n),\\
T\K.f&:=&tf(\mathfrak{m}_m\theta_m)+f^*\mathfrak{m}_n\theta (f), \\
T\A_e.f&:=&tf(\theta_m)+\omega f(\theta_n), \\
T\K_e.f&:=&tf(\theta_m)+f^*\mathfrak{m}_n\theta (f). 
\end{eqnarray*}

\

\noindent
{\bf Determinacy.} Let $\Gcal=\A$ or $\K$. 
A map-germ $f$ is {\it $\Gcal$-finitely determined} if 
there is some $k$ so that if $j^kg(0)=j^kf(0)$ then $g \sim_\Gcal f$. 
Finite determinacy is equivalent to that 
the orbit $\Gcal.f$ has finite codimension, i.e., 
$\dim_\C \, \mathfrak{m}_m\theta(f)/T\Gcal.f < \infty$ 
($\Leftrightarrow \dim_\C \, \theta(f)/T\Gcal_e.f < \infty$).  
Then, the  process for $\Gcal$-classification of finitely determined map-germs 
 is reduced to the level of jets (Taylor polynomials): 
We may replace $\E_0(m,n)$ and $\Gcal$ by 
jet spaces $J^k(m,n)$and  $J^k\Gcal$, respectively, 
which are finite dimensional and the action is algebraic. 

\

\noindent
{\bf  Stability.} 
 $f: \C^m, 0 \to \C^n,0$ is  a {\it stable germ} if 
 any infinitesimal deformation of $f$ is recovered by 
 some infinitesimal deformations of source and target coordinate changes 
 (not necessarily preserving the origin), that is, 
$$\theta(f)=T\A_e.f.$$
By the Malgrange preparation theorem this condition is equivalent to that 
$$\textstyle
\theta(f)=T\K_e.f+\oplus_{i=1}^n \C \frac{\rd}{\rd y_i}. $$ 

It is known that 
$f \sim_\K g$ if and only if  $f \sim_\A g$ for stable germs $f, g$. 
Namely,  for a stable germ $f$, 
$$\A.f = \{\mbox{\small Stable germs}\} \cap \K.f.$$

\

\noindent
{\bf  Jet extension.}  
Intuitively, a stable germ $f$ means that 
for any small perturbation of any representative $f: U \to V$, 
still the same type singularity remains at some  point nearby the origin. 
This is justified by the transversality of jet extension. 
A representative $f:U \to V$ produces a map 
$$\bar{f}: U \to V \times \E_0(m,n), \; p \mapsto \mbox{germ of $f(x+p)$ at $x=0$},$$
then the image of the derivative of this map at $0$ is just 
the linear subspace 
$tf(T_0U)$ of $\theta(f)=\omega f(T_0V) \oplus \mathfrak{m}_m\theta(f)$. 
Note that 
$$T\A_e.f= tf(T_0U)+\omega f(T_0V)+ T\A.f$$ 
and $\omega f(T_0V)+ T\A.f$ is regarded as 
the tangent space of $V \times \A.f$. Thus 
we have 
\begin{center}
$\theta(f)=T\A_e.f$ $\Longleftrightarrow$ $\bar{f}$ is transverse to $V\times \A.f$ at $0$ 
\end{center}
Also this is equivalent to that 
$\bar{f}$ is transverse to $V\times \K.f$ at $0$, 
using the interpretation of the stability in terms of $T\K_e.f$. 

Precisely saying, we should state the transversality (the right hand side) on the level of jets: 
Let $J(TM, TN)$ be the jet bundle over $M \times N$ 
(with fiber $J(m,n)$ of order high enough $(\ge n+1)$ and group $\A$), 
and denote by  $jf$ the jet extension which assigns to points $x \in M$ the pair of 
$f(x) \in N$ and the jet of the germ $f: M, x \to N, f(x)$: 
$$
\xymatrix{
&J(TM, TN) \ar[d]\\
M \ar[ur]^{jf} \ar[r]_{(id, f)} & M \times N 
}
$$
\begin{center}
\begin{tabular}{l}
$f: M,x \to N,  f(x)$ is stable \\
$\Longleftrightarrow$ $jf:M \to J(TM, TN)$ is transverse  to the $\A$-orbit at $x$. \\
$\Longleftrightarrow$ 
$jf:M \to J(TM, TN)$ is transverse  to the $\K$-orbit at $x$.
\end{tabular}
\end{center}

\

\noindent
{\bf  Versal unfolding.}  
An unfolding of $f: \C^m,0 \to \C^n, 0$ is a map-germ 
$$F: \C^m\times \C^k, (0,0) \to \C^n\times \C^k, (0,0), \; \; 
F(x, u)=(f_{u}(x), u)$$
 so that $F(x,0)=(f(x),0)$ (i.e., $f_0=f$).  
 Note that $f$ itself is regarded as an unfolding without parameters ($k=0$). 
 Two unfoldings $G, F$ of $f$ with $k$ parameters are {\it equivalent} if 
there are unfoldings of identity maps $id_m$ of $\C^m$ and $id_n$ of $\C^n$, 
$$\Phi: \C^m \times \C^k, 0 \to \C^m \times \C^k, 0, \quad 
\Psi: \C^n \times \C^k, 0 \to \C^n \times \C^k, 0, $$
respectively, so that  $F\circ \Phi=\Psi \circ G$. 
An unfolding of $f$ is {\it trivial} if it is equivalent to  
the product $(f\times id_k)(x, u):= (f(x), u)$. 

Given a map $h: \C^\ell, 0 \to \C^k,0$, 
the {\it induced unfolding $h^*F$ from $F$ via the base-change $h$} is defined by 
the unfolding $h^*F(x,v) := (f_{h(v)}(x), v)$. 

We say that $F$ is an {\it $\A _e$-versal unfolding of $f$} 
if any unfolding of $f$ is equivalent to some unfolding induced from $F$. 
The so-called {\it versality theorem} says that 
$F$ is $\A _e$-versal if and only if it holds that 
$$\textstyle \theta(f)=T\A_e.f + \sum_{i=1}^k \C \cdot \frac{\rd}{\rd u_i} f_u \bigr{|}_{u=0}.$$
This identity means that 
the map 
$$U\times W \to V \times \E_0(m,n), \; (p,u) \mapsto \mbox{germ of $f_u(x+p)$ at $x=0$}$$
is transverse to $V \times \A.f$ at $(p,u)=(0,0)$, 
where $F: U \times W \to V \times W$ is a representative. 

The {\it $\A_e$-codimension} of $f$ is defined to be 
$\dim_\C \theta(f)/T\A_e.f$, 
which is the minimum number of parameters required 
for constructing an $\A_e$-versal unfolding of $f$. 
In particular, 
\begin{center}
\begin{tabular}{c}
$f$ is a stable germ 
$\Longleftrightarrow$ $\A_e\mbox{-}{\rm codim}(f)=0$\\
 $\Longleftrightarrow$ 
$f$ itself is $\A_e$-versal 
$\Longleftrightarrow$  any unfolding of $f$ is  trivial. 
\end{tabular}
\end{center}

\

\subsection{Basics in intersection theory} \label{Chern} 
Basic references are,  e.g., 
 \cite{MS, Hirzebruch, Fulton, Suwa, Kleiman}. 

\

\noindent
{\bf  Homology.} 
Throughout, $H^*$ and $H_*$ stand for 
the singular cohomology ring (with cup product) 
and the {\it Borel-Moore} homology group 
(=the closed supported homology 
=the homology of locally finite chains),  respectively.  
\begin{itemize}
\item 
$H^*$ is a contravariant functor: 
the pullback $f^*: H^*(Y) \to H^*(X)$ is a ring homomorphism 
defined for a continuous map $f: X \to Y$, and it holds that $(g\circ f)^*=f^*\circ g^*$. 
\item
$H_*$ is covariant for {\it proper} maps: 
the pushforward $f_*: H_*(X) \to H_*(Y)$ is a group homomorphism 
defined for a proper continuous map $f: X \to Y$, 
and it holds that $(g\circ f)_*=g_*\circ f_*$. 
For compact spaces  it is the same as the usual homology group. 
\end{itemize}
There is a natural pairing (cap product): $\frown: H^k(X) \times H_m(X) \to H_{m-k}(X)$. 
The two maps $f^*$ and $f_*$ are related by the useful {\it projection formula}:
$$f_*(f^*\alpha \frown c)=\alpha \frown f_*c $$
for  $\alpha \in H^*(Y)$ and $c\in H_*(X)$. 
For a (possibly non-compact) complex irreducible variety $X$ of dimension $m$, 
there always exists the {\it fundamental class} 
$[X] \in H_{2m}(X)$: For any regular point $x \in M$, 
the class generates $H_{2m}(X, X-x) \simeq \Z$ 
being compatible with the complex orientation. 
If $M$ is a complex manifold, 
it yields the well-known {\it Poincar\'e duality} isomorphism 
$$H^i(M) \simeq H_{2m-i}(M), \;\; \omega \mapsto \omega \frown [M].$$
We denote by $\Dual c$ the Poincar\'e dual to $c \in H_*(M)$ 
but often omit this notation when it would not cause any confusion. 

For proper maps $f: M \to N$ between manifolds of 
relative codimension $\kappa=\dim N-\dim M$, 
the Gysin homomorphism is 
defined by 
the dual to the homology pushforward (we abuse the notation): 
$$f_*=\Dual \circ f_* \circ \Dual^{-1}: H^*(M) \to H^{*+\kappa}(N).$$ 
For instance, $f_*(1)=\Dual\, f_*[M]$. 

\begin{prop} \label{transverse}
If $f: M \to N$ between complex manifolds is transverse to 
a closed subvariety $Y \subset N$, then the pullback of $[Y]$ is expressed by 
the preimage of $Y$ via $f$, 
$f^*\Dual\, [Y] = \Dual [f^{-1}(Y)]  \in H^*(M)$. 
\end{prop}

\noindent
{\bf  Chow group.} 
In the  context of algebraic geometry, instead of $H_*$, 
we may take the Chow group $A_*$ of 
algebraic cycles under rational equivalence \cite{Fulton}: 
The group of algebraic $k$-cycles on a variety $M$ is 
freely generated by symbols $[V]$ associated to 
$k$-dimensional closed irreducible subvarieties $V \subset M$, and 
two algebraic $k$-cycles are {\it rationally equivalent} 
if they are joined by a family of cycles parametrized by $\Proj^1$ 
(such a family forms an algebraic $(k+1)$-cycle on $M \times \Proj^1$). 
The pushforward $f_*: A_*(M) \to A_*(N)$ 
 is defined for proper algebraic morphisms $f: M \to N$ 
 by $f_*[V]=\deg(f|_V)\cdot [f(V)]$ if $\dim V=\dim f(V)$, and $0$ otherwise. 
If $M$ is non-singular and of dimension $m$,  
we put 
$$A^*(M)=\oplus\,  A^k(M), \quad A^k(M):=A_{m-k}(M).$$
The intersection product of two algebraic cycles 
is generally defined 
(\cite[\S 20]{Fulton}, \cite{Kleiman}), 
that put on $A^*(M)$ a ring structure; then it is called the Chow ring of $M$. 
The pullback $f^*: A^*(N) \to A^*(M)$ for a morphism between algebraic manifolds 
is defined by taking a scheme theoretic preimage, 
i.e., the intersection product of the graph of $f$ and the cartesian product  $M$ 
times subvarieties of $N$. 
Over the ground field $\C$, 
there is a ring homomorphism, called the {\it cycle map},  
$$cl: A^*(M) \to H^{2*}(M) $$
sending an algebraic cycle to the dual to the fundamental class of 
the underlying analytic space: 
$cl$ is compatible with 
the pullback and the Gysin homomorphism (pushforward). 
In particular, 
$$cl([V]\cdot [W])=cl([V]) \cdot cl([W]),$$
hence, the algebraic intersection number of cycles (in $A^*$) 
coincides with the topological intersection number 
defined by the cup product (in $H^*$). 

\

\noindent
{\bf Chern classes.} 
A  {\it complex vector bundle} $p: E \to M$ of rank $n$ is a locally trivial fibration 
with fiber $\C^n$ and structure group $GL_n(\C)$: 
$E$ is called the {\it total space}, $M$ the {\it base space} 
and $\C^n$ the {\it fiber}, and 
the {\it zero section} $Z \subset E$ is 
the subvariety consisting of all zero vectors of fibers. 
The pullback induced by the projection map $p$ provides a canonical isomorphism 
$$p^*: H^*(M) \stackrel{\sim}{\longrightarrow} H^*(E).$$

The {\it trivial} $n$-bundle $\epsilon^n$ means that it is globally trivialized, i.e., 
 isomorphic to the product $M \times \C^n \to M$. 
To measure `non-triviality' of a given vector bundle  $p: E \to M$, 
the most basic invariant is the {\it top Chern class} of $E$ 
defined by the fundamental class of the zero section: 
$$c_n(E):=(p^*)^{-1}\Dual [Z]  \;\; \in H^{2n}(M; \Z).$$
For a section $s: M \to E$ (i.e., $p\circ s=id_M$), 
we have $s^*=(p^*)^{-1}$, and 
if $s: M \to E$ is transverse to $Z$, 
then by Proposition \ref{transverse} 
the top Chern class is represented by the {\it degeneracy locus} (zero locus) of $s$: 
$$c_n(E)=s^*\Dual [Z] =\Dual [s^{-1}(Z)].$$
The top Chern class is regarded as a cohomological obstruction for the existence of 
a trivial line sub-bundle of $E$: That means that 
if such a trivial sub-bundle exists, then there is a section $s$ nowhere zero, 
i.e., $Z(s)=\emptyset$,  thus $c_n(E)=0$. 
In the same manner  the lower Chern class $c_{i}(E)$ is introduced 
as a certain obstruction for the existence of a trivial  sub-bundle of rank $n-i+1$.  
So for the trivial bundle $\epsilon^n$,  all Chern classes $c_i(\epsilon^n)$ vanish. 

The Chern classes are also formulated in the following intrinsic way: 
Let $\pi: \Proj(E) \to M$ be the projectivized bundle of lines in $E$, 
then there is an exact sequence 
$$0 \longrightarrow L_E \longrightarrow \pi^*E \longrightarrow Q_E \longrightarrow 0$$
where $L_E$  is the tautological line bundle over  $\Proj(E)$; 
let $\Ost_E(1):=L_E^*$ denote the bundle dual to $L_E$ 
and put $t=c_1(\Ost_E(1))$ (top Chern class). 
Then $H^*(\Proj(E))$ naturally has a $H^*(M)$-module structure via $\pi^*$ 
generated by $t$: 
In fact one can define the Chern class $c_i(E) \in H^{2i}(M)$ by the identity
$$t^n+\pi^*c_1(E)t^{n-1}+\cdots + \pi^*c_n(E)=0 \;\; \in H^{2n}(\Proj(E))$$
which actually generates  the relation $I$ of $H^*(\Proj(E))=H^*(M)[t]/I$.  
In particular, in case that $M=pt$, 
this implies that 
$H^*(\Proj^{n-1})=\Z[t]/(t^n)$. 

\begin{exam} 
{\rm 
 ({\bf Poincar\'e-Hopf}) 
The {\it Chern class of a complex manifold $M$} means $c(TM)$ 
of the tangent bundle. 
If $M$ is compact, 
the top Chern class corresponds to the Euler characteristic of $M$
 $$c_n(TM) \frown [M]=\chi(M) \cdot [pt] \;\; \in H_0(M), $$
that is {\it the Poincar\'e-Hopf theorem} 
for a vector field $v$ on $M$ (a section of $TM$)  
 $$c_n(TM)=\Dual [v^{-1}(Z)]=\sum Ind(v, p) \stackrel{\mbox{\tiny P.H.}}{=}\chi(M).$$
 }
\end{exam}

\noindent
{\bf Axiom.} 
Chern classes satisfy the following axiom  
which is quite useful for actual computation: 
\begin{itemize}
\item $c_0(E)=1$ and $c_i(E)=0\;\; (i>n={\rm rank}\, E)$, i.e., 
$$c(E):=\sum_{i\ge 0} c_i(E)=1+c_1(E)+\cdots + c_n(E)$$ 
which called the {\it total Chern class} of $E$.
\item   ({\it naturality}) For pullback via $f: M' \to M$, 
$$c(f^*E)=f^*c(E).$$
 \item  ({\it Whitney  formula})
For any short exact sequence of vector bundles 
$0 \to E' \to E \to E'' \to 0$, it holds that 
$c(E) = c(E') \cdot c(E'')$, i.e., 
$$ c_k(E) = \sum_{i+j=k} c_i(E')c_j(E'').$$
 \item 
 ({\it normalization}) 
 $c_1(\mathcal{O}_{\mathbb{P}^1}(1))$ equals 
the divisor class $a \in H^2(\mathbb{P}^1)$. 
\end{itemize}

\

For instance, it follows that 
\begin{itemize}
\item[-] 
{\it Trivial  bundle}: For the trivial $n$-bundle, $c(\epsilon^n)=c(\oplus \epsilon^1)=1$. 
\item[-] 
{\it Additive group law}: For tensor product of line bundles $\ell_1, \ell_2$ over $M$: 
$$c(\ell_1 \otimes \ell_2)=1+c_1(\ell_1) + c_1(\ell_2).$$
\end{itemize}

If $E$ splits into line bundles, $E=\ell_1\oplus \cdots \oplus \ell_n$, 
$$c(E)=1+c_1(E)+\cdots + c_n(E)=\prod (1+a_i), $$
where $a_i=c_1(\ell_i)$ called the {\it Chern roots} of $E$. 
So the $i$-th Chern class $c_i(E)$ is nothing but 
the {\it $i$-th elementary symmetric function in $a_1, \cdots, a_n$}. 
This computation is formally allowed for any {\it non-split} vector bundle $E$ 
by regarding it  {\it virtually}  as the sum of line bundles, 
that is  the {\it splitting principle}. 
For instance, the product $E \otimes F$ is virtually regarded as 
the sum of products $\ell_i \otimes \ell_j'$ of  line bundles, 
hence by additive group law, 
$$c(E \otimes F)=\prod c(\ell_i\otimes \ell'_j)=\prod (1+a_i+b_j),$$
where $a_i$ and $b_j$ are Chern roots of $E$ and $F$, respectively. 
The calculus on Chern classes 
is essentially the same as the combinatorics of elementary symmetric functions. 

\

\noindent
{\bf Quotient Chern class.} 
 To measure in a formal way the difference between 
 two vector bundles $E$ and $F$ of rank $m, n$ over the same base space, 
 we define the quotient Chern class 
$$
 c(F-E)=\sum_{i\ge 0} c_i(F-E) := \frac{1+c_1(F)+ \cdots }{1+c_1(E)+\cdots} 
 =\frac{\prod (1+b_j)}{\prod (1+ a_i)} 
$$
by using formal expansion $\frac{1}{1-a}=1+a+a^2+\cdots$. 
Obviously, 
if $F=E\oplus E'$, then $c(F-E)=c(E')$. 

Let $P$ be a polynomial in components $c_i(E)$ and $c_j(F)$ ($i, j=1, 2, \cdots)$ 
i.e., 
$P=P(a_1,\cdots, a_m, b_1, \cdots, b_n)$ is symmetric 
in both variables $a_i$ and $b_j$. 
It is known that $P$ is written as 
a polynomial in quotient Chern classes $c_i(E-F)$ if and only if 
$P$ is {\it supersymmetric}, that is, 
$$P(a_1,\cdots, a_{m-1}, t, b_1, \cdots, b_{n-1}, t)$$
does not depend on $t$ (A. Lascoux).

The $K$-group $K_0(M)$ is the group completion of 
the monoid generated by isomorphism classes of vector bundles with 
the Whitney sum operation $\oplus$. 
Then the Chern class operation $c_*: K_0(M) \to H^*(M)$ is well-defined. 
Moreover, the {\it Chern character} of $E$ is defined by 
$ch(E)=\sum \exp{a_i}$ using Chern roots $a_i$ of $E$, 
and it produces a natural transformation $ch: K_0(M) \to H^*(M)$ as ring homomorphism 
(where $K_0(M)$ is a commutative ring with the tensor product $\otimes$).

 \begin{exam}
{\rm 
 ({\bf B\'ezout's theorem}) 
Let $\ell = \mathcal{O}_{\mathbb{P}^2}(1)$ be the dual tautological line bundle of 
the projective plane $\Proj^2$. 
Put $a=c_1(\ell) \in H^2(\Proj^2)$, the dual to a line. 
A homogeneous polynomial $P(x,y,z)$ of degree $d$ assigns 
to each point $[L] \in \Proj^2$ the function 
$L \to \C$ given by  $t\bv \mapsto P(\bv)t^d$, 
which gives a section of the line bundle tensorred $d$ times 
$\mathcal{O}_{\mathbb{P}^2}(d):=\ell^{\otimes d}$. 
The zero locus of this section is nothing but 
the projective plane curve defined by $P=0$. 
Since $c(\ell^{\otimes d})=1+ d\cdot a$, 
the fundamental class of the plane curve $P=0$ is represented by 
the top Chern class $d\cdot a$.  
For two projective plane curves of degree $d$ and $d'$ without common factor, 
the sum of algebraic intersection numbers corresponds to 
the cup product of their fundamental classes, 
$c_1(\ell^{\otimes d})\cdot c_1(\ell^{\otimes d'})=dd'\cdot a^2 
\in H^4(\Proj^2)=\Z$ via the cycle map $cl$. 
This means  {\it classical B\'ezout's theorem}. 
}
\end{exam}

\section{Chern class for singular varieties} 

\subsection{Singular Chern classes}  
As seen in the previous section, 
the Chern class of an $n$-dimensional complex manifold $X$ 
is the total cohomology class 
$$c(TX)=1+c_1(TX)+\cdots + c_n(TX) \; \in H^*(X).$$
Note that $c_n(TX) \frown [X] = \chi(X)$ and $1\frown [X]=[X]$. 
For a singular variety $X$, there is no longer the tangent bundle, 
so $c(TX)$ does not make sense at all. 
However we may have a chance to find some substitute to $TX$, 
for instance 
by taking a reasonable partial desingularization $p: \widehat{X} \to X$ 
(e.g. Nash blowing-up, which will be described below) 
or a deformation to smooth varieties $X_t$ if it exists. 
Then we consider Chern classes of the substitute on $\widehat{X}$ or $X_t$. 
According to the direction of `arrow' $p$, 
we switch to homology and take the image of the Chern class via 
the pushforward  $p_*: H_*(\widehat{X}) \to H_*(X)$ 
or the specialization map $sp_*: H_*(X_t) \to H_*(X)$, 
that provide a kind of ``singular Chern classes" defined in $H_*(X)$. 
The {\it Chern-Schwartz-MacPherson class} (CSM class) 
is a typical one: It is the most useful `singular Chern class' from the functorial viewpoint, 
which we briefly introduce in this section. In particular, 
the CSM class of a (compact, irreducible) possibly singular variety $X$ 
is a total homology class of the form 
$$c^{\SM}(X)=\chi(X) \cdot [pt] + \cdots + [X] \; \in H_*(X).$$

Throughout this section, (Borel-Moore) homology $H_*$ can be replaced by Chow group $A_*$.

\subsection{Chern-Schwartz-MacPherson class}  \label{CSM}
Let $X$ be a complex algebraic variety of dimension $n$. 
For a subvariety $W \subset X$,  
we denote by $\jeden_W: X \to \Z$ 
the characteristic function which takes value $1$ on  points of $W$, otherwise $0$. 
Then a {\it constructible function} $\alpha: X \to \Z$ is 
a function on $X$ 
given by a finite sum $\alpha=\sum n_i \jeden_{W_i}$ with $n_i \in \Z$, 
$W_i$ subvarieties of $X$. 
Let $\F(X)$ be the abelian group of constructible functions on $X$. 
The {\it integral of $\alpha$} is defined by 
$$\int_X \alpha  := \sum \, a_i\, \chi(W_i),$$
where $\chi$ means the Euler characteristics 
using the Borel-Moore homology of underlying analytic spaces.  
Furthermore, for morphisms $X\to Y$, the {\it pushforward} is defined by  
$$f_*: \F(X)\to \F(Y), \quad f_*(\alpha)(y):=\int_{f^{-1}(y)} \alpha \quad 
(y \in Y).$$
Note that $\int_X\alpha=pt_*\alpha \in \Z=\F(pt)$ where $pt: X \to pt$. 
It holds that 
$$(f\circ g)_*=f_*\circ g_*.$$
Also the pullback $f^*: \F(Y) \to \F(X)$ is defined by $f^*\alpha:=\alpha \circ f$. 

The group of constructible functions  
$\F$ and the Borel-Moore homology $H_*$ define covariant functors  
$Var \to Ab$ 
from the category of complex algebraic varieties and proper morphisms 
to the category of abelian groups. 
\begin{thm}\label{C_*} \cite{Mac}
There is a unique natural transformation 
$$C_* : \F(X) \longrightarrow H_*(X)$$ 
between these functors 
so that $C_*(\jeden_X)=c(TX)\frown [X]$  if $X$ is non-singular. 
\end{thm}
Naturality means that 
\begin{itemize}
\item $C_*(\alpha+\beta)=C_*(\alpha)+C_*(\beta)$ (additive homomorphism) 
\item 
$C_*f_*(\alpha)=f_*C_*(\alpha)$ for proper morphisms $f: X \to Y$. 
\end{itemize}
In particular,  
if $pt: X\to pt$ is proper, we have 
$$pt_*C_*(\alpha)=C_*pt_*(\alpha)=\int_X \alpha$$ 
(where $C_*: \F(pt)=H_0(pt)$), 
hence for $\alpha=\jeden_X$, 
the $0$-dimensional component of $C_*(\jeden_X)$ 
corresponds to $\chi(X)$. 
For irreducible $X$,  
the top dimensional component of $C_*(\jeden_X)$ 
is the fundamental class $[X]$, as see below. 

\begin{definition}
{\rm 
The {\it Chern-Schwartz-MacPherson class} of $X$ 
is defined by $c^{\SM}(X):=C_*(\jeden_X)$. 
For non-reduced scheme $X$, we define 
$c^{\SM}(X):=c^{\SM}(X_{Red})$ of the underlying reduced scheme. 
}
\end{definition}

In the later sections, we  consider the CSM class of $X$ 
in  ambient smooth space $M$; 
We  often write $C_*(\alpha) \in H^*(M)$ for $\alpha \in \F(M)$ 
without the notation $\Dual$, 
that would not cause any confusion. 

\begin{rem}{\rm 
({\bf Schwartz class}) 
Historically earlier than MacPherson's paper \cite{Mac},   
M. Schwartz had defined  an obstruction class in the relative cohomology 
$H^*(M, M-X)$ 
 for the existence of {\it radial} vector frames over 
 a neighborhood of $X$ in an ambient manifold $M$, 
that can be seen as a special kind of {\it degeneracy loci class} 
(for frames controlled 
 in a tubular neighborhood of each stratum of a fixed Whitney stratification of $X$).  
The Schwartz class coincides with $C_*(\jeden_X)$ via the Alexander duality 
$H^*(M, M-X) \simeq H_{2m-*}(X)$, that was proved in Brasselet-Schwartz \cite{BS}. 
}
\end{rem}

\begin{rem} {\rm 
({\bf Nash blow-up and Chern-Mather class}) 
We briefly explain about MacPherson's original construction of $C_*$ 
in \cite{Mac} using a specified desingularization - the Nash blow-up. 
Assume that $X$ is embedded in an ambient manifold $M$ and 
is of equidimension $n$. 
Let $\nu_M: Gr(TM,n) \to M$ be  
the Grassmannian bundle of $n$-dimensional linear subspaces in $TM$. 
Then, there is a unique section $\rho$ over the regular locus $X_{Reg}$ of $X$ 
which sends $x \in X_{Reg}$ to the tangent space $T_xX$. 
We denote by $\widehat X$ the closure of the image $\rho(X_{Reg})$ and 
by $\nu_X:\widehat X \to X$ the natural projection, 
which is called the  {\it Nash blow-up of $X$}. 
The Nash tangent bundle $\widehat {TX}$ 
is defined by the restriction of the tautological $n$-bundle of the Grassmaniann to $\widehat{X}$. 
Then we define  the {\it Chern-Mather class}  
$$ c^\Ma(X) := \nu_*(c(\widehat {TX})\frown [\widehat X]) \quad \in H_*(X),$$
which is known to be independent of the choice of the embedding. 
This is the main ingredient for defining $C_*$.  
The second one is 
the {\it local Euler obstruction function} $Eu_W \in \F(X)$, 
 which is originally defined using the  obstruction theory for radial vector fields. 
It satisfies that $Eu_W(x)=1$ for nonsingular points $x \in W_{reg}$, and that 
  $\F(X)$ is freely generated by $Eu_W$ of subvarieties $W$ of $X$. 
Then $C_* : \F(X) \to H_*(X)$ is defined by 
$C_*(Eu_W) := \iota_*c^\Ma(W)$, $\iota:W \to X$ being the inclusion, 
and by extending it linearly. 
In fact, $c^\Ma$ and $Eu$ behave in a similar way for pushforward, 
that imply the functoriality of $C_*$. 
If $X$ is smooth, then $c^\SM(X)=c^\Ma(X)=c(TX)\frown [X]$. 
}
\end{rem}

\begin{rem}{\rm ({\bf Motivic type description of $C_*$}) 
There is an alternative convenient description of $C_*$ 
using Hironaka's resolution of singularities. 
Notice that 
any constructible function $\alpha \in \F(X)$ admits a finite sum expression 
$$\alpha=\sum  a_i \, \rho_{i*}\jeden_{M_i},$$
where $a_i \in \Z$ and 
$\rho_0: M_0 \to X$ is a proper surjective birational morphism 
and $\rho_i: M_i \to X \; (1\le i \le s)$ 
is a proper birational morphism mapped to a subvariety 
of  dimension smaller than $\dim X$, 
and all $M_i\; (0\le i \le s)$ are non-singular. 
That is easily shown by using 
resolution of singularities and the induction of the dimension of 
 supports of constructible functions. 
Then, by properties of $C_*$, we see that  
$$C_*(\alpha)=C_*\left(\sum a_i \, \rho_{i*}\jeden_{M_i}\right) 
=\sum a_i \,  \rho_{i*}(c(TM_i) \frown [M_i]).$$
Now let $M^+(X)$ be the free abelian group generated by 
all equivalence classes of {\it proper} morphisms 
$f: M \to X$ with  {\it non-singular} $M$ 
(morphisms $f_1, f_2$ mapped to $X$ are equivalent if $f_1=f_2\circ \sigma$ 
by some isomorphism of sources), and 
define the additive homomorphisms 
$e$ and $\mathfrak{c}_*$ by 
linear extensions of  
$$e[f:M\to X]:=f_*\jeden_M, $$ 
$$\mathfrak{c}_*[f:M\to X]:=f_*(c(TM)\frown[M]), $$ 
respectively. Note that $e$ is surjective. 
Then,  MacPherson's Chern class transformation is expressed by
 $C_*=\mathfrak{c}_* \circ e^{-1}$: 
$$
\xymatrix{
&M^+(X) \ar[dl]_e \ar[dr]^{\mathfrak{c}_*} &\\
\F(X) \ar[rr]_{C_*} && H_*(X)
}
$$ 
We may replace $M^+(X)$ by 
the relative Grothendieck group $K_0(Var/X)$ of 
varieties over $X$, that enables us to deal with  
motivic integrations and stringy Chern classes \cite{Aluffi}, 
and more generally, the (singular) Hirzebruch classes \cite{BSY}. 
}
\end{rem}

\subsection{Segre-SM classes} \label{SSM}
Let $M$ be a complex algebraic manifold. 
We define the {\it Segre-Schwartz-MacPherson class} 
of a closed embedding $\iota: X \hookrightarrow M$ by 
$$s^{\SM} (X,M):= c(\iota^*TM)^{-1}\frown  c^{\SM}(X) \;\; \in H_*(X).$$
We regard the class $s^{\SM} (X,M)$ in $H^*(M)$ via the pushforward $\iota_*$ and $\Dual$. 
Also we set for $\alpha \in \F(M)$ 
$$s^{\SM} (\alpha,M) := c(TM)^{-1}\cdot  C_*(\alpha)\in H^*(M).$$  
Notice that if $X$ is a closed submanifold of $M$ with the normal bundle 
$\nu = \iota^*TM - TX$,  
then the Segre-SM class is nothing but 
the {\it inverse normal Chern class} for $X \hookrightarrow M$: 
$$s^{\SM} (\jeden_X,M)=\iota_* c(-\nu) \; \; \in H^*(M).$$ 

\begin{rem}\label{sign_convention}
{\rm ({\bf Sign convention}) 
We should remark that we follow the sign convention of the Segre class 
due to Fulton \cite{Fulton}. 
The other convention corresponds to the  {\it dual} version $\iota_*c(-\nu^*)$ 
for smooth embeddings. 
An advantage of our convention is 
that it is easy to switch between $C_*$ and $s^{\SM}$ 
via multiplying by the ambient Chern class $c(TM)$. 
It could be possible to follow the other convention, 
which fits the {\it positivity property} especially, but to do this 
we had to correct the normalization condition of $C_*$ 
so that $C_*(\jeden_X)=c(T^*X) \frown [X]$ 
for non-singular $X$. 
This causes the change of signs of each component  
$C_i$ by multiplying $(-1)^i$.
}
\end{rem}

\begin{rem}\label{segre_covariance}{\rm 
({\bf Fulton's Chern class}) 
The ordinary {\it Segre covariance class} $s(X, M)$ 
of a closed embedding $X \hookrightarrow M$ 
is defined using the blowing-up $M$ along $X$,  
and it is totally different from our Segre-SM class in general: 
The difference concentrates on the singular locus, and in fact 
these two Segre classes coincide if $X$ is non-singular. 
Our definition of Segre-SM class is just an analogy to 
{\it Fulton's Chern class} \cite{Fulton} defined by 
$$c^{\rm F}(X):=c(\iota^*TM)\frown s(X, M).$$
The difference between these two homology Chern classes  
is an important invariant of singularities of $X$,  called the {\it Milnor class}: 
$$\mathcal{M}(X):=(-1)^{\mbox{\tiny ${\dim X}$}}(c^{\rm F}(X)-c^{\SM}(X)).$$
}
\end{rem}

The Segre-SM class has an expected nice property for transverse pullback 
like as the fundamental class in Proposition \ref{transverse}. 

\begin{prop}\label{segre} 
Let  $f:M \to N$ be a map between complex manifolds, and 
let $Y$ be a closed singular subvariety of $N$.  
Assume that $f$ is transverse to (a Whitney stratification of) $Y$. 
Then it holds that 
$$f^*s^{\SM}(Y,N)=s^{\SM}(f^{-1}(Y),M) \;\;\;\in \;\;H^*(M).$$
\end{prop} 

In fact  the formula holds in $H_*(X)$. 
This is a special case of  the Verdier type Riemann-Roch formula, see \cite[Cor. 0.1]{Joerg} 
based on micro-local techniques. 
Here, for the sake of completeness, 
we give an elementary proof.  

\

{
\proof \\
\noindent
(Step 1)  
Assume that $f:M^m \to N^n$ is a closed embedding ($\kappa=n-m\ge 0$). 
Put $E=f^*TN/TM$, the normal bundle of $M$ in $N$, and 
Let $i: Y \hookrightarrow N$ a closed embedding and $p=\dim Y$.
Let $\nu_Y:\widehat{Y} \to Y$ be the Nash blowing-up of $Y$ defined 
in the Grassmaniann bundle $\mu_N: Gr(TN, p) \to N$. 
Let $i': X:=f^{-1}(Y) \hookrightarrow M$, the transverse intersection of $M$ with $Y$  
($\dim X=p-\kappa$), 
and $\nu_X:\widehat{X} \to X$ the Nash blowing-up of $X$. 

Let $\{S^\alpha\}$ be a Whitney stratification of $Y$. 
By the assumption, $\{M \cap S^\alpha\}$ is a Whitney stratification of $X$ 
so that $T(M \cap S^\alpha)_x = TM_x \cap TS^\alpha_x$. 
In particular, if $S^\alpha$ is a top dimensional stratum and 
$S^\beta$ is a nearby stratum, then 
a limiting tangent of  $Y$  at $x \in M \cap S^\beta$, 
$\lambda_x = \lim TS^\alpha_{x_i}$ with $x_i \to x$ ($x_i \in S^\alpha$) 
corresponds in 1-to-1 to 
a limiting tangent of $X$ at $x$,  
$\lambda_x'= \lim T(M\cap S^\alpha)_{x_i} = TM_x \cap \lambda_x$; 
indeed, $TS^\beta_x \subset \lambda_x$  by the $a$-regularity, 
and $TS^\beta_x$ is transverse to $TM_x$ by the assumption, 
hence $\lambda_x$ is so. 
Thus $\widehat{X}$ is canonically identified with 
the restriction of $\widehat{Y}$ over $M\cap Y$, 
so we have the fiber square where 
$f$ and $\bar{f}$ are regular embeddings with normal bundles $E$ and $\nu_X^*E$: 
$$
\xymatrix{
\widehat{X} \ar[r]^{\bar{f}} \ar[d]_{\nu_X} & \widehat{Y} \ar[d]^{\nu_Y}\\
X \ar[r]_f & Y
}
$$

Note that $\bar{f}^*\widehat{TY} = \widehat{TX}\oplus \nu_X^*E$. 
By properties of the {\it refined Gysin pullback} 
in  \cite[Thm.6.2, Prop. 6.3]{Fulton}, we have 
$$f^* (\nu_X)_*=(\nu_Y)_* \bar{f}^* \;\; \mbox{and} \;\; 
\bar{f}^*[\widehat{Y}]=[\widehat{X}], $$
and hence 
\begin{eqnarray*}
f^* c^\Ma(Y) &=& f^*  (\nu_Y)_* (c(\widehat{TY})\frown [\widehat{Y}]) \\
&=& (\nu_X)_* \bar{f}^* (c(\widehat{TY})\frown [\widehat{Y}])\\
&=&(\nu_X)_*(\nu_X^*c(E) \cdot c(\widehat{TX}) \frown [\widehat{X}]) \\
&=& c(E) \frown (\nu_X)_*(c(\widehat{TX})\frown [\widehat{X}])\\
&=& c(E) \frown c^\Ma(X). 
\end{eqnarray*}
Thus we have $f^* s^\Ma(Y,N) =s^\Ma(X,M)$. 

\noindent
(Step 2) General case:  Let $f: M \to N$ be a map transverse to $Y$. 
Put $\Delta:M \to M\times M$  the diagonal map, 
and consider the graph embedding 
$$g= (id_M\times f) \circ \Delta: M \to M\times N.$$
The normal bundle of $g$ is isomorphic to $f^*TN$. 
Let $Y':=M\times Y$, then $X:=f^{-1}(Y)=g^{-1}(Y')$. 
Since $f$ is transverse to $Y$, the embedding $g$ is transverse to $Y'$, 
hence as seen just above, 
$$g^*c^{\Ma}(Y')=c(f^*TN)\frown c^{\Ma}(X).$$
On one hand, since since $C_*$ commutes with homology cross product, 
$$c^{\Ma}(Y')=c^M(M\times Y)=c^{\Ma}(M)\times c^{\Ma}(Y).$$ 
For a manifold,  $c^{\Ma}(M)=c(TM) \frown [M]$, therefore 
\begin{eqnarray*}
g^*c^{\Ma}(Y')&=&\Delta^*\circ (id_M\times f)^*(c^{\Ma}(M)\times c^{\Ma}(Y))  \\
&=& c(TM) \frown f^*c^{\Ma}(Y). 
\end{eqnarray*}
It then follows that $f^*s^{\Ma}(Y,N)=s^{\Ma}(X,M)$. 

\noindent
(Step 3) Write $\jeden_Y=\sum_S n_S Eu_S$ 
for some subvarieties $S$ of $Y$ 
 (including $Y$ itself; $n_Y=1$),  where $S_{reg}$ are strata of 
a Whitney stratification of $Y$. 
Since $f$ is  transverse to each stratum $S_{reg}$, 
we obtain  $\jeden_X =\sum_S n_S Eu_{M\cap S}$ 
by a property of the Euler obstruction for transverse intersections \cite{Mac}. 
Hence, putting $E=f^*TN-TM$, 
\begin{eqnarray*}
f^* c^{\SM}(Y)&=&\sum \,  n_S \, f^* c^{\Ma}(S)\\
&=&\sum n_S \,  c(E) \frown c^\Ma(M\cap S)\\
&=& c(E)\frown C_*(\sum n_S Eu_{M\cap S}) \\
&=& c(E)\frown c^{\SM}(X). 
\end{eqnarray*} 
Thus $f^* s^{\SM}(Y,N)=s^{\SM}(X,M)$. 
This completes the proof. \qed

}

\subsection{Equivariant Chern/Segre-SM class} \label{equivCSM} 
There has been established the theory of equivariant CSM class 
by the author \cite{Ohmoto06, Ohmoto12}, 
which is based on the equivariant intersection theory \cite{EG}. 
In the latter sections, however,  we avoid technical matters in the theory 
as much as possible, so readers may skip most of this subsection,  
and, instead,  take Definition \ref{tp_G^SM}  and Theorem \ref{segre2} below
 as the starting point for reading the following sections. 

To state theorems precisely, 
we briefly explain about the {\it algebraic Borel construction} \cite{EG}. 
The idea is classical and simple: 
Let $G$ be a  complex linear algebraic group of dimension $g$. 
Take a Zariski open subset $U$ in  
an $\ell$-dimensional linear representation of $G$ 
so that $G$ acts on $U$ freely. Then the quotient variety $U_G:=U/G$ exists, 
and the inductive limit of the quotient map $U \to U_G$ taken 
over all representations of $G$ (with respect to inclusions) 
is regarded as an algebro-geometric counterpart of 
the universal principal bundle $EG \to BG$ in topology. 

\begin{exam}{\rm 
For the algebraic torus $T=\C^*=\C-\{0\}$, 
the quotient map $U=\C^N-\{0\} \to \Proj^N=U_T$ with dimension $N$ large enough  
is the substitute to $ET \to BT$. For the general linear group $G=GL_n$, 
let $U$ be an open set in $\Hom(\C^n, \C^N)$ 
consisting of injective linear maps, 
then $U_G$ is the Grassmaniann of $n$-planes in $\C^N$ 
and the quotient map approximates $EGL_n \to BG L_n$. 
}
\end{exam}

Let $X$ be an algebraic variety with a $G$-action. Then 
the diagonal action of $G$ on $X \times U$ is also free, 
hence the mixed quotient $X_U:=(X\times U)/G$ exists so that  
the projection $p_U: X_U \to U_G$ is a fiber bundle with fiber $X$ and group $G$. 
Then the {\it $G$-equivariant cohomology} of $X$ is given as the projective limit 
$$H_G^*(X)=H_G^*(EG\times_G X)=\varprojlim H^*(X_U).$$ 
This becomes a contravariant functor:  
the pullback of a $G$-morphism $f$  is denoted by $f^*_G$.  

Let $\xi$ be a $G$-equivariant vector bundle $E \to X$ 
(i.e., $E$, $X$ are $G$-varieties and the projection is $G$-equivariant 
so that the action on $E$ preserves fibers linearly). 
Then we have a vector bundle $E_U \to X_U$ over the mixed quotient for each $U$, 
denoted by $\xi_{U}$, and define 
the {\it $G$-equivariant Chern class} 
$c^G(\xi) \in H_G^*(X)$ 
to be the projective limit of Chern classes $c(\xi_{U})$. 

We define {\it the $i$-th equivariant homology group}  to be 
the inductive limit 
$$
H^G_i(X) = \varinjlim  H_{i+2(\ell-g)}(X_U). 
$$
(in fact, the right hand side is stabilized for large $\ell$). 
This group is trivial for $i>2n$, but unlikely the non-equivariant case, 
it is nontrivial for $i<0$ in general. 
The direct sum is denoted by $H_*^G(X)=\oplus_{i \in \Z}  H^G_i(X)$. 
For a proper $G$-morphism $f:X \to Y$, 
the pushforwad $f^G_*$ is defined by taking limit 
of $(f_U)_*:X_U\to Y_U$; thus 
$H^G_*$ becomes a covariant functor. 

The (Borel-Moore) fundamental class  $[X_U]$  tends to 
a unique element of $H_{2n}^G(X)$, denoted by $[X]_G$, which is 
called {\it the $G$-equivariant fundamental  class} of $X$. 
It induces a homomorphism 
$$
\frown [X]_G: H_G^{i}(X) \to H_{2n-i}^G(X), \quad 
a \mapsto r_U(a) \frown [X_U]
$$
where $r_U$ denotes the restriction to $X_U$. 
If $X$ is nonsingular, this is isomorphic, called  
{\it the $G$-equivariant Poincar\'e dual}. 
The inverse is denoted by $\Dual_G$. 

We are now ready to state the equivariant version of Therorem \ref{C_*}. 
Let $\F^G_{inv}$ denote the group of $G$-invariant constructible functions. 
We define 
$$
C_i^G(\jeden_X):=p_U^*c(TU_G)^{-1}\frown C_{i+\ell-g}(\jeden_{X_U}). 
$$

\begin{thm}  \cite{Ohmoto06, Ohmoto12} 
For $G$-varieties and proper $G$-morphisms, 
there is a unique natural transformation 
$$C_*^G: \F^G_{inv}(X) \to H^G_*(X)$$ 
so that $C_*^G(\jeden_X)=c^G(TX) \frown [X]_G$ if $X$ is non-singular. 
\end{thm}

\begin{rem}\label{degree_CSM}
{\rm 
Each dimensional component of the equivariant CSM class 
has its support on an invariant algebraic cycle in $X$ \cite[\S 4.1]{Ohmoto06}. 
In particular, the lowest and highest terms are as follows: 
if $X$ is of equidimension $n$, 
the top term is the fundamental class:  
$C_n^G(\jeden_X) = [X]_G$. 
If $X$ is compact, the degree is equal to the weighted Euler characteristics 
 (the pushforward of $pt: X \to pt$): 
$pt_*^G C_0^G(\alpha)=\int_X \alpha$.  
}
\end{rem}

Next, we introduce the {\it degeneracy loci formula} associated to the CSM class 
which has been formulated in \cite{Ohmoto06}. 
Let $V=\C^n$ on which $G$ acts linearly, and 
identify 
$$H_G^*(V)=H^*(BG)$$ 
via the pullback of the projection $pt: V \to 0$. 
For this purpose, 
the right object is the Segre-SM class rather than the CSM class: 

\begin{definition} \label{tp_G^SM} 
For any invariant function $\alpha \in \F_{inv}^G(V)$, we define 
$$
tp_G^{\SM}(\alpha):=c^G(TV)^{-1} \cdot \Dual_G\, C^G_*(\alpha) \; \in H^*(BG).  
$$
We set $tp_G^{\SM}(W):=tp_G^{\SM}(\jeden_W)$ 
for invariant subvarieties $W$ of $V$. 
\end{definition}

We have the following: 

\begin{thm} 
 \label{segre2} \cite{Ohmoto06} 
Let $V=\C^n$ be a $G$-vector space with the fixed point $0 \in V$.  
Let $W$ be a $G$-invariant affine (irreducible) subvariety 
with the inclusion $\iota: W \to V$, and $\alpha \in \F^G_{inv}(V)$ 
an invariant constructible function.  
Then, 
\begin{enumerate}
\item The leading term of $tp^{\SM}(W)$ is the $G$-fundamental class: 
$$tp_G^{\SM}(W)=\Dual_G\, \iota_*^G[W]_G + h.o.t.$$
\item The $G$-degree of $C_*^G(\alpha)$ expresses the integral of $\alpha$: 
$$\Dual_G C_0^G(\alpha) 
=[c^G(TV)\cdot tp_G^{\SM}(\alpha)]_n=\left(\int_V \alpha\right) \cdot c_n^G(TV). $$ 
\item  For any $G$-morphism $\Psi:V' \to V$ which is transverse to $W$, 
it holds that 
$$tp_G^{\SM}(\Psi^{-1}(W))=\Psi^*\, tp_G^{\SM}(W).$$
\item {\rm (Degeneracy loci formula)}  
Given a vector bundle $\pi: E \to M$ over a complex manifold $M$ 
with fiber $V$ and structure group $G$,  
let  $W(E) \to M$ be the fiber bundle with the fiber $W$ and group $G$. 
Then, for any holomorphic section $s: M \to E$  transverse to $W(E)$, 
it holds that 
$$\Dual\, s^{\SM}(W(s), M)=\rho^* tp_G^{\SM}(W) \in H^*(M)$$
where $W(s):=s^{-1}(W(E))$ and $\rho$ is the classifying map for  $E \to M$. 
\end{enumerate}
\end{thm}

\proof 
(1) is obvious since the top term of 
$C_*^G(\jeden_W) \in H_*^G(V)$ is 
the equivariant fundamental class $\iota_*^G [W]_G$. 
(3) follows from Proposition \ref{segre} and 
(4) is just \cite[Thm. 5.11]{Ohmoto06}. 
As for (2), we take the maximal torus $T$ of $G$: 
Since $H^*_G(pt) \to H^*_T(pt)$ is injective (the splitting lemma), 
we may think of the degree via the $T$-action, instead.  
We embed 
$$V \hookrightarrow \Proj^n=\Proj(V\oplus \C)$$ 
equivariantly with respect to the $T$-action ($T$ acts on the second factor $\C$ trivially) 
and compute in two ways the $T$-degree $pt_*^T C_0^T(\alpha) \in H^*_T(pt)$ 
where $pt: \Proj^n \to pt$ is the natural map. 
As mentioned in Remark \ref{degree_CSM}, 
the degree is equal to $\int_{\Proj^n} \alpha$. 
Since the support of $\alpha$ is in $V$, we have
$$pt_*^T C_0^T(\alpha)=\int_V \alpha.$$
Note that $\{0\}$ is a connected component of 
the $T$-fixed point set $(\Proj^n)^T$,  whose normal bundle is 
$T_0V=V$ with the $T$-action. 
Put $j: \{0\} \to \Proj^n$ the inclusion. 
We then apply the Atiyah-Bott localization formula \cite[\S 3 (3.8)]{AB}; 
it can be seen that 
only the contribution from the fixed point $0$ remains,  
i.e., the contribution from fixed point sets in $\Proj^{n-1}=\Proj^n-V$ becomes zero, 
hence we have 
$$pt_*^T C_0^T(\alpha)=\frac{j^* C_0^T(\alpha)}{c_n^T(TV)}.$$ 
Thus (2) is proved (cf. Weber \cite[\S 6]{Weber}).  \qed


\section{Thom polynomials for singularities of maps}  

\subsection{Main Theorems} \label{thom}
Two germs with the same relative codimension, 
say $f:\C^{m+s},0 \to \C^{n+s},0$ 
and $g:\C^m,0 \to \C^{n},0$,  
is called to be {\it stably $\K$-equivalent} 
if $f$ is $\K$-equivalent to  
the trivial unfolding $g\times id_s$ with $s$ parameters. 

Let $\eta$ be a $\K$-singularity type in $\E_0(m, n)$ ($\kappa=n-m$). 
For a holomorphic map $f: M \to N$ with relative codimension $\kappa$, we set 
$$\eta(f):= \{\; x \in M \; | \; 
\mbox{the germ $f$ at $x$ is stably $\K$-eq. to $\eta$} \; \}.$$
If $f$ is a stable map, then the jet extension $jf$ is transverse to 
the $\K$-orbit and $\eta(f)$ consists of stable singularities of type $\eta$. 
We call the (analytic) closure $\overline{\eta(f)} \subset M$ 
the {\it $\eta$-type singular locus} of $f$. 

We are concerned with the simplest primary obstruction for the existence 
of the $\eta$-type singular point for stable map $f$, e.g., 
\cite{Thom, Porteous, Ronga, Damon, 
FR02, FR, Kaz97, Kaz03, Kaz06,  Ohmoto94, Rimanyi}. 

\begin{thm} \label{tp} 
For a stable singularity type $\eta$ as above, 
there exists a unique polynomial 
$tp(\eta) \in \Z[c_1, c_2, \cdots]$
so that for any stable map $f: M \to N$ of  relative codimension $\kappa$, 
 the singular locus of type $\eta$ is expressed by 
the polynomial evaluated by the quotient Chern class $c_i=c_i(f)=c_i(f^*TN-TM)$: 
$$\Dual [\overline{\eta(f)}] = tp(\eta)(c(f)) \;\;   \in H^{2\, \mbox{\tiny $\codim \eta$}}(M).$$
\end{thm}

\begin{definition}{\rm 
We call $tp(\eta)$ {\it the Thom polynomial of 
stable singularity type $\eta$}.  
}\end{definition}

As an advanced version, 
the theory of {\it Thom polynomials for stable multi-singularities} 
has been developed by M. Kazarian \cite{Kaz03, Kaz06}, 
that merges multiple point formulas 
(developed by Kleiman \cite{Kleiman81, Kleiman90}) 
 and the above Thom polynomials for mono-singularities together 
from the viewpoint of cobordism theory (also see \cite{Rimanyi02, Szucs}). 
That is briefly reviewed in Section \ref{tp_mult}. 

\begin{rem}{\rm 
A major problem is  
to determine the precise form of $tp(\eta)$ 
for a given contact type $\eta$. 
A traditional algebro-geometric method for the computation is to construct  
a suitable embedded resolution of the $\eta$-type singular locus 
$X=\overline{\eta(f)} \subset M$ using flag bundles 
\cite{Porteous, Damon, Ronga} 
or to find 
a suitable projective resolution of the structure sheaf $\Ost_X$ 
(cf. \cite{FW00, FW03}) 
but it usually becomes a very hard task. 
On the other hand, 
a more effective new method, 
called the {\it restriction} or {\it interpolation} {\it method}, 
has been  introduced  by R. Rim\'anyi \cite{Rimanyi}. 
It enables us to compute many $tp$ for stable singularities in nice dimensions,  
see Section \ref{symmetry}. 
Also the Atiyah-Bott type localization formula and the iterated residue formula are 
very useful for computation of $tp$'s, 
about which the reader should be referred to 
B\'erczi-Szenes \cite{BercziSzenes} and 
Feh\'er-Rim\'anyi \cite{FR}. 
As for another interesting questions, 
the positivity of Thom polynomials has firstly been dealt in Pragacz-Weber \cite{PW07}, 
and for applications to Schubert calculus, see e.g.  \cite{FP, FR03, Kaz06}. 
}
\end{rem}

As mentioned in the Introduction, 
it is natural to expect a similar universal expression  
not only for the fundamental class but also for 
some other distinguished cohomology classes supported on the singular locus 
$X=\overline{\eta(f)}$.  
For example, 
if the locus $X$ is a closed submanifold of $M$ 
with the inclusion  $\iota$, 
e.g., $\overline{A_k}$ for Morin maps, 
then 
the Gysin homomorphism image 
$\iota_*c(TX) \in H^*(M)$ of the total Chern class 
would be a reasonable candidate;  
Indeed Ando \cite{Ando} and Levine \cite{Levine} 
partially studied such classes in the case of Morin maps. 
However,  the orbit closure $\overline{\eta}$ is  
singular along some orbits of 
more complicated singularities,  
therefore the $\eta$-type singular locus may be singular. 
So $c(TX)$ does not make sense in general. 

Instead, 
our strategy is to incorporate  
the theory of Chern-Schwartz-MacPherson classes 
into the theory of Thom polynomials.  
There always exists 
$$\iota_*c^{\SM}(X) = C_*(\jeden_X) \in H_*(M) = H^*(M),$$
 and  if $X$ is smooth,  then it equals $\iota_*c(TX)$. 
The right object is rather the {\it Segre-SM class} $s^{\SM}(X, M)$ 
obtained by multiplying $c(TM)^{-1}$ to the CSM class.  
Then, the SSM class admits the following Thom polynomial type expresson: 
\begin{thm}\cite{Ohmoto06, Ohmoto08}. \label{tp^SM}
For $\eta$ as above, 
there is a unique universal power series $tp^{\SM}(\overline{\eta}) \in \Z[[c_1,c_2, \cdots ]]$ 
so that  for any stable map 
$f : M \to N$ of relative codimension $\kappa$ 
it holds that 
$$
\Dual s^{\SM}(\overline{\eta(f)}, M) =  tp^{\SM}(\overline{\eta})(c(f)) \;\; \in H^*(M). 
$$
In particular, if $M$ is compact, the Euler characteristic of the $\eta$-type singular locus 
is given by the degree of  $C_*(\jeden_{\overline{\eta(f)}})$, which 
has a universal expression
$$
\chi(\overline{\eta(f)}) = \int_M \;  c(TM) \cdot   tp^{\SM}(\overline{\eta})(c(f)). 
$$
Furthermore, $tp^{\SM}(\alpha) \in \Z[[c_1,c_2, \cdots ]]$ 
is defined for any $\K$-invariant constructible function $\alpha$ 
in some jet space $J(m, m+\kappa)$ 
so that $tp^{\SM}(\jeden_{\overline{\eta}})=tp^{\SM}(\overline{\eta})$. 
\end{thm}

\begin{definition}{\rm 
We call $tp^{\SM}(\overline{\eta})$ the {\it higher Thom polynomial} for 
the orbit closure $\overline{\eta}$ 
with respect to the Segre-SM class. 
}
\end{definition}

The class $tp^{\SM}(\overline{\eta})$ is actually a power series, but 
do not confuse it with the terminology {\it Thom series} in \cite{FR} 
which is a different notion. 

Since the top term of the homology Chern class $c^{\SM}(X)$ is 
the fundamental class $[X]$, it immediately follows from the above definition that 
switching to the cohomology, 
$$tp^{\SM}(\overline{\eta})=tp(\eta) + \mbox{higher degree terms},$$
i.e., the leading term is just the Thom polynomial. 
The power series $tp^{\SM}(\overline{\eta})$  theoretically exists uniquely, 
 but  it is almost hopeless to find the explicit form of the series in general, 
 because the closure $\overline{\eta}$ contains infinitely many boundary strata   
 of high codimension. 
To compute low degree terms, we use Rim\'anyi's restriction method 
 together with embedded resolution techniques, see \S \ref{symmetry}.

\begin{rem}{\rm 
A prototype of Theorem \ref{tp^SM} can be seen in Parusinski-Pragacz \cite{PP}: 
They actually considered $c^{\SM}(\overline{\Sigma^k})$ 
of the first order Thom-Boardman  strata $\Sigma^k$ 
as a generalization of $tp(\Sigma^k)$, i.e., the degeneracy loci class arising in 
the Thom-Porteous formula \cite{Porteous}. 
In order to make a general statement as above, 
 we appeal to the equivariant theory of  CSM class reviewed in the previous section. 
In particular,  theorems can also be formulated 
appropriately in the context of algebraic geometry 
over an algebraically closed field of  characteristic $0$ 
using Chow groups under rational equivalence. 
}
\end{rem}

\begin{rem}{\rm 
In the same way, higher Thom polynomials 
with respect to the other Segre classes (by using blowing-up, conormal sheaves, etc)
can be defined.  It would be interesting to study 
the difference between these higher Thom polynomials with respect to 
different Segre classes, 
that will be discussed somewhere else. 
}
\end{rem}

\subsection{Proof} \label{proof} 
Essential is Theorem \ref{segre2}. 
Consequently, 
Theorem \ref{tp^SM} for $tp^{\SM}$ is proved  
in entirely the same way as the standard proof  of Theorem \ref{tp} for $tp$. 
Here let us see the common proof of Theorem \ref{tp} and \ref{tp^SM} 
along the argument given in  \cite[\S 7.2]{FR}. 

By finite determinacy, 
we may assume that $\eta \subset J(m,n)$, 
the corresponding $\K$-orbit in a jet space of sufficiently high order. 
Since $\eta$ is also $\A$-invariant, 
there is the sub-bundle of the fiber bundle $J(TM, TN)\to M\times N$ with fiber $\eta$, 
denoted by $\eta(M, N)$. 
For stable maps $f: M^m \to N^n$, by the definition 
$\eta(f)=jf^{-1}(\eta(M, N))$:  
$$
\xymatrix{
& &J(TM, TN) \ar[d] &  \;\; \overline{\eta(M,N)} \ar@{_{(}->}[l]\\
 \overline{\eta}(f) \;\; \ar@{^{(}->}[r] & M \ar[ur]^{jf} \ar[r]_{(id, f)} & M \times N &
}
$$
In particular, by Proposition \ref{transverse} 
 $$\Dual [\overline{\eta(f)}] = jf^*\Dual[\overline{\eta(M,N)}] 
 \in H^*(M).$$
We then apply Section \ref{equivCSM} to this setting: 
$$G:=J\K_{m,n}, \;\; V:=J(m,n), \;\; W:=\overline{\eta}$$

By Theorem \ref{segre2} (3), there is a universal class for  
 the degeneracy loci class $\Dual[\overline{\eta(f)}]$: 
$$tp_G^{\SM}(\overline{\eta}) \in H_G^*(J(m,n)).$$
Note that $J(m,n)$ is contractible and  
$G=J\K_{m,n}$ is homotopic to the $1$-jets $J^1\K_{m,n}=GL_m \times GL_n$. 
Thus 
$$H_G^*(J(m,n))=H_G^*(pt)=H^*(BGL_m)\otimes H^*(BGL_n),$$
that is generated by Chern classes of source and of target: 
In terms of Chern roots 
$a_1, \cdots , a_m$ and $b_1,  \cdots , b_n$  for the source and  target, respectively, 
$$H_G^*(J(m,n))=
\Z[a_1, \cdots, a_m, b_1, \cdots, b_n]^{\mathfrak{S}_m \times \mathfrak{S}_n}.$$
We show that $tp_G^{\SM}(\overline{\eta})(a,b)$ 
is actually written in terms of quotient Chern classes
$$c=1+c_1+c_2+\cdots =
\frac{\prod_{j=1}^n (1+b_j)}{\prod_{i=1}^m (1+ a_i)}.$$
The following key lemma is easily checked: 
\begin{lem} \label{stabilization} \cite{Ohmoto94, FR}. 
The natural embedding of jet spaces 
$$\Psi: J(m,n) \to J(m+s, n+s), \quad \Psi(jg(0)):=j(g\times id_{{s}})(0)$$
is transverse to any $\K$-orbits in $ J(m+s,n+s)$. 
 \end{lem}
Consider the group $G':=G \times GL_s \subset J\K_{m+s, n+s}$ 
which naturally acts on the jet space $J(m+s, n+s)$ and also acts on 
$J(m, n)$ by forgetting the $GL_s$-part 
so that $\Psi$ is $G'$-equivariant. 
Notice that the pullback $\Psi^*$ for $G'$-equivariant cohomology 
is the same as the identity map of $H^*(BG')$. 
Put 
$$\eta_s:=\K_{m+s, n+s}.\Psi(\eta) \subset  J(m+s, n+s), $$ 
then the closure $\overline{\eta_s}$ is also $G'$-invariant, 
$\Psi^{-1}(\overline{\eta_s})=\overline{\eta}$ 
and $\Psi$ is transverse to $\overline{\eta_s}$ by Lemma \ref{stabilization}. 
Hence Theorem \ref{segre2} (2) shows that 
$$ tp_{G'}^{\SM}(\overline{\eta_s})= 
tp_{G'}^{\SM}(\overline{\eta}) \;\; \in H^*(BG').$$
By the definition, the $G'$-SSM class $tp_{G'}^{\SM}(\overline{\eta_s})$ 
is written in Chern roots $a_1, \cdots, a_m, b_1, \cdots , b_n$ and $t_1, \cdots, t_s$ 
but the above formula implies that the SSM class does not depend on $t$-variables, 
in other words, it is supersymmetric, thus is written in quotient Chern classes. 
This completes the proof.  \qed

\subsection{Symmetry of singularities}\label{symmetry}
To compute the precise form of $tp(\eta)$, 
there is an effective method due to R.~ Rim\'anyi \cite{Rimanyi}, 
 called the {\it restriction method}. 
This method is also applicable for computing $tp^{\SM}(\overline{\eta})$ up to 
a certain degree. 
Below we demonstrate how to compute $tp^{\SM}$ for $A_2$ 
in  case of $\kappa=0$: 
$$tp^{\SM}(\overline{A_2}) = \sum_{i\ge 2} \;\; tp^{\SM}_i(\overline{A_2}) \;\; 
\in \Z[[c_1, c_2, \cdots]], 
\qquad \deg tp_i^{\SM}=i.$$ 

\

\noindent 
{\bf Leading term=Tp  (degree two).} 
First, let us consider $tp _2^{\SM}(\overline{A_2})=tp(A_2)$. 
It has the form
$$tp(A_2)=A c_1^2+Bc_2$$
in quotient Chern classes $c_i=c_i(\mbox{target}-\mbox{source})$ 
and our task is to determine the unknown coefficients  $A, B$. 

The key point is a simple fact that 
 weighted homogeneous germs admit a natural torus action $T=\C^*=\C-\{0\}$: 
The normal form of stable type $A_2$  is given by a polynomial map 
$$A_2: \C^2 \to \C^2, \;\; (x,y) \to  (x^3+yx, y),$$
and the torus actions on the source and the target are diagonal: 
$$\rho_0=\alpha\oplus \alpha^{\otimes 2}, \quad 
\rho_1=\alpha^{\otimes 3}\oplus \alpha^{\otimes 2} \qquad (\alpha \in T)$$
so that $A_2 \circ \rho_0=\rho_1 \circ A_2$. 

Take the dual tautological line bundle 
$\ell = \Ost_{\Proj^N}(1)$ over 
a projective space $\Proj^N$ of large dimension $N \gg 0$
(or the classifying space $BT=\Proj^\infty$ of the torus $T$).  
Define 
two vector bundles of rank $2$ 
$$E_0 \;(=E_0(A_2)) :=\ell \oplus \ell^{\otimes 2}, \qquad 
E_1\;(=E_1(A_2)):= \ell^{\otimes 3} \oplus \ell^{\otimes 2}.$$ 
That is, let $\{U_i\}$ be  an open cover of $\Proj^N$ giving a local trivialization of $\ell$ 
with $g_{ij}: U_i \cap U_j  \to T$, 
then the glueing maps $U_i \cap U_j  \to GL_2(\C)$ for $E_0$ and $E_1$ are given by 
$\rho_0\circ g_{ij}$ and $\rho_1\circ g_{ij}$, respectively. 

Since the normal form of $A_2$ is invariant under 
the torus action,  
we can glue together the product maps $id_{U_i} \times A_2: U_i \times \C^2 \to U_i \times \C^2$. 
The resulting map $f_{A_2}:E_0 \to E_1$ is a stable map between the total spaces 
$E_0$ and $E_1$ so that the following diagram commutes and 
the restriction to each fiber 
$$\C^2=(E_0)_x \longrightarrow (E_1)_x=\C^2 \qquad (x \in \Proj^N)$$
 is $\A$-equivalent to the normal form of $A_2$. 
 We call $f_{A_2}$  the {\it universal map for $A_2$}. 
$$\xymatrix{
E_0 \ar[rr]^{f_{A_2}} \ar[dr]_{p_0} & & 
E_1\ar[dl]^{p_1} 
\\
& \Proj^N &
}
$$
The loci $A_2(f_{A_2})$ and  $f(A_2(f_{A_2}))$) are 
just  the zero sections of $E_0$ and  of $E_1$, respectively.

\begin{figure}
\includegraphics[clip, width=10cm]{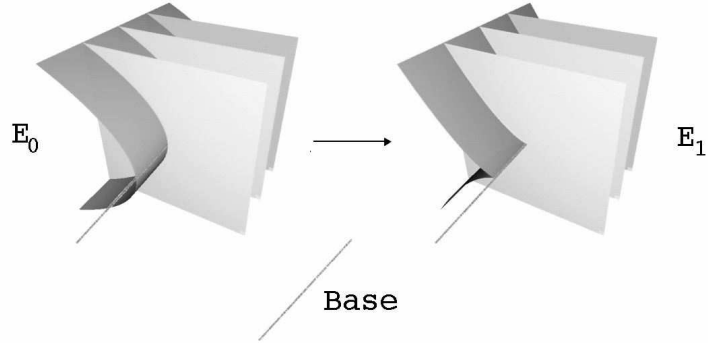}
\caption{\small Universal map for a singularity type}
\end{figure}

Put $a=c_1(\ell)$ and then 
$$H^*(\Proj^N)=\Z [a]/(a^{N+1})\;\;\;\; (N \gg 0), $$
and Chern classes of these vector bundles are written by 
\begin{eqnarray*}
c(E_0)&=&c(\ell \oplus \ell^{\otimes 2})=(1+a)(1+2a), \\
c(E_1)&=&c(\ell^{\otimes 3} \oplus \ell^{\otimes 2})=(1+ 3a)(1+2a).
\end{eqnarray*}
In the following argument, we always identify cohomology rings such as 
$$H^*(E_0) =H^*(\Proj^N) = H^*(E_1)$$ 
through the pullback $p_0^*$ and $p_1^*$. 
For instance, 
since the $A_2$-locus in the total space $E_0$ is the zero section, 
the top Chern class of the pullback bundle $p_0^*E_0$ represents the locus 
in $H^*(E_0)$; So we regard it as  
$$\Dual [\overline{A_2}(f_{A_2})]= c_2(p_0^*E_0)=c_2(E_0) = 2a^2. $$

The tangent bundles $TE_0$ and $TE_1$ of the total spaces canonically 
 split into the vertical and  horizontal components, 
 $$TE_i=p_i^*(E_i\oplus T\Proj^N)\;\;\; (i=0,1),$$
 thus we have in the $K$-group $K_0(E_0)$ 
$$f_{A_2}^*TE_1 - TE_0 = p_0^*(E_1-E_0).$$ 
Therefore, again through the identification $H^*(\Proj^N)=H^*(E_0)$ via $p_0^*$,  
the quotient Chern class for $f_{A_2}$ is written as follows: 
\begin{eqnarray*}
c(f_{A_2})&:=&c(f_{A_2}^*TE_1 - TE_0)\\
&=& c(E_1-E_0)
=\frac{c(E_1)}{c(E_0)}=\frac{1+3a}{1+a} =1+2a-2a^2+\cdots.
\end{eqnarray*}
The first and second degree terms are $c_1(f_{A_2})=2a$, $c_2(f_{A_2})=-2a^2$, 
so we have 
$$tp(A_2)(f_{A_2})=A c_1^2+B c_2 = (4A-2B)a^2.$$
By Theorem \ref{tp} it holds that 
$$tp(A_2)(f_{A_2}) = \Dual [\overline{A_2}(f_{A_2})],$$
hence $(4A-2B)a^2=2a^2$. Thus we have $2A-B=1$. 

\

Next we  apply $tp$ to the universal map 
of adjacent singularities. Let us take the normal form 
$$A_1: \C \to \C, \;\; x \mapsto x^2$$
and the associated universal map $f_{A_1}: E_0 \to E_1$, 
where $E_0=E_0(A_1)=\ell$ and $E_1=E_1(A_1)=\ell^{\otimes 2}$ 
 in the same way as above. 
Obviously, the universal map does not have $A_2$-singularity: 
$A_2(f_{A_1})=\emptyset$, 
thus by Theorem \ref{tp} again, we have 
$$tp(A_2)(f_{A_1}) = \Dual[\emptyset]= 0.$$
Since $c(f_{A_1})=\frac{1+2a}{1+a}=1+a-a^2+\cdots$, 
we have $A-B=0$. 
 
These two linear equations in $A, B$ have a unique solution $A=B=1$, 
thus we conclude that 
$$tp(A_2)=c_1^2+c_2.$$

\

\noindent 
{\bf Degree three term.}
The next term in $tp^{\SM}(\overline{A_2})$ is of degree $3$. 
Put 
$$tp^{\SM}_3(\overline{A_2})= Ac_1^3+Bc_1c_2+C c_3,$$
and determine unknown coefficients. 
We need to restrict this class to more complicated singularities than $A_2$. 

Consider $A_3$-singularity: the stable germ has the normal form 
$$A_3: (x,y,z) \mapsto (x^4+yx^2+zx, y,z).$$
The $T$-action on the source and target spaces are, respectively, 
$$\rho_0=\alpha \oplus \alpha^{\otimes 2} \oplus \alpha^{\otimes 3}, 
\quad \rho_1=\alpha^{\otimes 4}\oplus \alpha^{\otimes 2} \oplus \alpha^{\otimes 3}, $$
which produce the universal map $f_{A_3}: E_0 \to E_1$ over $\Proj^N$ for $A_3$-singularity.  
Then $c(f_{A_3})=\frac{1+4a}{1+a}=1+3a-3a^2+3a^3-\cdots$, 
and hence  
$$tp^{\SM}_3(\overline{A_2})(f_{A_3})=(27A-9B+3C)a^3.$$

The ${A_2}$-locus in the source $\C^3$ is 
a smooth curve tangent to the $x$-axis at $0$ 
and is invariant under the $T$-action,  
thus  $\iota: \overline{A_2}(f_{A_3}) \hookrightarrow E_0$  
is a closed submanifold of codimension $2$.  
The normal bundle is isomorphic to the pullback $\pi^*\nu$ of 
 $\nu=\ell^{\otimes 2} \oplus \ell^{\otimes 3}$ via 
 $\pi=p_0 \circ \iota: \overline{A_2}(f_{A_3}) \to \Proj^N$. 
Since $c(\nu)=(1+2a)(1+3a)$, 
the fundamental class of the locus in  $E_0$ is 
$$\iota_*(1)=c_2(p_0^*\nu) = 6a^2.$$ 

Recall that $tp^{\SM}$ is a universal expression of the Segre-SM class $s^{\SM}$, 
and  for a closed submanifold $X \stackrel{\iota}{\hookrightarrow} M$, 
it is the pushforward of the inverse normal Chern class: 
$$s^{\SM}(X, M)=\iota_*c(-\nu_{M/X}) \in H^*(M).$$
In our case, $X=\overline{A_2}(f_{A_3})$ and $M=E_0$, so 
$$\iota_*(c(-\pi^*\nu))=\iota_*(\iota^*c(-p_0^*\nu))
=c(-p_0^*\nu)\, \iota_*(1)= p_0^*(c(-\nu)c_2(\nu)).$$ 
Thus through the identification via $p_0^*$, 
$$\textstyle  
tp^{\SM}(\overline{A_2})(f_{A_3})
=c_2(\nu)c(-\nu)
=\frac{6a^2}{(1+2a)(1+3a)}=6a^2-30a^3+\cdots.$$
Compare the degree $3$ terms, then we obtain $27A-9B+3C=-30$. 

Again, we restrict $tp^{\SM}$ to adjacent singularities $A_1$ and $A_2$. 
For the universal map $f_{A_2}$, 
$$tp_3^{\SM}(\overline{A_2})(f_{A_2})=(8A-4B+2C)a^3,$$
because we have already seen that $c(f_{A_2})=1+2a-2a^2+2a^3-\cdots$. 
Since the locus $A_2(f_{A_2})$ is the zero section of $E_0=E_0(A_2)$, 
the pushforward of the inverse normal Chern class is 
$$c_2(E_0)c(- E_0)=2a^2-6a^3+\cdots.$$
Comparing the degree $3$ terms, we have $4A-2B+C=-3$. 

 For the universal map $f_{A_1}$, 
$$tp_3^{\SM}(\overline{A_2})(f_{A_1})=0, $$   
 since $A_2(f_{A_1})=\emptyset$. 
 Thus $A-B+C=0$. 
 
These three linear equations have a unique solution: $A=-2$, $B=-3$, $C=-1$, i.e., 
$$tp^{\SM}_3(\overline{A_2})= -(2c_1^3+3c_1c_2+ c_3).$$

\

\noindent 
 {\bf Degree four term.} 
Let us consider the degree $4$ term. 
Using the restriction to $A_k$-singularities ($k=1,2,3,4$) 
we get 
$$tp^{\SM}_4(\overline{A_2})=3c_1^4+6c_1^2c_2+4c_2^2+c_4+A\cdot tp(I_{2,2})$$ 
where $A \in \Z$ is unknown and 
$tp(I_{2,2})=c_2^2-c_1c_3$ for the singularity type 
$$ I_{2,2}:  (x,y, u,v) \mapsto (x^2+2uy, y^2+2vx, u,v).$$
This singularity type is of corank $2$ and the Milnor number is $3$.  
In order to determine $A$, we restrict $tp^{\SM}$ to $ I_{2,2}$. 

The $\overline{A_2}$-locus of the polynomial map $I_{2,2}$ 
is a surface in the source space $\C^4$ having an isolated singular point at $0$ 
(it is defined by $xy-uv=x^2-uy=y^2-vx=0$, so it is not a complete intersection). 
Note that  $\chi(\overline{A_2})=1$. 

Let us consider the $T$-action 
with weights $(1,1,1,1)$ and degrees $(2,2,1,1)$ for the map $I_{2,2}$, 
which produces the universal map $f_{I_{2,2}}: E_0 \to E_1$ 
(where $E_0$ and $E_1$ has rank $4$). 
Then $c(f_{I_{2,2}})=1+2a-a^2+a^4+\cdots$, 
and we substitute them into $tp^{\SM}_4(\overline{A_2})$ described above. 
Since $c(E_0)=(1+a)^4$, 
the CSM class of the $\overline{A_2}$-locus is written by 
$$c(E_0)\cdot tp^{\SM}(\overline{A_2})(f_{I_{2,2}}) = 3a^2+2a^3+(7+A)a^4+\cdots.$$
Now we use Theorem \ref{segre2} (2): 
The degree of the CSM class is  
$$(7+A)a^4=\chi(\overline{A_2})\cdot c_4(E_0) = 1\cdot a^4.$$
Thus $A=-6$, and we have  
$$tp^{\SM}_4(\overline{A_2})=3c_1^4+6c_1^2c_2-2c_2^2-6c_1c_3+c_4.$$ 

In order to seek for higher terms of degree greater than four, 
we need more finer information about the $\overline{A_2}$-locus 
for $I_{2,2}$ and also for more complicated singularity types. 
Here we should combine the restriction method just as described above 
with a traditional method using some $T$-equivariant desingularization 
of the $\overline{A_2}$-locus. 

\

\noindent 
 {\bf Summary.} 
In entirely the same way, 
we compute 
the truncated polynomials of $tp^{\SM}(\overline{\eta})$ up to degree $4$ 
(in case $\kappa=0$): 
{\small 
\begin{eqnarray*}
tp^{\SM}(\overline{A_1})&\equiv& 
c_1-c_1^2+c_1^3-c_1^4+P_1  \\
tp^{\SM}(\overline{A_2}) &\equiv& 
c_1^2+c_2 - (2c_1^3+3c_1c_2+c_3)+3c_1^4+6c_1^2c_2+4c_2^2+c_4+P_2  \\
tp^{\SM}(\overline{A_3}) &\equiv& 
c_1^3+3c_1c_2+2c_3-(3c_1^4+12c_1^2c_2+15c_2^2+6c_4)+P_3  \\
tp^{\SM}(\overline{A_4}) &\equiv& 
c_1^4+6c_1^2c_2+2c_2^2+9c_1c_3+6c_4  \\
tp^{\SM}(\overline{I_{2,2}}) &\equiv& c_2^2-c_1c_3  
\end{eqnarray*}
}
where  
$$P_i=t_i  \cdot tp(I_{2,2}), \qquad t_1=1, \;\; t_2=-6, \;\; t_3=14.$$

As an observation, each term of 
the above $tp^{\SM}(\overline{A_k})$ for Morin maps (i.e. letting $P_i=0$) 
satisfies the positivity both in the Chern monomial basis and  in the Schur polynomial basis
after correcting the sign convention mentioned before, 
i.e., all coefficients are non-negative after multiplying $\pm 1$ accordingly to dimensions. 
But the general form including $P_i$ does not satisfy this property. 

Another observation is concerning the Milnor number constructible function. 
Define $\mu: J(m,m) \to \Z$ by assigning to 
a (jet of) finitely determined germ $f: \C^m,0 \to \C^m, 0$ its Milnor number $\mu(f)$ 
(the value $0$, otherwise). 
This is a constructible function invariant under the $\K$-action and is written by  
\begin{eqnarray*}
\mu&=&
\textstyle 1\jeden_{A_1}+2\jeden_{A_2}+3\jeden_{A_3}+4\jeden_{A_4}+3\jeden_{I_{2,2}}+\alpha \\
&=&
\textstyle
\jeden_{\overline{A_1}}+\jeden_{\overline{A_2}}+\jeden_{\overline{A_3}}+\jeden_{\overline{A_4}} + \alpha'
\end{eqnarray*}
where $\alpha$ and $\alpha'$ are some constructible functions 
having the support of codimension greater than $4$. 
Here, $\jeden_{A_2}$ means the constant function on the $A_2$-orbit and 
$\jeden_{\overline{A_2}}=\jeden_{A_2}+\jeden_{A_3}+\cdots$ 
is the constant function on the orbit-closure. 
Then, summing up $tp^{\SM}(\overline{\eta})$,  we observe a cancellation of several terms 
at least up to degree four: 
\begin{eqnarray*}
tp^{\SM}(\mu)&=&tp^{\SM}(\overline{A_1}) + \cdots + tp^{\SM}(\overline{A_4}) + tp^{\SM}(\alpha') \\
&=&c_1+c_2+c_3+c_4+\cdots.
\end{eqnarray*}
In fact, this is a consequence of a more general property of $tp^{\SM}$ 
for the Milnor number of isolated complete intersection germs 
(=$\K$-finite germs in $\kappa \le 0$), 
which will be discussed in detail somewhere else.

\subsection{Thom polynomials in $\A$-classification} \label{tp^A}

As seen above, 
the Thom polynomial $tp$ for $\K$-classification of map-germs 
is a polynomial in quotient Chern classes 
$c_i(\mbox{source}-\mbox{target})$. 
On the other hand,  Lemma \ref{stabilization} does not hold for $\A$-orbits, 
thus, $tp$ for $\A$-classification is just a polynomial in Chern classes 
of source and that of target. 

A relevant geometric setting for $\A$-classification is described as follows. 
Consider the commutative diagram 

$$\xymatrix{ 
X \ar[rr]^{f} \ar[dr]_{p_0} & & 
Y\ar[dl]^{p_1} 
\\
& B &
}
$$
where $X, Y, B$ are complex manifolds, 
$p_0: X \to B$ and $p_1: Y \to B$ are submersions 
of constant relative dimension $m$ and $n$, respectively.  
For each $x \in X$, 
the germ at $x$ of $f$ restricted to the fiber  is defined:
$$f|_{p_0^{-1}(p_0(x))}: \C^m, 0 \to\C^n,0$$ 
(local coordinates centered at $x$ and $f(x)$). 
Given an $\A$-finite singularity type $\eta$ of maps $\C^m \to \C^n$, 
the {\it singularity locus} $\eta(f) \subset X$ 
and the {\it bifurcation locus} $B_\eta(f)=p_0(\eta(f)) \subset B$ 
are defined. 
It is not difficult to show the following theorem \cite{SSO14}:

\begin{thm} 
Let $\eta$ be an $\A$-finite singularity type. 
For  generic maps $f: X \to Y$, 
$\Dual [\overline{\eta}(f)] \in H^*(X)$ is expressed 
by  a universal polynomial $tp^{\A}(\eta)$ 
in the Chern class $c_i=c_i(T_{X/B})$ and $c_j=c_j(T_{Y/B})$
of relative tangent bundles. 
$\Dual [\overline{B_\eta}(f)] \in H^*(B)$ is also 
expressed by the pushforward $p_{0*}tp^{\A}(\eta)$. 
\end{thm}

$$\xymatrix{
\overline{\eta}(f) \ar[d]_{p_0} \;\ar@{^{(}->}[r] \; & X \ar[rr]^{jf\qquad \quad} 
\ar[d]_{p_0}  &&  J(T_{X/B}, f^*T_{Y/B})
\\
\overline{B_\eta}(f) \; \ar@{^{(}->}[r] \;& B &&
}
$$
\begin{rem}{\rm 
The case of maps between families of curves (e.g., families of rational functions) 
has extensively been studied by Kazarian-Lando \cite{KL04, KL07} 
for the study of Hurwitz numbers. 
}
\end{rem}

\begin{exam}
{\rm 
($\A$-classification of $\C^2, 0 \to \C^2,0$) \\
 Let us see Table 1: the list of $\A$-simple germs of plane-to-plane maps 
 up to $A_e$-codimension $2$ \cite{RR}. 
 For each $\A$-orbit $\eta$, the Thom polynomial is defined to be 
$$tp^\A(\eta) \in \Z[c_1, c_2, c_1', c_2']$$
where $c_i, c_i'$ are Chern classes of relative tangent bundles of source and target, respectively. 
 
 \begin{table}[htb]
$$
\begin{array}{l | c | l }
\mbox{type}  & \codim & \;\mbox{miniversal unfolding}\\
\hline
\hline
\mbox{lips(beaks)} & 3& \;  (x^3+ xy^2+ax,y) \\
\mbox{swallowtail} & 3 & \; (x^4+xy+ax^2,y) \\
\mbox{goose} & 4 & \; (x^3+ xy^3+axy+bx,y) \\
\mbox{gulls}  &4 & \;  (x^4+xy^2+x^5+axy+bx,y) \\
\mbox{butterfly } & 4 &  \; (x^5+ xy+x^7+ax^3+bx^2,y) \\
\mbox{sharksfin } (I_{2,2}^{1,1})  & 4 &  \; (x^2+y^3+ay,  y^2+x^3+bx) \\
\hline
\end{array}
$$
\caption{}
\end{table}

Note that for each of swallowtail, butterfly and $I_{2,2}^{1,1}$, 
the $\A$-orbit is an open dense subset of its $\K$-orbit in $J(2,2)$, 
thus the closures of the $\A$ and $\K$-orbits coincide. 
That means that 
the corresponding $tp^\A$ coincides with $tp$ 
for its $\K$-orbits. 

For other singularities in the list, 
the $\A$-orbit has positive codimension in its $\K$-orbit. 
For instance, look at the case of lips $(x^3+ xy^2, y)$. 
It is $\K$-equivalent to the cusp $A_2$ but not $\A$-equivalent. 
The $\A$-miniversal unfolding $\C^2\times \C, 0 \to \C^2\times \C,0$ with one parmeter $a$ 
gives a $3$-dimensional normal slice 
of the $\A$-orbit of lips type in jet space $J(2,2)$. 
The intersection of the slice with 
the $\K$-orbit of $A_2$ form a smooth curve in the source $\C^2\times \C$ of the unfolding;
The curve is mapped to the cuspidal edge of the critical value set in the target so that 
it is tangent to the plane $\C^2\times \{0\}$ and 
transverse to $\C^2\times \{a\}\,\, (a\not=0)$, 
see Fig. \ref{Lips}.  
\begin{figure}
\includegraphics[clip, width=7.5cm]{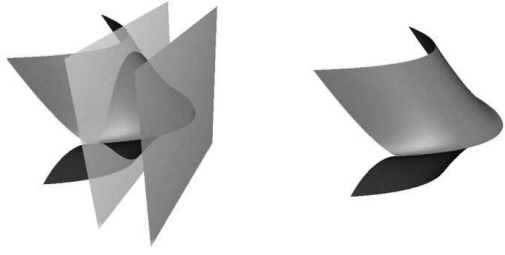} 
\caption{\small Lips and Cuspidal edge}
\label{Lips}
\end{figure}

By the restriction method, we can compute $tp^\A$ for lips, gulls and goose \cite{SSO14}.
There are applications of these formulas on projective algebraic geometry of surfaces. 
Here the normal form of gulls is not weighted homogeneous, but it suffices to consider its $4$-jet 
for computing $tp^\A$, because the closure of the $\A$-orbit is determined by the $4$-jet. 
Note that they can not be expressed in terms of quotient Chern classes. 
On one hand, $tp^\A$ for swallowtail, butterfly and $I_{22}$ are also obtained, 
that coincide with $tp$ for their $\K$-types so that of $1+c_1(f)+\cdots=\frac{1+c_1'+c_2'}{1+c_1+c_2}$.

 \begin{table}[htb]
$$
\begin{array}{l  | l}
\hline 
\mbox{lips} &\;\;   
{-2c_1^3+5c_1^2c_1'-4c_1c_1'^2-c_1c_2+c_2c_1'+c_1'^3}\\
\hline 
\mbox{gulls} &\;\;   
{6c_1^4 -c_1^2c_2-4c_2^2-17c_1^3c_1'+4c_1c_2c_1'+17c_1^2c_1'^2-3c_2c_1'^2} \\
&\;\;
{-7c_1c_1'^3+c_1'^4+2c_1^2c_2'+6c_2c_2'-4c_1c_1'c_2'+2c_1'^2c_2'-2c_2'^2 } \\
 \hline 
 \mbox{goose} &\;\;   
{2c_1^4 +5c_1^2c_2+4c_2^2-7c_1^3c_1'-10c_1c_2c_1'+9c_1^2c_1'^2+5c_2c_1'^2} \\
&\;\;
{-5c_1c_1'^3+c_1'^4-2c_1^2c_2'-6c_2c_2'+4c_1c_1'c_2'-2c_1'^2c_2'+2c_2'^2 } \\
 \hline 
\end{array}
$$
\caption{$tp^\A$ for plane-to-plane germs}
\end{table}
}
\end{exam}

\

\subsection{Thom polynomials for stable multi-singularities}\label{tp_mult}

This subsection is a quick introduction to 
M. Kazarian's theory on Thom polynomials for multi-singularities \cite{Kaz03, Kaz06, Kaz08}. 

\begin{definition}{\rm 
A {\it  multi-singularity} means an ordered set 
$\underline{\eta}:=(\eta_1, \cdots, \eta_r)$ 
of mono-singularities $\eta_i$ of map-germs $\C^m, 0 \to \C^n,0$ 
(especially, we distinguish the first entry $\eta_1$ from others). 
In case of $\kappa=n-m\le 0$, we assume that 
the collection  $\underline{\eta}$ contains no submersion-germs. 
}
\end{definition}

\begin{exam}{\rm 
For instance,  in case of  $(m,n)=(3,3)$,  there are four non-mono stable types; 
Double folds $A_1^2:=A_1A_1$, Triple folds $A_1^3:=A_1A_1A_1$ 
and intersections of fold and cusp $A_1A_2$ and $A_2A_1$. 
The last two types have 
different meanings in {\it source space} 
but the same  in {\it target}, that is indicated by Fig. \ref{A1A2}.

\begin{figure}
\includegraphics[clip,  width=11cm]{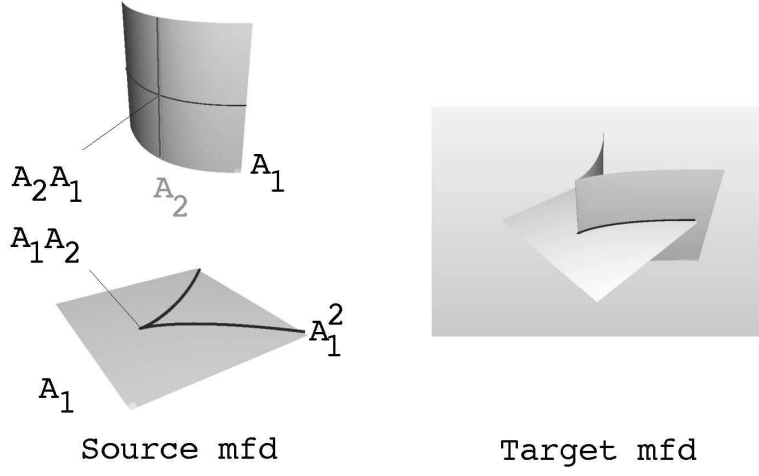} 
\caption{\small $A_1A_2$ and $A_2A_1$ in case of $\kappa=0$}
\label{A1A2}
\end{figure}
}
\end{exam}

For  a stable map $f: M \to N$,  
we set 
$$\underline{\eta}(f)
:=\left\{ \; x_1 \in \eta_1(f) \; \Bigl| \;
\begin{array}{l} 
\exists x_2, \cdots, x_r \in f^{-1}f(x_1)-\{x_1\} \; s.t. \;  x_i\not=x_j \\
(i\not=j) \; \mbox{and  $f$ at $x_i$ is of type $\eta_i$ } 
\end{array}  \right\}$$ 
and call  its analytic closure $\overline{\underline{\eta}(f)} \subset M$ 
the {\it multi-singularity locus of type $\underline{\eta}$ in source};  
The image is  
$$f(\underline{\eta}(f))
:=\left\{ \; y \in N \; \Bigl|  \;
\begin{array}{l} 
\exists x_1, \cdots, x_r \in f^{-1}(y)  \; s.t. \;  x_i\not=x_j  \\
(i\not=j) \; \mbox{and  $f$ at $x_i$ is of type $\eta_i$} 
\end{array}  \right\}$$
and we call the closure $\overline{f(\underline{\eta}(f))} \subset N$  
the {\it  multi-singularity locus of $\underline{\eta}$ in target}. 

The restriction map 
$$f: \overline{\underline{\eta}(f)} \to \overline{f(\underline{\eta}(f)})$$ 
 is finite-to-one: let $\deg_1 \underline{\eta}$ be 
the degree of this map, then 
$$\deg_1 \underline{\eta}
= \mbox{the number of $\eta_1$ appearing in the tuple $\underline{\eta}$. }$$
For instance, $\deg_1 A_1^3=3$, $\deg_1A_1A_1A_2=2$. 

 \begin{rem}{\rm 
For instance, in case of $m=n$, 
$A_1^2(f)$ contains $A_1^k(f)$ of $k \le n$, and  
$\overline{f(A_1^2(f))}-f(A_1^2(f))$ consists of $f(A_1A_2(f))$ and $\overline{f(A_3(f))}$, 
and so on. 
This notional convention might not be so common, 
but it is convenient (economical) for our purpose. 
This is not essential: we usually take the closure in any cases. 
}
\end{rem}

\begin{definition}{\rm 
The {\it Landweber-Novikov class} for proper maps $f: M \to N$ multi-indexed 
by $I=(i_1 i_2\cdots)$  is defined by 
$$s_I=s_I(f)=f_*(c_1(f)^{i_1}c_2(f)^{i_2}\cdots) \;\; \in H^*(N)$$
where $c_i(f)=c_i(f^*TN-TM)$, e.g., 
$$s_0=f_*(1), \; s_i=f_*(c_1^i), \; s_{ij}=f_*(c_1^ic_2^j),  \; 
s_{ijk}=f_*(c_1^ic_2^jc_3^k),  \;  \cdots.$$
For simplicity we often denote $s_I$ to stand for 
its pullback $f^*s_I \; \in H^*(M)$  (i.e., omit the letter $f^*$) 
unless it causes a confusion. 
}
\end{definition}

The following statement has first appeared in M. Kazarian \cite{Kaz03} 
with a topological justification 
using complex cobordism, $h$-principle and Vassiliev's spectral sequence, 
but there has not yet been any rigorous proof  up to the present, as far as the author knows    
-- the proof should be achieved in the context of intersection theory 
of algebraic geometry. 
So precisely saying, this is still a conjecture, see also Remark \ref{remark_Kaz} below.  
On one hand, there are some concrete results  supporting this statement 
 in restrictive cases. Those are mostly due to S. Kleiman's school in 80's 
with techniques using Hilbert schemes and the iteration method, 
see Kleiman \cite{Kleiman81, Kleiman90}, 
also see an unpublished note by Kazarian \cite{Kaz08}: 
For projective maps only with corank one singularities, 
a certain algorithm for computing the multi-singularity loci (stationary multiple point loci) 
has been presented with some actual computations in  small (co)dimensions, 
see  the dissertation of S. Colley \cite{Colley}, for instance. 
There is however a very hard technical difficulty to extend directly this approach to 
general maps having  singularities of corank greater than one. 
Anyway,  in latter chapters, 
we will make use of some concrete computations and arguments 
only for (multi) singularities of $A_k$-types in particularly low dimensions 
 (Example \ref{exam33} below).

\begin{"thm"} {\rm ({\bf Conjecture}  \cite{Kaz03, Kaz06, Kaz08})} \label{kaz_thm}
Given a stable multi-singularity type $\underline{\eta}$ of 
$\C^m \to \C^{m+\kappa}$,  
there exists a unique polynomial in $c_i$ and $s_I$ 
$$tp(\underline{\eta}) \; \in \; 
 \Q [\; c_i,  s_I ; \;  i \ge 1, \; I =(i_1 i_2\cdots)\; ]$$
so that 
 for any proper stable map $f: M \to N$ of relative codimension $\kappa$, 
the locus in source is expressed by 
the polynomial evaluated by 
$c_i=c_i(f)=c_i(f^*TN-TM)$
and  $s_I=s_I(f)=f^*f_*(c^I(f))$: 
$$\Dual [\overline{\underline{\eta}(f)}]
=tp(\underline{\eta}) \;\;  \in H^*(M; \Q).$$
Also the locus in target 
is expressed by a universal polynomial  in $s_I(f)$
$$\textstyle 
\Dual [\overline{f(\underline{\eta}(f)})] = 
tp_{\mbox{\tiny target}}(\underline{\eta}) := 
\frac{1}{\deg_1\underline{\eta}}\, f_*tp(\underline{\eta}) \;\;  \in H^*(N; \Q).$$
\end{"thm"}

\begin{definition}{\rm 
We call $tp(\underline{\eta})$ 
{\it the Thom polynomial of a stable multi-singularity type $\underline{\eta}$} 
and $tp_{\mbox{\tiny target}}(\underline{\eta})$ 
{\it the Thom polynomial of  $\underline{\eta}$ in target}. 
}
\end{definition}

\begin{exam}\label{exam33} {\rm  
In case of relative codimension $\kappa=0, 1$, 
Thom polynomials for multi-singularities of stable maps in low dimensions 
are given in the following Tables \ref{tp_codim0} and \ref{tp_codim1}  
\cite{Kaz03, Kaz08} -- Rim\'anyi's restiction method is also effective 
for computing these polynomials $tp(\underline{\eta})$. 
Those polynomials are also computed in e.g. \cite{Colley} 
within an entirely different approach. 
%
\begin{table}[htb]
$$
\begin{array}{l  | c | l}
\hline 
\mbox{type} & \codim & tp\\
\hline 
A_1 & 1 & c_1 \\
A_2 & 2 & c_1^2 + c_2 \\
A_1^2 & 2 & c_1s_1-4c_1^2-2c_2\\
A_3 & 3 & c_1^3+3c_1c_2+2c_3\\
A_1^3 & 3 & \frac{1}{2}\left(
\begin{array}{l}
c_1s_1^2-4c_2s_1-4c_1s_2-2c_1s_{01}-8c_1^2s_1\\
 +40c_1^3+56c_1c_2+24c_3
 \end{array}
 \right)\\
A_1A_2 & 3 &  c_1s_2+c_1s_{01}-6c_1^3-12c_1c_2-6c_3\\
A_2A_1 & 3 & c_1^2s_1+c_2s_1 -6c_1^3-12c_1c_2-6c_3\\
\hline
\end{array}
$$
\caption{\small $\kappa=0$}
\label{tp_codim0}
\end{table}
\begin{table}[htb]
$$
\begin{array}{l  | c | l}
\hline 
\mbox{type} & \codim & tp\\
\hline 
A_0^2 & 1 & s_0-c_1 \\
A_1 & 2 & c_2 \\
A_0^3 & 2 & \frac{1}{2}(s_{0}^2-s_1-2s_0c_1+2c_1^2+2c_2)\\
A_0A_1 & 3 & s_{01}-2c_1c_2-2c_3\\
A_1A_0 & 3 &  s_0c_2-2c_1c_2-2c_3 \\
A_0^4 & 3 &  
\frac{1}{3!}\left(
\begin{array}{l}
s_0^3-3s_0s_1+2s_2+2s_{01}-3s_0^2c_1+3s_1c_1\\
+6s_0c_1^2+6s_0c_2-6c_1^3-18c_1c_2-12c_3
\end{array}
\right)\\
\hline
\end{array}
$$
\caption{\small $\kappa=1$}
\label{tp_codim1}
\end{table}
}
\end{exam}

\begin{rem}\label{remark_Kaz}
{\rm 
The above `theorem' infers  
a sort of manifestation for an expected modern enumerative theory of singularities -- 
the full theory should involve algebraic cobordisms and relative Hilbert schemes 
within intersection theory. 
In fact, this touches a deep issue: 
For instance, 
the G\"ottsche conjecture (now theorem) states the existence of universal polynomials 
of Chern classes for counting nodal curves on a given projective surface, 
that is actually regarded as a typical example of 
muti-singularity Thom polynomials for $A_1^k$; 
Kontsevich's formula counting rational plane curves (Gromov-Witten invariants) 
also relates to counting curves with some prescribed singularities, 
see \cite{Kaz03, Kaz08}. 
}
\end{rem}

\section{Computing $0$-stable invariants of map-germs}

\subsection{Stable perturbation}\label{0stable}
\begin{figure}
\includegraphics[clip, width=9cm]{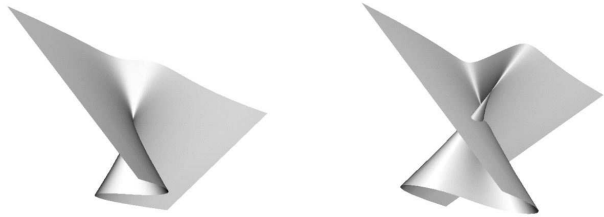}
\caption{\small $H_2$-singularity $(x^3, x^5+xy,y)$ - 
its stable perturbation has two crosscaps and one triple point. }
\label{H2}
\end{figure}

Let $f: \C^m,0 \to \C^n,0$ be a finitely determined map-germ, and 
$\eta$ a stable (mono/multi-)singularity type  of codimension $n$ in the target 
(equivalently, of codimension $m$ in source). 
Take a {\it stable perturbation} 
$$f_t: U \to \C^n\;\; (t\in \Delta \subset \C, \; 0 \in U \subset \C^m)$$
so that $f_0$ is a  representative of $f$ and $f_t$ for $t\not=0$ is a stable map. 
Then $\eta(f_t)$ for $t\not=0$ consists of finitely many isolated points (Fig. \ref{H2}): 
the number is constant for non-zero $t$ and does not depend on the choice of 
stable perturbation (note that if $\eta$ is a mono-stable singularity type,  
it is enough to assume that $f_0$ is $\K$-finite, while 
for multi-singularity type, we need $\A$-finiteness of $f_0$). 
The number of $\eta(f_t)$ is usually called an {\it $0$-stable invariant} 
 of the original germ $f$.

\begin{figure}[H]
\includegraphics[clip, width=8cm]{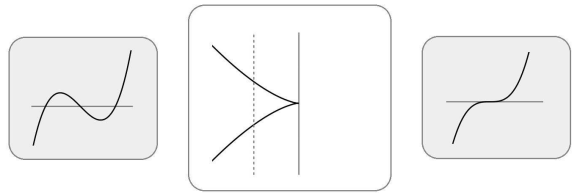}
\caption{\small Discriminant of universal unfolding of $A_2$}
\label{x_cubed}
\end{figure}

Our problem is to compute such a local invariant of map-germs. 
A major prototype is the famous theorem of J. Milnor in the function case ($n=1$): 
The number of Morse singularities arising in a stable perturbation of 
$f: \C^m,0 \to \C,0$  is given by 
the length of the Milnor algebra: 
$$\#  A_1(f_t) = \dim_\C \Ost_{\C^m, 0}/J_f$$
where $J_f$ is the Jacobi ideal. 
For instance, take $f_0=x^3$, then the number of $A_1$-points is 
$\dim_\C \Ost/\langle x^2 \rangle =2$, and 
this is just the degree of the discriminant of the universal unfolding 
$(x,u) \to (x^3-ux,u)$ (Fig. \ref{x_cubed}). 

\begin{rem}{\rm 
In  Fukuda-Ishikawa \cite{FI87} and Gaffney-Mond \cite{GM1, GM2}, 
the formula has been generalized to  
the case of plane-to-plane germs  
for counting the numbers of 
cusps and double fold points in generic perturbation. 
Since then,  several authors,  
Nu\~no Ballesteros, Saia, Fukui, Jorge Perez, Miranda \cite{FBS, MMR, VictorHugo, Miranda} etc, 
have been developing this direction further 
for higher dimensional cases. 
The strategy is as follows. 
 For a mono stable singularity type $\eta$ (e.g. a Thom-Boardman type), 
the first task is to describe the defining ideal of 
the Zariski closure of the corresponding $\K$-orbit (or TB stratum) 
in a jet space of certain order. 
The second task is to determine when the ideal is Cohen-Macaulay: 
if the ideal is CM, 
the algebraic intersection number of the Zariski closure $\overline\eta$ 
and the jet extension $jf_0$ can easily be computed by the length of 
an associated algebra because the higher torsion sheaves vanish.  
If not, one need more tasks 
to deal with the sygyzy for the ideal. 
Counting stable multi-singularities is more involved and indirect. 
The multiple point schemes are studied using Fitting ideals, and 
usually one assume that  the original germ $f$ is of corank one 
in order to make it possible to handle. 
}
\end{rem}

\subsection{Thom polynomial approach} \label{0stable_via_Tp}
We propose a new topological method based on Thom polynomials 
for computing stable invariants for 
{\it weighted homogeneous} map-germs. 
This provides a significantly simpler computation without any corank condition 
and a transparent perspective for the counting problem 
in weighted homogeneous case. 
We consider the non-negative codimensional case, $\kappa=n-m \ge 0$. 

Let $f: \C^m, 0 \to \C^n, 0$ be a weighted homogeneous germ 
with  weights $w_1, \cdots, w_m$ and  degrees $d_1, \cdots, d_n\in \Z_{>0}$, i.e., 
there are diagonal representations of $T=\C^*$ in source and target spaces 
$$
\rho_0=\alpha^{w_1}\oplus \cdots \oplus \alpha^{w_m}, \quad 
\rho_1=\alpha^{d_1}\oplus \cdots \oplus\alpha^{d_n}$$
which stabilizes the map-germ:  $f=\rho_1 \circ f \circ \rho_0^{-1}$. 

Suppose that $f$ is finitely determined. 
Then its $\A_e$-versal unfolding 
$$F: \C^{m+k},0 \to \C^{n+k},0$$ is also weighted homogeneous 
(e.g., see \cite{Wall}). 
Let $r_1, \cdots , r_k$ be the weights of unfolding parameters. 
Note that by the torus action, $f$ and $F$ can be regarded 
as polynomial maps on affine spaces $\C^m$ and $\C^{m+k}$ respectively. 
Let $i_0: \C^m \times \{0\} \hookrightarrow \C^{m+k}$ and 
$\iota_0: \C^n  \times \{0\} \hookrightarrow \C^{n+k}$ be 
natural inclusions. 

Consider a stable mono/multi-singularity type $\underline{\eta}$ of codimension $n$ 
in the target. Of course, $F$ itself is a stable map, 
so we have the singularity loci in source and target of $F$: 
$$\begin{array}{ccccc}
&\C^m &\stackrel{f}{\longrightarrow}& \C^n&\\
& \mbox{\footnotesize $i_0$}\downarrow \;\;&  &\;\; \downarrow \mbox{\footnotesize $\iota_0$}& \\
\underline{\eta}(F) \subset &\C^{m+k} &\stackrel{F}{\longrightarrow}&   \C^{n+k}
& \supset F(\underline{\eta}(F))
\end{array}
$$
Take a generic (non-equivariant) perturbation $\iota_t$ of $\iota_0$ 
by $t \in \C$ sufficiently close to $0$ 
so that $\iota_t$ ($t\not=0$) is transverse to the critical value set of $F$. 
For instance, this is achieved by taking a generic affine transition of the subspace 
$\C^n\times \{0\}$ in  $\C^{n+k}$. 
The fiber product of $\iota_t$ and $F$ defines
 a perturbation of the embedding $i_0$ of the source space, 
 say $i_t: \C^m \to \C^{m+k}$,   
and it hence gives a stable perturbation $f_t$ of the original map $f_0=f$ 
so that $F \circ i_t=\iota_t \circ f_t$. 
The $\underline{\eta}$-locus of $f_t$ in target is  
the intersection of $\iota_t$ with  $F(\underline{\eta}(F))$,  
which consists of finitely many points 
because of the assumption that the codimension of $\underline{\eta}$ 
and the above construction of maps are complementary. 
 
Now, thanks to the torus action, we deal with the global setting 
associated to the above diagram of polynomial maps. 
We introduce three vector bundles over $BT=\Proj^\infty$ 
(or large dimensional projective space) 
by sums of tensor powers of the canonical line bundle $\ell = \Ost(1)$:  
$$E_0:=\oplus_{i=1}^m \Ost(w_i), \quad 
E_1 :=\oplus_{j=1}^n \Ost(d_j), 
\quad 
E'=\oplus_{i=1}^k \Ost(r_i),$$ 
which correspond to representations of $T=\C^*$ 
on the source, target and parameter spaces, respectively. 
Then our weighted homogenous polynomial maps $f$ and $F$ 
yield well-defined universal maps between the total spaces of these vector bundles. 
For simplicity, we denote these universal maps  
by the same notations $f$, $F$, $\iota_0$, $i_0$, 
that would not cause any confusion: 

$$\begin{array}{ccccc}
&E_0 &\stackrel{f}{\longrightarrow}&  E_1 & \\
&\mbox{\footnotesize $i_0$}\downarrow \;\;& & \;\; \downarrow  \mbox{\footnotesize $\iota_0$}& \\
\underline{\eta}(F) \subset  & E_0\oplus E' &\stackrel{F}{\longrightarrow}&  E_1\oplus E' &\supset F(\underline{\eta}(F))
\end{array}
$$

\

\begin{figure}
\includegraphics[clip, width=7.5cm]{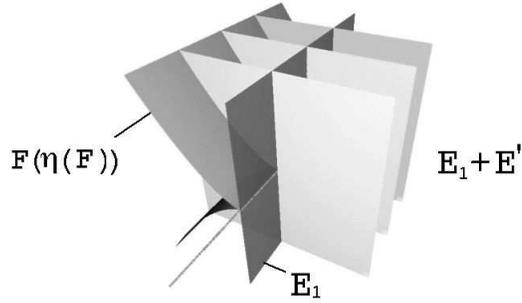}
\caption{\small The target space of the universal stable map $F$}
\label{target_universal_map}
\end{figure}

Perturb the embedding $\iota_0$ of $E_1$ in order to yield a desired  perturbation 
$f_t: E_0 \to E_1$ of the original map $f_0=f$. For instance,  
this is achieved by taking a section $s \in \Gamma(E')$ 
and  
$$\iota_t: E_1 \to E_1 \oplus E', \;\;  
\iota_t(p, v):=\iota_0(v)+t\cdot s(p)$$  
for $p \in BT$, $v \in (E_1)_p$. 
For generic $s$,  the shifted embedding $\iota_t$  is transverse 
to the critical value locus of $F$ in the total space $E_1 \oplus E'$ 
over an open dense set of $BT$. 
The fiber product of $F$ and $\iota_t$ defines deformations 
$f_t: E_0 \to E_1$ and $i_t: E_0 \to E_0\oplus E'$  so that 
$F\circ i_t=\iota_t\circ f_t$ and that 
$f_t:(E_0)_p \to (E_1)_p$ is a stable map for allmost all $p \in BT$.

\begin{figure}
\includegraphics[clip, width=9cm]{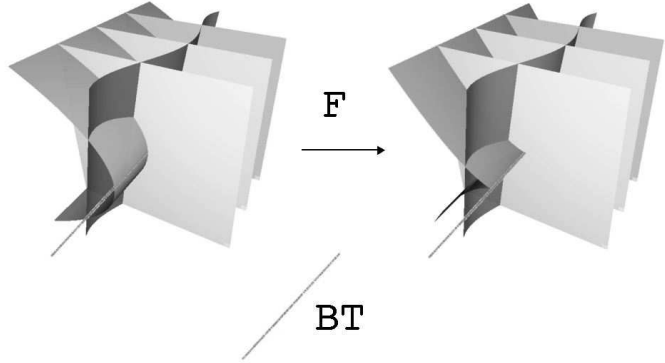}
\caption{\small Perturbation of $i_0$ and $\iota_0$}
\end{figure}

By the pullback via $p: E_1\oplus E' \to BT$ 
we identify 
$$H^*(E_1\oplus E'; \Q)= H^*(BT; \Q)=\Q[[a]],$$
where $a=c_1(\ell)$, the first Chern class of the canonical line bundle. 
The $\underline{\eta}$-type (multi-)singularity loci of $F$ 
defines an $n$-dimensional cocycle in the target total space $E_1\oplus E'$ 
which is expressed by the target Thom polynomial associated to $\underline{\eta}$: 
$$ [\overline{F(\underline{\eta}(F))}] 
= tp_{\mbox{\tiny target}}(\underline{\eta})(F) 
=  h \cdot a^n \qquad (\exists \, h \in \Z). $$
On one hand, 
$E_1 \stackrel{\iota_0}{\hookrightarrow} E_1 \oplus E'$ 
is an embedding of the total spaces with the normal bundle $p^*E'$, 
hence the fundamental cycle defines a $k$-dimensional cocycle in $E_1\oplus E'$ 
which is expressed by the top Chern class of the normal bundle:  
$$\Dual [E_1] =\iota_{0*}(1)=c_{k}(p^*E')=r_1\cdots r_k \cdot a^k.$$
Now our perturbation $\iota_t$ is  transverse to 
the $\underline{\eta}$-locus of stable map $F$ 
and $\iota_t$ is homotopic to $\iota_0$, 
thus the intersection cocycle represents the cohomology cap product 
in $H^*(E_1 \oplus E')$  
$$[\overline{F(\underline{\eta}(F))} \cap \iota_t(E_1)] 
= [\overline{F(\underline{\eta}(F))}] \cdot \Dual [E_1]. $$
Since the intersection cocycle has codimension $m+k$, 
the cycle must be an integer multiple of 
the class represented by the zero section of $E_1 \oplus E'$, 
i.e., the top Chern class $c_{n+k}(E_1\oplus E')$. 
The multiplicity is equal to 
the degree of the projection 
$$p': \overline{F(\underline{\eta}(F))} \cap \iota_t(E_1) \to BT.$$
Looking at generic fiber of $p'$, 
the degree coincides with $\#  \underline{\eta}(f_t)$ in the local setting 
(this number is well-defined by the assumption). Hence we have 
$$\# \underline{\eta}(f_t)=
\frac{tp_{\mbox{\tiny target}}(\underline{\eta})(F) \cdot \iota_{0*}(1)}
{c_{n+k}(E_1\oplus E')}
=\frac{h \cdot r_1\cdots r_k }{d_1\cdots d_n\cdot r_1\cdots r_k}
=\frac{h}{d_1\cdots d_n}$$
(consequently, $h$ is divisible by the product of degrees). 

\begin{rem}\label{quotientChern_counting}{\rm 
Note that  
the quotient Chern classes $c(f_0)$ and $c(F)$ 
are the same (by cancelation of the $E'$ factor): 
$$c(F)=c(f_0)=1+c_1(f_0)+c_2(f_0)+\cdots 
=\frac{\prod (1+d_j a)}{\prod (1+ w_i a)} \; \in \Q[[a]].$$ 
For a mono-singularity $\underline{\eta}=\eta$, 
the Thom polynomial $tp(\eta)$ for $F$ is a polynomial in $c_i(f_0)$, 
so it is computed in terms of weights and degrees. 
For  a multi-singularity $\underline{\eta}$, 
the Thom polynomial 
$tp(\underline{\eta})$ for $F$ 
is a polynomial in $c_i(f_0)$ and $s_I(f_0)$. 
Since $f_0: E_0 \to E_1$ is a proper map (we assume that $m \le n$), 
the (co)homology pushforward $f_{0*}$ is defined. 
The zero locus of $E_0$ is mapped via $f_0$ identically to the zero locus of $E_1$, 
hence 
$$f_{0*}(c_m(E_0))=c_m(E_0)f_{0*}(1)=c_n(E_1)$$
 (\cite[Lem. 4.1]{Kaz06}), 
so we have 
$$
s_0(f_0)=f_{0*}(1)=\frac{d_1\cdots d_n}{w_1\cdots w_m}a^{n-m}, \quad 
s_I(f_0)=c^I(f_0)s_0(f_0). $$ 
Hence $tp(\underline{\eta})$ 
(and  $tp_{\mbox{\tiny target}}(\underline{\eta})$) is written by weights and degrees. 
}
\end{rem}

Thus the following theorem is proved: 

\begin{thm}\label{counting_0_stable_inv}
Let $m \le n$ and 
let $f_0: \C^m,0\to \C^n,0$  be 
an $\A$-finitely determined weighted homogeneous map-germ 
with weight $w_i$ and degree $d_j$. 
Given a stable mono/multi-singularity $\underline{\eta}$ of codimension $n$ in target, 
the corresponding $0$-stable invariant of $f_0$ is computed by 
$$\# \underline{\eta}(f_t) 
= 
\frac{tp_{\mbox{\tiny target}}(\underline{\eta})}
{ d_1\cdots d_n}
= \frac{tp(\underline{\eta})}{{ \deg_1 \underline{\eta}} \cdot w_1\cdots w_m}$$
where numerators stand for the coefficient of $a^m$ and $a^n$ 
of the Thom polynomials in source and target applied to the universal map 
$f_0: E_0 \to E_1$, respectively.  
In particular, for the case of mono-singularity $\underline{\eta}=(\eta_1)$, 
we have $\#{\eta_1}(f_t) = tp(\eta_1)/ w_1 \cdots w_m$. 
\end{thm}

\begin{rem}{\rm 
As seen, we restrict 
the Thom polynomial $tp(\underline{\eta})$ to a more complicated singularity $f=f_0$. 
The resulting class in $H^*(BT)$ is a sort of {\it incidence class} 
introduced by Rim\'anyi \cite{Rimanyi}. 
}
\end{rem}

\begin{rem}{\rm 
If $m >n$, then $f_{0}$ is not proper, 
so the argument about $s_I(f_0)$ in Remark \ref{quotientChern_counting} 
is not available. 
Instead, 
since the restriction of $f_0$ to the critical point set is generically one-to-one, 
hence proper,  the pushforward of the restricted map is defined and computable. 
Then a similar formal computation of Thom polynomials works,  
as pointed out in  \cite[\S 4]{Kaz06}. 
}
\end{rem}

\begin{rem}{\rm 
Not only the $0$-stable invariant but also 
higher stable invariants are defined by 
the degree of the subvariety $\underline{\eta}(f_t)$ 
which has positive dimension. 
Our theorem can also be generalized for computing such stable invariants 
for finite weighted homogeneous germs. 
}
\end{rem}

\subsection{Computation}\label{0stable_computation}

Computing the $0$-stable invariants for $f$ via Tp 
is simply reduced to elementary polynomial algebra, 
i.e., we compute  
$$\#  \eta(f_t)=\frac{tp(\underline{\eta})}{\deg_1 \underline{\eta} \cdot w_1\cdots w_m}$$ 
by substitution.  
Below we demonstrate some computations

 \begin{exam}\label{counting22}{\rm 
$(m,n)=(2,2)$: 
Tp of stable singularities of codimension $2$ are 
$$tp(A_2)=c_1^2+c_2, \quad tp(A_1^2)=c_1s_1-4c_1^2-2c_2.$$
Let  $f: \C^2,0\to \C^2,0$ be a finitely determined weighted homogeneous 
germ with weights $w_1, w_2$ and degrees $d_1, d_2$. 
The quotient Chern class is 
$$c(f)=\frac{(1+d_1a)(1+d_2a)}{(1+w_1a)(1+w_2a)},$$ 
so we get 
\begin{eqnarray*}
c_1&=&(d_1+d_2-w_1-w_2)a, \\
c_2&=&(d_1d_2-d_1w_1-d_2w_1+w_1^2-d_1w_2-d_2w_2+w_1w_2+w_2^2)a^2, \\
s_0&=&\textstyle 
\frac{d_1d_2}{w_1w_2}, \\
s_1&=&
\textstyle s_0 c_1=\frac{d_1d_2}{w_1w_2}(d_1+d_2-w_1-w_2)a. 
\end{eqnarray*}
Substitute them into 
$$\frac{tp(A_2)}{w_1w_2}, \quad \frac{tp(A_1^2)}{2w_1w_2},$$ 
we obtain the $0$-stable invariants of cusp and double folds for $f$:  
\begin{eqnarray*}
&&\# A_2
=\frac{1}{w_1w_2}
\left(
\begin{array}{l}
d_1^2+d_2^2+2w_1^2+3d_1(d_2-w_1-w_2)\\
+3w_1w_2+2w_2^2-3d_2(w_1+w_2)
\end{array}
\right)
\\
&&\# A_1^2=
\frac{1}{2w_1^2w_2^2}
\left(
\begin{array}{l}
d_1d_2(d_1+d_2-w_1-w_2)^2
-4w_1w_2(d_1+d_2\\
-w_1-w_2)^2
-2w_1w_2\{w_1^2+w_1w_2+w_2^2\\
+d_1(d_2-w_1-w_2)-d_2(w_1+w_2)\}
\end{array}
\right).
\end{eqnarray*}
These coincide with Gaffney-Mond's results \cite{GM1}. 
}
 \end{exam}

 \begin{exam}\label{counting23}{\rm 
$(m,n)=(2,3)$: Tp of stable singularities of $\codim 2$ in source are 
$$\textstyle 
tp(A_1)=c_2, \quad tp(A_0^3)=\frac{1}{2}(s_0^2-s_1-2c_1s_0+2c_1^2+2c_2).$$
Expand 
$$c(f)=\frac{(1+d_1a)(1+d_2a)(1+d_3a)}{(1+w_1a)(1+w_2a)},$$ and 
substitute terms into 
$$\frac{tp(A_1)}{w_1w_2}, \quad \frac{tp(A_1^3)}{3w_1w_2},$$
then we obtain the $0$-stable invariants of crosscap and  triple point for $f$:  
\begin{eqnarray*}
&&\# A_1
=\frac{1}{w_1w_2}
\left(
\begin{array}{l}
d_1d_2+(d_1+d_2)d_3-(d_1+d_2+d_3)w_1+w_1^2\\
-(d_1+d_2+d_3-w_1)w_2+w_2^2
\end{array}
\right)
\\
&&\# A_0^3=
\frac{1}{6w_1^3w_2^3}
\left(
\begin{array}{l}
d_1^2d_2^2d_3^2-3d_1d_2d_3w_1w_2(d_1+d_2+d_3\\
-w_1-w_2)+2w_1^2w_2^2\{d_1d_2+(d_1+d_2)d_3\\
-(d_1+d_2+d_3)w_1+w_1^2\\
+(d_1+d_2+d_3-w_1-w_2)^2\\
-(d_1+d_2+d_3-w_1)w_2+w_2^2\}
\end{array}
\right).
\end{eqnarray*}
These numbers coincide with the result in Mond \cite{Mond} obtained 
by a completely different method. 
}
 \end{exam}

 \begin{exam}\label{counting33}{\rm 
$(m,n)=(3,3)$: Tp for stable (multi-)singularities are 
\begin{eqnarray*}
&& tp(A_3)= c_1^3+3c_1c_2+2c_3,  \\
&& tp(A_1A_2)=c_1s_2+c_1s_{01}-6c_1^3-12c_1c_2-6c_3, \\
&&\textstyle 
tp(A_1^3)= \frac{1}{2}\left(
\begin{array}{l}
c_1s_1^2-4c_2s_1-4c_1s_2-2c_1s_{01}-8c_1^2s_1\\
 +40c_1^3+56c_1c_2+24c_3
 \end{array}
 \right).
 \end{eqnarray*}
The corresponding $0$-stable invariants 
for weighted homogeneous finite germs $f: \C^3,0 \to \C^3,0$ 
are computed below. 
Note that our method is valid  for germs $f$ of {\it  any corank}.

\begin{table}
{\small     
\begin{eqnarray*}
&&\textstyle \# A_3=
\frac{1}{w_1 w_2 w_3}((d_1+d_2+d_3-w_1-w_2-w_3)^3+3 (d_1+d_2+d_3-w_1\\&&
-w_2-w_3)(d_1 d_2+(d_1+d_2) d_3-(d_1+d_2+d_3) w_1+w_1^2-(d_1+d_2+d_3\\&&
-w_1)w_2+w_2^2-(d_1+d_2+d_3-w_1-w_2) w_3+w_3^2)+2 (d_1 d_2 d_3-(d_2d_3\\&&
+d_1 (d_2+d_3)) w_1+(d_1+d_2+d_3) w_1^2-w_1^3-(d_1 d_2+(d_1+d_2)d_3\\&&
-(d_1+d_2+d_3) w_1+w_1^2) w_2+(d_1+d_2+d_3-w_1) w_2^2-w_2^3-(d_1d_2\\&&
+(d_1+d_2) d_3-(d_1+d_2+d_3) w_1+w_1^2-(d_1+d_2+d_3-w_1) w_2+w_2^2)w_3\\&&
+(d_1+d_2+d_3-w_1-w_2) w_3^2-w_3^3)). \\&&
\\&&
\textstyle
\# A_1A_2 = \frac{1}{w_1^2 w_2^2 w_3^2}(d_1^4 d_2 d_3+d_1^3 (4 d_2^2 d_3+4 d_2d_3 (d_3-w_1-w_2-w_3)\\&&
-6 w_1 w_2 w_3)-6 w_1 w_2 w_3 (d_2^3+d_3^3-4w_1^3-8 w_1^2 w_2-8 w_1 w_2^2-4 w_2^3\\&&
+5 d_2^2 (d_3-w_1-w_2-w_3)-8 w_1^2 w_3-13w_1 w_2 w_3-8 w_2^2 w_3-8 w_1 w_3^2\\&&
-8 w_2 w_3^2-4 w_3^3-5 d_3^2 (w_1+w_2+w_3)+d_3(8 w_1^2+8 w_2^2+13 w_2 w_3+8 w_3^2\\&&
+13 w_1 (w_2+w_3))
+d_2 (5 d_3^2+8w_1^2+8 w_2^2+13 w_2 w_3+8 w_3^2+13 w_1 (w_2\\&&
+w_3)-13 d_3 (w_1+w_2+w_3)))
+d_1^2(4 d_2^3 d_3+9 d_2^2 d_3(d_3-w_1-w_2-w_3)\\&&
+30 w_1 w_2 w_3 (-d_3+w_1+w_2+w_3)
+d_2(4 d_3^3-30 w_1 w_2 w_3-9 d_3^2 w_1\\&&
+w_2+w_3)+d_3 (5 w_1^2+5 w_2^2+9w_2 w_3+5 w_3^2+9 w_1 (w_2+w_3))))\\&&
+d_1 (d_2^4 d_3+4 d_2^3 d_3(d_3-w_1-w_2-w_3)-6 w_1 w_2 w_3(5 d_3^2+8 w_1^2+8 w_2^2\\&&
+13 w_2 w_3+8w_3^2+13 w_1 (w_2+w_3)-13 d_3 (w_1+w_2+w_3))\\&&
+d_2^2 (4 d_3^3-30 w_1 w_2w_3-9 d_3^2 (w_1+w_2+w_3)+d_3 (5 w_1^2+5 w_2^2+9 w_2 w_3\\&&
+5 w_3^2+9 w_1(w_2+w_3)))+d_2 (d_3^4-4 d_3^3 (w_1+w_2+w_3)\\&&
+78 w_1 w_2 w_3 (w_1+w_2+w_3)+d_3^2(5 w_1^2+5 w_2^2+9 w_2 w_3+5 w_3^2\\&&
+9 w_1 (w_2+w_3))-d_3 (2 w_1^3+2 w_2^3+5w_2^2 w_3+5 w_2 w_3^2+2 w_3^3\\&&
+5 w_1^2 (w_2+w_3)+w_1 (5 w_2^2+87 w_2 w_3+5w_3^2)))))
\end{eqnarray*}
}
\caption{\small $0$-stable invariants (Swallowtail and Fold+Cuspidal edge) for $\C^3,0 \to \C^3,0$. }
\label{table_counting33_1}
\end{table}

\begin{table}
{\small     
\begin{eqnarray*}
&&\textstyle
\#A_1^3=\frac{1}{6 w_1^3 w_2^3 w_3^3}
(d_1^5 d_2^2 d_3^2+3 d_1^4 d_2 d_3 (d_2^2d_3+d_2 d_3 (d_3-w_1-w_2-w_3)\\&&
-4 w_1 w_2 w_3)
-8 w_1^2 w_2^2 w_3^2(-5 d_2^3-5 d_3^3+15 w_1^3+32 w_1^2 w_2+32 w_1 w_2^2+15 w_2^3\\&&
-22 d_2^2 (d_3-w_1-w_2-w_3)+32w_1^2 w_3+54 w_1 w_2 w_3+32 w_2^2 w_3+32 w_1 w_3^2\\&&
+32 w_2 w_3^2+15 w_3^3
+22d_3^2 (w_1+w_2+w_3)-2 d_3 (16 w_1^2+16 w_2^2+27 w_2 w_3\\&&
+16 w_3^2+27 w_1 (w_2+w_3))-2d_2 (11 d_3^2+16 w_1^2+16 w_2^2+27 w_2 w_3+16 w_3^2\\&&
+27 w_1 (w_2+w_3)-27d_3(w_1+w_2+w_3)))+d_1^3 (3 d_2^4 d_3^2+6 d_2^3 d_3^2 (d_3-w_1\\&&
-w_2-w_3)-42d_2 d_3 w_1 w_2 (d_3-w_1-w_2-w_3) w_3+40 w_1^2 w_2^2 w_3^2\\&&
+3 d_2^2d_3 (d_3^3-14 w_1 w_2 w_3-2 d_3^2 (w_1+w_2+w_3)+d_3 (w_1+w_2+w_3)^2))\\&&
+d_1^2(d_2^5 d_3^2+3 d_2^4 d_3^2 (d_3-w_1-w_2-w_3)-176 w_1^2 w_2^2 w_3^2 (-d_3+w_1+w_2\\&&
+w_3)
+3d_2^3 d_3 (d_3^3-14 w_1 w_2 w_3-2 d_3^2 (w_1+w_2+w_3)+d_3 (w_1+w_2+w_3)^2)\\&&
-2d_2 w_1 w_2 w_3 (21 d_3^3-88 w_1 w_2 w_3-45 d_3^2 (w_1+w_2+w_3)+3 d_3(8 w_1^2+8 w_2^2\\&&
+15 w_2 w_3+8 w_3^2+15 w_1 (w_2+w_3)))+d_2^2 d_3 (d_3^4-3d_3^3 (w_1+w_2+w_3)\\&&
+90 w_1 w_2 w_3 (w_1+w_2+w_3)+3 d_3^2 (w_1+w_2+w_3)^2-d_3(w_1^3+3 w_1^2 (w_2\\&&
+w_3)
+(w_2+w_3)^3+3 w_1 (w_2^2+32 w_2 w_3+w_3^2))))+2d_1 w_1 w_2 w_3 (-6 d_2^4 d_3\\&&
-21 d_2^3 d_3 (d_3-w_1-w_2-w_3)+8 w_1w_2 w_3 (11 d_3^2+16 w_1^2+16 w_2^2+27 w_2 w_3\\&&
+16 w_3^2+27 w_1 (w_2+w_3)-27d_3 (w_1+w_2+w_3))-d_2^2 (21 d_3^3-88 w_1 w_2 w_3\\&&
-45 d_3^2 (w_1+w_2+w_3)+3d_3 (8 w_1^2+8 w_2^2+15 w_2 w_3+8 w_3^2\\&&
+15 w_1 (w_2+w_3)))-3 d_2 (2d_3^4-7 d_3^3 (w_1+w_2+w_3)\\&&
+72 w_1 w_2 w_3 (w_1+w_2+w_3)+d_3^2 (8 w_1^2+8w_2^2+15 w_2 w_3+8 w_3^2\\&&
+15 w_1 (w_2+w_3))-d_3 (3 w_1^3+3 w_2^3
+8 w_2^2w_3+8 w_2 w_3^2+3 w_3^3\\&&
+8 w_1^2 (w_2+w_3)+w_1 (8 w_2^2+87 w_2 w_3+8 w_3^2)))))
\end{eqnarray*}
}
\caption{\small $0$-stable invariant (Triple folds) for  $\C^3,0 \to \C^3,0$. }
\label{table_counting33_2}
\end{table}

The iterated Jacobian ideal $J_{111}$ 
defining the $A_3$-locus (i.e, $\Sigma^{1,1,1}$) 
is not Cohen-Macaulay  along $\Sigma^2$
(communication with Nu\~no-Ballesteros, also see \cite{FBS, VictorHugo, Miranda}). 
So the commutative algebra approach requires more hard works, 
while our topological approach is straightforward and gives the right answer. 
For instance,  consider the following map-germ of corank $2$
$$f(x,y,z)=(x^2+y^2+xz, xy, z).$$ 
Substitute weights $(1,1,1)$ and degrees $(1,2,2)$
 into a bit long formula of $A_3$ as noted above, 
 then it returns the correct answer $2$. Namely, 
 this germ has exactly two $A_3$ points in any stable perturbation.  
 On one hand, the length computation gives a wrong number ($\dim \Ost/J_{111}(f)=4$). 
For the same germ, the remaining two formulas in Table answer the number to be $0$,  
that is,  both $A_1A_2$  and  $A_1^3$ points do not appear in stable perturbation. 

For another example of corank $2$, 
$$f(x,y,z)=(x^9+y^2+xz, xy, z),$$ 
 we have 
 $\# A_3=16$, $\# A_1A_2=105$, $\# A_1^3=98$. 
Those numbers coincide with the result in \cite{Miranda}. 
For counting mono-singularity, our formula is valid also for $\K$-finite germs. 
For instance, the germ $(x^2, y^2, z^2)$ has 23 $A_3$ points in its stable perturbation,  
while there has been no way to compute such a number for germs of corank $3$ so far. 
On the other hand,  applying our formula of $A_1A_2$ or $A_1^3$ to non-$\A$-finite germ 
does not make sense.

For germs $f$ of corank one, the counting formula for each singularity 
has a significantly simpler form. 
Put $w_1=d_1$, $w_2=d_2$ and use $w_0, d$ instead of $w_3, d_3$, 
then we recover a result in Marar-Montaldi-Ruas \cite{MMR}: 
\begin{eqnarray*}
&&\# A_3
=\frac{(d-w_0)(d-2w_0)(d-3w_0)}{w_0w_1w_2},\\
&&\#  A_1A_2
=\frac{(d-w_0)(d-2w_0)(d-3w_0)(d-4w_0)}{w_0^2w_1w_2},\\
&&\#  A_1^3
=\frac{(d-w_0)(d-2w_0)(d-3w_0)(d-4w_0)(d-5w_0)}{6w_0^3w_1w_2}. 
\end{eqnarray*}

We emphasize that the most convenient and well-organized expression for general cases 
is the formula in Theorem \ref{counting_0_stable_inv}. 
 }
 \end{exam}


 \begin{exam}\label{counting34}{\rm 
$(m,n)=(3,4)$: 
Tp for stable quadruple points is 
$$\textstyle 
tp(A_0^4)=
\frac{1}{6}\left(
\begin{array}{l}
s_0^3-3s_0s_1+2s_2+2s_{01}-3s_0^2c_1+3s_1c_1\\
+6s_0c_1^2+6s_0c_2-6c_1^3-18c_1c_2-12c_3
\end{array}
\right).$$
The corresponding $0$-stable invariants 
 is given in Table \ref{table_counting34}. We omit other singularity types. 
  
 \begin{table}[b]
 {\small 
\begin{eqnarray*} 
&&\# A_0^4=\textstyle
\frac{1}{6 w_1^4 w_2^4 w_3^4}
(d_1^3 (d_2^3 d_3^3 d_4^3-6 d_2^2 d_3^2 d_4^2w_1 w_2 w_3+11 d_2 d_3 d_4 w_1^2 w_2^2 w_3^2\\&&
-6 w_1^3 w_2^3 w_3^3)-6w_1^3 w_2^3 w_3^3 (d_2^3+d_3^3+d_4^3-6 d_4^2 w_1+11 d_4 w_1^2-6 w_1^3-6 d_4^2w_2\\&&
+17 d_4 w_1 w_2-11 w_1^2 w_2+11 d_4 w_2^2-11 w_1 w_2^2-6 w_2^3\\&&
+6 d_3^2(d_4-w_1-w_2-w_3)+6 d_2^2 (d_3+d_4-w_1-w_2-w_3)-6 d_4^2 w_3\\&&
+17 d_4w_1 w_3-11 w_1^2 w_3+17 d_4 w_2 w_3-17 w_1 w_2 w_3\\&&
-11 w_2^2 w_3+11 d_4w_3^2
-11 w_1 w_3^2-11 w_2 w_3^2-6 w_3^3\\
&&+d_2 (6 d_3^2+6 d_4^2+11 w_1^2+17 w_1w_2+11 w_2^2
+17 d_3 (d_4-w_1-w_2-w_3)\\
&&+17 w_1 w_3+17 w_2 w_3+11 w_3^2-17 d_4(w_1+w_2+w_3))+d_3 (6 d_4^2+11 w_1^2\\
&&+11 w_2^2+17 w_2 w_3+11 w_3^2+17 w_1(w_2+w_3)-17 d_4 (w_1+w_2+w_3)))\\
&&-6 d_1^2 w_1 w_2 w_3 (d_2^3 d_3^2d_4^2
-6 w_1^2 w_2^2 w_3^2 (-d_3-d_4+w_1+w_2+w_3)\\
&&+d_2^2 d_3 d_4 (d_3^2d_4+d_3 d_4 (d_4-w_1-w_2-w_3)-5 w_1 w_2 w_3)\\&&
+d_2 w_1 w_2 w_3(-5 d_3^2 d_4+6 w_1 w_2 w_3+5 d_3 d_4 (-d_4+w_1+w_2+w_3)))\\
&&+d_1w_1^2 w_2^2 w_3^2 (11 d_2^3 d_3 d_4+6 d_2^2 (5 d_3^2 d_4+5 d_3 d_4 (d_4-w_1-w_2-w_3)\\
&&-6w_1 w_2 w_3)-6 w_1 w_2 w_3 (6 d_3^2+6 d_4^2+11 w_1^2+17 w_1 w_2+11 w_2^2\\
&&+17d_3 (d_4-w_1-w_2-w_3)+17 w_1 w_3+17 w_2 w_3+11 w_3^2\\
&&-17 d_4 (w_1+w_2+w_3))+d_2(11 d_3^3 d_4+30 d_3^2 d_4 (d_4-w_1-w_2-w_3)\\
&&+102 w_1 w_2 w_3 (-d_4+w_1+w_2+w_3)+d_3(11 d_4^3-102 w_1 w_2 w_3\\
&&-30 d_4^2 (w_1+w_2+w_3)
+d_4 (19 w_1^2+19 w_2^2+30w_2 w_3+19 w_3^2\\&&
+30 w_1 (w_2+w_3)))))).
\end{eqnarray*}
}
\caption{\small Quadruple points for $\C^3,0 \to \C^4,0$}
\label{table_counting34}
\end{table}

For example, consider the map-germ of corank $2$ 
$$\Hat{A}_k: (x, y^k+xz+x^{2k-2}y, yz, z^2+y^{2k-1}), $$ 
then the number of quadruple points is $\frac{8}{3}(k-1)^2(k^3-5k^2+9k-6)$. 

For germs $f$ of corank one, it holds that 
$$\# A_0^4
=\frac{(d_1-w_0)(d_1-2w_0)(d_1-3w_0)(d_2-w_0)(d_2-2w_0)(d_2-3w_0)}{6w_0^4w_1w_2}. 
$$

}
\end{exam}

\section{Image and discriminant Chern classes}

\subsection{Izumiya-Marar formula} \label{Izumiya_Marar}

\begin{figure}
\includegraphics[clip, width=10cm]{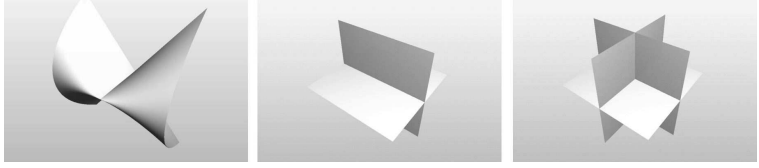} 
\caption{\small Crosscap, double points and triple points in the target space 
of $2$-to-$3$ maps}
\label{image2-3}
\end{figure}

To grasp the main idea quickly, for a moment let us consider 
 a $C^\infty$ stable map 
 from a closed (real) surface $M$ into a (real) $3$-manifold $N$. 
Look at its image singular surface $f(M) \subset N$. 
Stable singularities are of type $A_1$, $A_0^2$ and $A_0^3$ (Fig.\ref{image2-3}).

\begin{thm} \label{IMthm} 
{\rm (Izumiya-Marar  \cite{IM}, cf. \cite{Sedykh})}  
For a $C^\infty$ stable map $f: M^2 \to N^3$, being $M$ compact without boundary, 
the Euler characteristic of the image singular surface 
satisfies the following  formula: 
$$\chi(f(M))=\chi(M)+\frac{1}{2}\#  C + \#  T$$
where $C$ and $T$ are the sets of crosscaps and of triple points in target, respectively. 
\end{thm}

\proof 
Recall that  in the source space $M$, 
\begin{eqnarray*} 
&&A_1(f)=\mbox{the critical point set of $f$}\\
&&A_0(f)=\mbox{the regular point set of $f$} \\
&&A_0^2(f)={\{}\; x \in A_0(f)\; | \;  
\exists\; x' \in A_0(f), \; x'\not=x, \; f(x)=f(x')\; {\}},\\
&&A_0^3(f)={\{}\; x \in A_0(f) \; | \; \exists \, x', x'' \in A_0(f) \cap f^{-1}f(x), \;  
\mbox{\small $x, x',x''$ distint }
 {\}}. 
\end{eqnarray*}
By the definition, $A_0^3 \subset A_0^2 \subset A_0$ and 
the closure $\overline{A_0^2}=A_0^2 \sqcup A_1$. 
Set 
$$A_0^{2\circ}:=A_0^2-A_0^3, \qquad A_0^\circ:=A_0-A_0^2,$$
and 
$$R:=f(A_0^\circ), \;\; 
D:=f(A_0^{2\circ}), \;\; T:=f(A_0^{3}), \;\; C:=f(A_1), $$
then $f$ is stratified by 
$$
\xymatrix{
M=A_0^\circ  \sqcup A_0^{2\circ} \sqcup A_0^3  \sqcup A_1 
\ar[r]^f & R \sqcup D \sqcup T \sqcup C=f(M).
}$$
Obviously, 
\begin{eqnarray*}
&&\jeden_{f(M)}
=\jeden_R+\jeden_D+\jeden_T+\jeden_C, \qquad  \\
&&\jeden_{\overline{A_0^2}}
=\jeden_{A_0^{2\circ}}+\jeden_{A_0^3}+\jeden_{A_1}, \\
&& f_*\jeden_M
=f_*(\jeden_{A_0^\circ}+\jeden_{A_0^{2\circ}}+\jeden_{A_0^3}+\jeden_{A_1})
=\jeden_R+2 \jeden_D+3\jeden_T+\jeden_C, 
\end{eqnarray*}
and a simple computation shows 
\begin{equation} \label{IMformula}
\textstyle
\jeden_{f(M)}=
f_*\left(\jeden_M-\frac{1}{2}\jeden_{\overline{A_0^2}}-\frac{1}{6}\jeden_{A_0^3}+\frac{1}{2}\jeden_{A_1}\right). 
\end{equation}
Take the integration of constructible functions: 
$$\textstyle
\chi(f(M))={\dis\int_N} \jeden_{f(M)} 
={\dis\int_M} 
\left(\jeden_M-\frac{1}{2}\jeden_{\overline{A_0^2}}-\frac{1}{6}\jeden_{A_0^3}+\frac{1}{2}\jeden_{A_1}\right). 
$$
Now we speak about real geometry: since $\overline{A_0^2}$ is a union of immersed curves 
whose double point set is just $A_0^3$, we have 
 $$\chi(\overline{A_0^2})+\chi(A_0^3)=\chi(\mbox{disjoint circles})=0.$$
Hence the integral is rewritten as follows: 
$$\textstyle 
\chi(f(M))=\chi(M)+(\frac{1}{2}-\frac{1}{6})\cdot 3\# T + \frac{1}{2}\# C
=\chi(M)+\# T+\frac{1}{2}\# C.$$ 
This competes the proof.  \qed 

\
 
Notice that the above equality (\ref{IMformula}) is shown 
by using only the combinatorics of adjacencies of singularities, 
thus it is valid for complex singularities as well. 
From now on, let us switch into the complex case. 
We assume that $M, N$ are 
compact complex manifolds of dimension $2, 3$, respectively, 
and $f: M \to N$ is a holomorphic map 
which admits only (mono/multi-)stable singularities 
(in other words, $f$ is a normalization of a singular surface in $N$ 
having ordinary singularities).  
Put 
$$\textstyle 
\imageconst:=\jeden_M-\frac{1}{2}\jeden_{\overline{A_0^2}}-\frac{1}{6}\jeden_{A_0^3}+\frac{1}{2}\jeden_{A_1} \in \mathcal{F}(M)$$
and apply the CSM class transformation to the equality (\ref{IMformula}) 
($f$ is now proper), then we have 
$$\textstyle
C_*(\jeden_{f(M)})
=f_*C_*(\imageconst).$$
We think of this class in $H^*(N)$ via the Poincar\'e dual 
and omit the notation $\Dual$. 
Note that 
$$
\chi(f(M))
 =\int_N C_*(\jeden_{f(M)})
= \int_N f_*C_*(\imageconst) 
= \int_M C_*(\imageconst).
$$
Look at each term in 
$$\textstyle
C_*(\imageconst)=C_*(\jeden_M)-{\frac{1}{2}} C_*(\jeden_{\overline{A_0^2}})-{\frac{1}{6}}C_*(\jeden_{A_0^3})+{\frac{1}{2}}C_*(\jeden_{A_1})\;\; \in H^*(M).$$

\begin{itemize}
\item the normalization of CSM class: 
$$C_*(\jeden_M)=c(TM),$$
\item  $A_1$-locus (crosscaps) is  finite: It is given by $tp$ for $A_1$ ($\kappa=1$)  
$$\textstyle C_*(\jeden_{A_1})=[A_1]=tp(A_1)=c_2 \; (=c_2(f^*TN-TM)),$$
\item Triple point locus in $M$ is also finite: It is given by $tp$ for $A_0^3$  ($\kappa=1$) 
$$\textstyle C_*(\jeden_{A_0^3}) =[A_0^3]=tp(A_0^3)
=\frac{1}{2}(s_0^2-s_1-2c_1s_0+2c_1^2+2c_2),$$
\item Double point curve $\overline{A_0^2}$ in $M$: 
The dual to the CSM class consists of $1$ and $2$-dimensional components in cohomology $H^*(M)$. 
The first component is the fundamental class of the curve, 
thus it is given by $tp$ for $A_0^2$  ($\kappa=1$), 
while the second component corresponds to the Euler characteristics, 
which is easily computed 
using the fact that the curve has only nodes at $A_0^3$-points: 
\begin{eqnarray*} 
&&C_*(\jeden_{\overline{A_0^2}})
=[\overline{A_0^2}]+h.o.t=tp(A_0^2)+h.o.t\\
&&=\textstyle (s_0-c_1)+ \left\{c_1(TM) (s_0-c_1) 
+\frac{1}{2}(-s_0^2-s_1+2c_1s_0+2c_2)\right\}.
\end{eqnarray*}
\end{itemize}

Summing up those classes, we obtain 
 a universal expression of complex version of the Izumiya-Marar formula:

\begin{prop}\label{complexIM}
Given a stable map $f: M^2 \to N^3$ of compact complex manifolds. 
Then it holds that 
$$\chi(f(M))=\frac{1}{6}  {\dis\int_M} \left( 
\begin{array}{l}
3c_1(TM) c_1+6c_2(TM)-3c_1(TM) s_0\\
-c_1^2-c_2-2c_1s_0+s_0^2+2s_1
\end{array}
 \right) $$ 
where 
$c_i=c_i(f^*TN-TM), \; s_0=f^*f_*(1), \; s_1=f^*f_*(c_1)$. 
\end{prop}

 \begin{exam}\label{enriques}{\rm ({\bf A classical formula of Enriques})  
 Let $X$ be a projective surface of degree $d$ in $\Proj^3$ 
 having only ordinary singularities, 
 i.e., crosscap ($A_1$) and normal crossings. 
Denote by $\# C$ the number of crosscaps, 
by $\# T$ the number of triple points, and by $\delta$ 
the degree of the double point curve  of $X \subset \Proj^3$. 
Let us take a normalization of $X$; 
then we have a proper stable map $f: M \to N=\Proj^3$ so that 
$M$ is non-singular and 
the image is just the singular surface $X$ (cf. \cite{Mather}). 
It follows from a classical formula of Enriques that 
the Chern numbers of $M$ are expressed by 
\begin{eqnarray*}
\textstyle \int_M c_1(TM)^2&=&d(d-4)^2 -(3d-16) \delta + 3 T - C, \\
\textstyle \int_M c_2(TM)&=&d(d^2-4d+6) - (3d-8) \delta + 3 C -2T,
\end{eqnarray*}
and 
$f_*c_1(TM)=(d(4-d)+2\delta)a^2$, 
where $a=c_1(\Ost(1))$ the divisor class 
(cf. \cite{Tsuboi}). 
Notice that these formulas are quite easily obtained from 
Thom polynomials: 
In fact, 
$$Ca^3=f_*tp(A_1), \;\; 3Ta^3=f_*tp(A_0^3), \;\; 2\delta a^2=f_*tp(A_0^2), $$
while the target Thom polynomials are written in Landweber-Novikov classes, 
hence their degrees are written by Chern numbers of $M$ and $d$; 
Therefore, the Chern numbers can be written by $C, T, \delta$ and $d$, that recovers the above classical formulas. 
Now let us substitute the Chern numbers into the formula in Proposition \ref{complexIM}, 
then we have 
$$\textstyle \chi(X)=d(d^2-4d+6)+2(2-d) \delta + T -\frac{3}{2}C.$$
 }
 \end{exam}

\subsection{Image Chern class  for stable maps} \label{Image}

Universal expression of the Euler characteristics of the image in Proposition \ref{complexIM}
should be given in a more general form 
for stable maps $f: M^m \to N^{m+1}$ ($m\ge 1$) between complex manifolds. 
In fact, 
our universal formula has a particularly {\it well-structured} form 
(Theorem \ref{imageSSM} and Corollary \ref{imageEuler} below). 

\

\noindent 
{\bf M\"obius inverse formula for the adjacency poset}: 
Recall the adjacency relation of multi-singularities both in source and target: 
The diagram of source multi-singularities of $m$-to-$(m+1)$ maps is 
$$\xymatrix{
A_0 \ar[r]  & A_0^2 \ar[r] \ar[rd] &A_0^3  \ar[r]  \ar[rd] & A_0^4  \ar[r] \ar[rd]& \cdots\\
& & A_1 \ar[r]  & A_1A_0\, \&\, A_0A_1\ar[r]& \cdots \\
}
$$
where the arrow $\underline{\eta} \to \underline{\xi}$ 
means that $\underline{\xi}$ is contained the closure of $\underline{\eta}$. 
That makes the set of all multi-singularity types to be a poset (partially ordered set). 

For a multi-singularity type $\underline{\eta}$ and a stable map $f:M \to N$, set 
$$\underline{\eta}^\circ(f)
:=\overline{\underline{\eta}(f)} - \sqcup \; \underline{\xi}(f) \;\; \subset M$$
where the union runs over all $\underline{\xi} \; (\not=\underline{\eta})$ 
with $\underline{\eta} \to \underline{\xi}$. 

The stratum $\underline{\eta}^\circ(f)$ is mapped to 
its image $f(\underline{\eta}^\circ(f))$ as 
a $\deg_1 \underline{\eta}$-to-one covering, and 
the image does not depend on the order of entries   
of the tuple $\underline{\eta}$, e.g., $f(A_0A_1)^\circ(f)=f(A_1A_0)^\circ(f)$. 
Then the source $M$ breaks into the disjoint union of strata $\underline{\eta}^\circ(f)$ 
and the target $N$ is decomposed into the corresponding image strata, 
that is, $f: M \to N$ is stratified by those locally closed multi-singularity loci 
in source and target. 

Then the constant function $\jeden_{f(M)}$ of the stable image  is 
 written by the sum of 
$f_*\jeden_{\underline{\eta}^\circ(f)}$ with some rational coefficients. 
Therefore,  by the exclusion-inclusion principle, 
the {\it M\"obius inverse formula} for this poset expresses 
the function $\jeden_{f(M)}$ by 
the pushforward via $f_*$ of a certain linear combination of 
constant functions of {\it the closure} 
$\overline{\underline{\eta}(f)} \; (=\overline{\underline{\eta}^\circ(f)})$ 
with rational coefficients. 
 Namely, extending the same procedure as in the proof of Theorem \ref{IMformula} 
 to more general case involving strata of higher codimension, 
 we obtain a constructible function on the source space $M$ 
 having a generalized form of  (\ref{IMformula}): 
\begin{eqnarray*}
\imageconst &=&\textstyle  
\jeden_{\overline{A_0}}
-\frac{1}{2}\jeden_{\overline{A_0^2}}
-\frac{1}{6}\jeden_{\overline{A_0^3}}+\frac{1}{2}\jeden_{\overline{A_1}} \\
&&\textstyle
\qquad  -\frac{1}{12}\jeden_{\overline{A_0^4}} 
+ \frac{1}{6}\jeden_{\overline{A_0A_1}}-\frac{1}{3}\jeden_{\overline{A_1A_0}}
+ \cdots 
\end{eqnarray*}
so that 
$$f_*(\imageconst)=\jeden_{f(M)}.$$
Notice that this constructible function depends only on the classification of stable multi-singularities. 

\begin{definition}
{\rm 
We call the CSM class 
$$C_*(\jeden_{f(M)})=f_*C_*(\imageconst) \in H^*(N)$$ 
{\it the image Chern class of stable maps $f:M \to N$. }
}
\end{definition}

For Morin maps $M^m \to N^{m+1}$, that is, 
stable maps having only corank one singularities,  
the local structures of $A_\mu$ and their multi-singularities are well-understood,  
e.g., stable maps with $m \le 5$ are Morin maps (cf. \cite{Colley, Kaz08}). 
In that case we can prove the following theorem -- the key point here is again 
the property of the Segre-SM class for the transverse pullback in Proposition \ref{segre}.  
Conjecturally the theorem would hold for any dimension 
 and for any stable maps, that is, there must be the Segre-SM class version of 
Theorem \ref{kaz_thm}, see Remark \ref{multising_SSM}. 

\begin{thm} \label{imageSSM}
There is a polynomial $tp^{\SM}(\imageconst)$ 
 in the quotient Chern class $c_i=c_i(f^*TN-TM)$ 
and  the Landweber-Noviknov class $s_I=f^*f_*(c^I)$ so that 
$$
C_*(\imageconst)
=c(TM) \cdot tp^{\SM}(\imageconst) \quad \in H^*(M)
$$
for any proper stable maps $M^m \to N^{m+1}\; (m\le 5)$: 
The low degree terms are given by 
\begin{eqnarray*}
tp^{\SM}(\imageconst)
&=&\textstyle 1+\frac{1}{2}(c_1-s_0)\\
&&\textstyle +\frac{1}{6}(s_0^2+2s_1-2c_1s_0-c_1^2-c_2)\\
&&\textstyle +\frac{1}{24}\left(
\begin{array}{l}
\textstyle 2c_1^3-10c_1c_2+2c_1^2s_0+2c_2s_0+3c_1s_0^2\\
\textstyle  -s_0^3+14s_{01}+5c_1s_1-5s_0s_1-6s_2 
\end{array}
\right) \\
&&+ \cdots.\\
\end{eqnarray*}
\end{thm}

\begin{rem}\label{rem_imageSSM}{\rm 
Note that for a stable map $f: M \to N$, 
$$tp^{\SM}(\imageconst) = s^{\SM}(\imageconst, M) \in H^*(M).$$
The above theorem implies that 
the Segre-SM class of the image $f(M)$ 
in the target space 
$$tp^{\SM}(\jeden_{f(M)})
:=s^{\SM}(\jeden_{f(M)}, N)  \;\; \in H^*(N)$$
is universally expressed  in terms of the Landweber-Novikov classes $s_I(f)$. 
In fact,  
\begin{eqnarray*}
tp^{\SM}(\jeden_{f(M)})&=& c(TN)^{-1}C_*(\jeden_{f(M)}) \\
&=& c(TN)^{-1}C_*f_*(\imageconst)\\
&=&c(TN)^{-1}f_*C_*(\imageconst)\\
&=&c(TN)^{-1}f_*(c(TM)\cdot tp^{\SM}(\imageconst))\\
&=& f_*(c(f)^{-1}\cdot tp^{\SM}(\imageconst)), 
\end{eqnarray*} 
hence 
\begin{eqnarray*}
tp^{\SM}(\jeden_{f(M)})
&=& \textstyle 
s_0 - \frac{1}{2}(s_0^2+s_1) \\
&&\textstyle +\frac{1}{6}(s_0^3-7s_{01}+3s_0s_1+2s_2) \\
&& \textstyle-  \frac{1}{24} \left(
\begin{array}{l}
\textstyle s_0^4+6 s_0^2 s_1 -28 s_0 s_{01}+8 s_0 s_2\\
\textstyle +24 s_{001}+3 s_1^2-30 s_{11}+6 s_3
\end{array}\right) + \cdots.
\end{eqnarray*} 
Note that 
$$C_*(\jeden_{f(M)})
=c(TN) \cdot tp^{\SM}(\jeden_{f(M)}) \; \in H^*(N)$$
is written in the target Chern class $c_i(TN)$ and the Landweber-Novikov classes. 
}
\end{rem}

\begin{definition}{\rm 
We call the universal Segre-SM classes 
$tp^{\SM}(\imageconst)$ and  $tp^{\SM}(\jeden_{{\rm image}})$ 
the {\it source and target higher Thom polynomials} 
for the image of stable maps, respectively. 
}
\end{definition}

In particular we obtain a more general statement of Proposition \ref{complexIM}: 
\begin{cor}\label{imageEuler}
The Euler characteristic of the image of $f: M^m \to N^{m+1}$ 
is expressed by 
$$
\chi(f(M))
=\int_M c(TM) \cdot  tp^{\SM}(\imageconst) 
= \int_N c(TN) \cdot  tp^{\SM}(\jeden_{f(M)}). 
$$
\end{cor}
 
\begin{rem}{\rm 
We emphasize that the above image Euler number formula (Corolloary \ref{imageEuler}) 
has a particularly {\it well-structured form}. 
The second degree term of $c(TM) \cdot  tp^{\SM}(\imageconst)$ 
is just the Euler characteristic of the image of stable maps from a surface into $3$-fold, 
that is exactly Proposition \ref{complexIM}, 
and the third degree term expresses the Euler characteristic of the image of 
 stable maps from $3$-fold into $4$-fold, ... and so on. 
 Classically, those invariants were separately considered, but 
they are in fact mutually related in a very convenient way. 
}
\end{rem}

Notice that 
\begin{eqnarray*}
tp^{\SM}(\imageconst)&=&\textstyle tp^{\SM}(\jeden_M-\frac{1}{2}\jeden_{\overline{A_0^2}}
-\frac{1}{6}\jeden_{\overline{A_0^3}}+\cdots)\\
&=&
\textstyle 1-\frac{1}{2}tp^{\SM}(\overline{A_0^2}) 
- \frac{1}{6}tp^{\SM}(\overline{A_0^3})+\cdots.
\end{eqnarray*}
Thus, to obtain the explicit form of $tp^{\SM}(\imageconst)$ 
in Theorem \ref{imageSSM}, 
we compute the Segre-SM classes $tp^{\SM}$ for 
the closure of individual singularity types 
$$
\jeden_{\overline{A_0^2}}, \;\; 
\jeden_{\overline{A_0^3}}, \;\; 
\jeden_{\overline{A_1}}, \;\; 
\jeden_{\overline{A_0^4}}, \;\; 
\jeden_{\overline{A_0A_1}}, 
\;\; \jeden_{\overline{A_1A_0}},  \;\; \cdots.$$
They are polynomials in $c_i$ and $s_I$, which are in Table \ref{table_tpSM}  
 up to degree $3$. 
 To get them, the method in \S \ref{symmetry} 
 is effective, see Example \ref{SSM_1_ex}. 
The locus of some singularity type in the source and target 
might be {\it non-reduced}, 
but the CSM class depends only on the underlying reduced scheme by definition. 

\begin{table}
\begin{eqnarray*}
tp^{\SM}(\overline{A_0^2})
&=& \textstyle (s_0 -c_1)+\frac{1}{2}(2c_2+2c_1s_0-s_0^2-s_1)\\
&&\textstyle + \frac{1}{6}\left(12 c_1c_2-3 c_1 s_0^2 
 -3 c_1 s_1-6 c_2 s_0+6 c_3+s_0^3\right.\\
&&\textstyle 
\left. 
+3 s_0 s_1-7 s_{01}+2 s_2\right) + \cdots\\
 tp^{\SM}(\overline{A_0^3})
 &=& 
 \textstyle
 \frac{1}{2}\left(2 c_1^2-2 c_1 s_0+2 c_2+s_0^2-s_1\right) 
 +\frac{1}{6} \left(-6 c_1^2 s_0\right.\\
&& \textstyle 
\left.-18 c_1 c_2+6 c_1s_0^2  
   -18 c_3-2 s_0^3+5 s_{01}+2 s_2\right) + \cdots \\
tp^{\SM}(\overline{A_1})
 &=& c_2 - (c_1c_2+c_3) + \cdots
 \\
tp^{\SM}(\overline{A_0^4})
&=&
\textstyle 
   \frac{1}{6} \left(-6 c_1^3+6 c_1^2 s_0-18 c_1 c_2
   -3 c_1 s_0^2+3 c_1 s_1+6 c_2 s_0  \right.\\
&& \textstyle 
\left.
  -12 c_3 +s_0^3-3 s_0 s_1+2 s_{01}+2 s_2\right) + \cdots \\
tp^{\SM}(\overline{A_0A_1})
&=&
\textstyle 
(s_{01}-2c_1c_2-2c_3) + \cdots\\
tp^{\SM}(\overline{A_1A_0})
&=& 
\textstyle 
(s_0c_2-2c_1c_2-2c_3) + \cdots. \\
\end{eqnarray*}
\caption{Universal SSM class for the closure of several singularity types 
in case of $\kappa=1$. }
\label{table_tpSM}
\end{table}

As a byproduct, 
other type image Chern classes, e.g.,  
$C_*(\jeden_{\overline{f(A_0^k(f))}})$ of the {\it $k$-th multiple point locus} in target, 
$C_*(\jeden_{\overline{f(A_1(f))}})$ of the singular value set, 
... etc are also obtained in entirely the same way. 

For instance, 
there is a constructible function $\imageconst(2)$ on the source 
\begin{eqnarray*}
\imageconst(2) &=& \textstyle 
\frac{1}{2}\jeden_{\overline{A_0^2}}
-\frac{1}{6}\jeden_{\overline{A_0^3}}
+\frac{1}{2}\jeden_{\overline{A_1}} \\
&&\qquad \textstyle 
-\frac{1}{12}\jeden_{\overline{A_1^4}}
+\frac{1}{6}\jeden_{\overline{A_0A_1}}
-\frac{1}{3}\jeden_{\overline{A_1A_0}}+\cdots
\end{eqnarray*}
so that 
$$f_*(\imageconst(2))=\jeden_{\overline{f(A_0^2(f))}}.$$
Hence we have the following theorem: 

\begin{thm}
The CSM class of 
the double point locus in the target manifold,  
$\overline{f(A_0^2(f))} \subset N$, of 
stable maps $f: M^m \to N^{m+1}$ is universally expressed by 
$$C_*(\jeden_{\overline{f(A_0^2(f))}})
=f_*(c(TM)\cdot tp^{\SM}(\imageconst(2))) \;\; \in H^*(N)$$
where 
\begin{eqnarray*}
&&tp^{\SM}(\imageconst(2))\\
&&=\textstyle
 \frac{1}{2}(s_0-c_1)+\frac{1}{6}(-c_1^2+5c_2+4c_1s_0-2s_0^2-s_1)\\
&&  \quad \textstyle 
+ \frac{1}{24}\left(
\begin{array}{l}2c_1^3+38c_1c_2+24c_3+2c_1^2s_0-22c_2s_0-9c_1s_0^2\\
+3s_0^3-14s_{01}-7c_1s_1+7s_0s_1+2s_2
\end{array}
\right) + \cdots. 
\end{eqnarray*}
In particular, the Euler characteristics is given by 
$$\chi(\overline{f(A_0^2(f))})=\int_M c(TM) \cdot tp^{\SM}(\imageconst(2)).$$ 
\end{thm}

\begin{exam}\label{SSM_1_ex}
{\rm 
To compute the universal SSM classes, 
the way described in \S \ref{symmetry} for mono-singularity types 
works also for multi-singularity types. 
As an example, let us compute  
the third degree term $tp_3^{\SM}(\overline{A_0^2})$ of $c_i$ and $s_I$ 
for the double point locus of stable maps with codimension $\kappa=1$. 
There are 11 unknown coefficients, and all of them are determined 
by restricting it to 
mono/multi-singularity types of codimension $3$ in the source space. 
For instance, we shall seek for the restriction equation  
at each of types $A_0A_1$ and $A_1A_0$ for $3$-to-$4$ maps. 
Take the pair $f=f_1 \coprod f_2$ of germs with the same target $\C^4$ 
$$f_1: (x,y,z) \mapsto (x,y^2,xy,z), \quad f_2: (u,v,w) \mapsto (u,v,w,0),$$
where $f_1$ is of type $A_1A_0$ and $f_1$ is of type $A_0A_1$. 
The $3$-dimensional torus $T=(\C^*)^3$ acts on the sources of $f_1$ and $f_2$ via 
the following representations $\rho_0^{(1)}$ and $\rho_0^{(2)}$, respectively, 
and on the common target  via $\rho_1$: 
\begin{eqnarray*}
&&\rho_0^{(1)}=\alpha \oplus \beta \oplus \gamma, \qquad 
\rho_1=\alpha \oplus \beta^2 \oplus \alpha\beta \oplus \gamma, \\
&&\rho_0^{(2)}=\alpha \oplus \beta^2 \oplus \alpha\beta, 
\qquad (\alpha, \beta, \gamma) \in T.
\end{eqnarray*}
Hence the quotient Chern classes of universal maps for $f_1$ and $f_2$ are 
$$c(f_1)=1+(a+2b)+ab -ab^2, \quad 
c(f_2)=1+c \;\;\; \in H^*(BT),$$
where $a,b,c$ are the first Chern classes for 
standard representations $\alpha, \beta, \gamma$ of $\C^*$. 
Also Landweber-Novikov classes are 
\begin{eqnarray*}
&& s_0(f)=f_{1*}(1)+f_{2*}(1)=2(a + b)+c, \\
&& s_1(f)=f_{1*}(c_1(f_1))+f_{2*}(c_1(f_2))=2(a+b)(a+2b)+c^2,  
\end{eqnarray*}
and so on. 
Note that in the $xyz$-space, 
the $\overline{A_0^2}$-locus is the union of two planes $x=0$ and $z=0$, 
while in the $uvw$-space, 
the locus is just the crosscap $u^2v=w^2$. 
Then, 
the SSM class for $\overline{A_0^2}$ 
applied to the universal map $f_1$ is given by 
$$tp^{\SM}(\overline{A_0^2})(f_1)
=\frac{a}{1+a}+\frac{c}{1+c}-\frac{ac}{(1+a)(1+c)}$$
using the exclusion-inclusion of SSM classes: 
the plane $x=0$ plus the plane $z=0$ minus the $y$-axis 
(For the plane $x=0$, the corresponding 
normal Chern class is $1+a$, 
hence the SSM class in the ambient space is $a(1+a)^{-1}$). 
This is the restriction equation  at $A_1A_0$. 
The SSM class applied to the universal map $f_2$ is 
actually the target SSM class for the image of 
$$A_1: (x,y) \mapsto (u,v,w)=(x,y^2,xy).$$
Since we have already known that 
$$\textstyle tp_3^{\SM}(\jeden_{{\rm image}})
=\frac{1}{6}(s_0^3-7s_{01}+3s_0s_1+2s_2)$$
 (Proposition \ref{complexIM}), 
the restriction equation at $A_0A_1$ is obtained by 
$$tp_3^{\SM}(\overline{A_0^2})(f_2)=
tp_3^{\SM}(\jeden_{{\rm image}})(f_1)=(a+b)(4a^2+9ab+8b^2).$$

}
\end{exam}

\

\subsection{Discriminant Chern class  for stable maps} \label{Dis}

Let us consider the case of $m \ge n$ 
and the {\it discriminant} of proper stable maps $f: M \to N$
$$D(f) := \overline{f(A_1(f))}.$$
\begin{definition}{\rm 
We call 
$C_*(\jeden_{D(f)})\in H^*(N)$ {\it the discriminant Chern class of $f$}. 
}
\end{definition}

\begin{figure}
\includegraphics[clip, width=10cm]{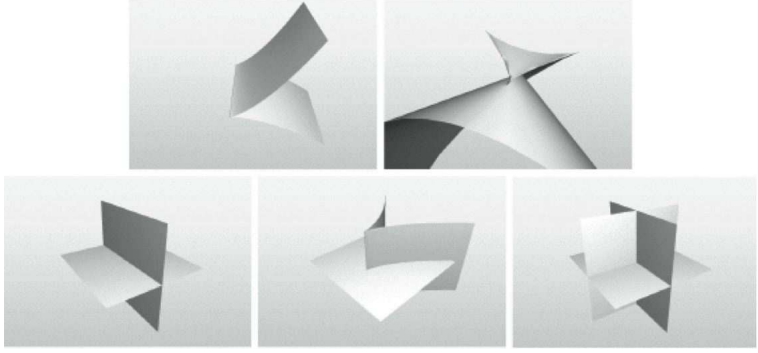}
\caption{\small Cuspidal edge ($A_2$), swallowtail ($A_3$) and 
stable multi-singularity loci in the target space of $3$-to-$3$ maps}
\label{image3-3} 
\end{figure}

For simplicity, we deal with the equidimensional case $m=n$ below. 
Stable singularities of codimension up to $3$ 
are $A_1$, $A_2$, $A_3$, $A_1^2$, $A_1A_2$, $A_2A_1$, $A_1^3$ 
(Fig. \ref{image3-3}).

The same procedure as in the case of image can be applied to 
the case of discriminant: 
There exists a constructible function on $M$ 
$$\textstyle \disconst
:=\jeden_{\overline{A_1}}-\frac{1}{2}\jeden_{\overline{A_1^2}}
-\frac{1}{6}\jeden_{\overline{A_1^3}}+\frac{1}{2}\jeden_{\overline{A_3}} 
+ \cdots \;\; \in \F(M)$$
so that 
$$f_*\disconst = \jeden_{D(f)}.$$
Since the local structures of $A_k$-singularities and 
$\K$-orbits in $\Sigma^2$ are well-understood, 
this constructible function can be explicitly written down up to a certain codimension. 
We can prove the following theorem: 

\begin{thm} \label{DisSSM}
There is a polynomial $tp^{\SM}(\disconst)$ 
in the quotient Chern class $c_i$ and the Landweber-Novkov class $s_I$ so that 
$$
C_*(\disconst)=c(TM)\cdot tp^{\SM}(\disconst) \; \in H^*(M) 
$$
for proper stable maps $f: M^n \to N^n$ in low dimension $(n <9)$.  
In fact, the low degree terms are given by 
\begin{eqnarray*}
tp^{\SM}(\disconst)
&=&{\textstyle c_1+\frac{1}{6}(6c_1^2+6c_2-3c_1s_1)}\\
&&{\textstyle +\frac{1}{6}\left(
\begin{array}{l}
{\textstyle c_1^3+11c_1c_2+6c_3-2c_1s_{01}-5c_1^2 s_1}\\
{\textstyle  -4 c_2s_1+c_1s_1^2+2c_1s_2 }
\end{array}
\right) + h.o.t}. 
\end{eqnarray*}
\end{thm} 

\begin{rem}{\rm 
We denote by $tp^{\SM}(\jeden_{D(f)})$ the universal Segre-SM class for 
the discriminant $D(f)$: 
\begin{eqnarray*}
tp^{\SM}(\jeden_{D(f)}) &:=& c(TN)^{-1} C_*(D(f))\\
&=& c(TN)^{-1}f_*C_*(\disconst) \\
&=& f_*(c(f)^{-1}\cdot tp^{\SM}(\disconst))\\
&=& 
\textstyle 
s_1+(s_{01}-\frac{1}{2}s_1^2) \\
& & 
\textstyle 
\;\; +(s_{001}-s_{01}s_1
+\frac{1}{6}s_1^3-\frac{1}{6}s_{11}+\frac{1}{6}s_3) + \cdots.
\end{eqnarray*}

}
\end{rem}

\begin{cor} 
The Euler characteristics of the discriminant of a proper stable map 
is universally expressed by 
$$
\chi(D(f))
=\int_M c(TM) \cdot  tp^{\SM}(\disconst)
=\int_N c(TN)\cdot tp^{\SM}(\jeden_{D(f)}).$$
\end{cor}

A reduced divisor $D$ in a complex manifold $N$ is called to be 
{\it free} (in the sense of Kyoji Saito) 
if the sheaf of germs of logarithmic vector fields ${\rm Der}_N(-\log D)$ 
is locally free. 
As for the CSM class of a free divisor $D$,  
the following equality was conjectured  by P. Aluffi, 
and was recently proved by X. Liao \cite{Liao}: 

\begin{thm}[{\bf CSM class of free divisors} \cite{Liao}] 
If $D$ is locally quasi-homogeneous (i.e., 
at each point, there is a weighted homogeneous defining equation  
in some local coordinates), it holds that 
\begin{center}
$c^{\SM}(N-D)=c({\rm Der}_N(-\log D))  \qquad  \in H^*(N)$
\end{center}
(in fact, the condition can be more weakened). 
\end{thm}

In our case, 
it is known that the discriminant $D(f)$ of a stable map $f: M \to N$ 
in Mather's nice dimension is  
a free divisor in $N$ which is locally quasi-homogeneous. 
We have seen that 
the CSM class of $D=D(f)$ in the ambient space $N$ is expressed  
using our target universal Segre-SM class: 
\begin{eqnarray*}
c^{\SM}(N-D)&=&c^{\SM}(N)-c^{\SM}(D(f))\\
&=&c(TN)(1-tp^{\SM}(\jeden_{D(f)})).
\end{eqnarray*}
Hence, the Chern class $c({\rm Der}_N(-\log D(f))$ 
is universally expressed in terms of $s_I$ and $c(TN)$. 
Namely, the meaning of 
our discriminant SSM class (written in $s_I$) becomes clearer: 

\begin{cor} 
The discriminant SSM class for proper stable maps 
is exactly the same as the quotient Chern class for the sheaf 
of  logarithmic vector fields and that of ambient vector fields of the target manifold, 
without the constant $1$: 
$$tp^{\SM}(\jeden_{D(f)})=1-c({\rm Der}_N(-\log D(f))-TN).$$
\end{cor}

\subsection{Generating function of multi-singularty SSM classes}  \label{multising_SSM} 
A better treatment of the universal SSM class for multi-singularities of stable maps 
may be as follows. 
This is due to a communication with M. Kazarian. 
Actually, this is parallel to his argument on multi-singularity Thom polynomials. 

Let $\underline{\eta}=(\eta_1, \cdots, \eta_r)$ be a multi-singularity type. 
Let $|{\rm Aut}(\underline{\eta}) |$ denote 
the number of permutations $\sigma \in \mathfrak{S}_r$ 
preserving the types of entries, $\eta_{\sigma(i)}=\eta_i \; (1\le i \le r)$,  
that is, if $\underline{\eta}$ consists of $k_i$ copies of mutually distinct mono-singularities, 
then $|Aut(\underline{\eta})|=k_1! \cdots k_s!$. 
Hence, in particular, 
$\deg_1 \underline{\eta} \cdot |{\rm Aut}(\eta_2, \cdots, \eta_r) | 
= |{\rm Aut}(\underline{\eta}) |$. 

For a stable map $f: M \to N$, 
let $M(\underline{\eta})(f)$ denote the closure of the locus of 
points $(x_1, \cdots, x_r, y) \in M^r\times N$ so that  
$f(x_1)=\cdots =f(x_r)=y$,  $x_i\not=x_j \; (i\not=j)$ 
and $f$ at $x_i$ is of type $\eta_i$. 
Put 
$$p_1: M^r \times N \to M, \;\; p': M^r \times N \to N$$ 
the projection 
to the first and the last factors, respectively. 
Then the {\it source and target multi-singularity constructible functions}  are defined by 
$$\textstyle 
\alpha_{\underline{\eta}}
:=p_{1*}\jeden_{M(\underline{\eta})(f)} \in \F(M), \quad 
\beta_{\underline{\eta}}
:=p'_*\jeden_{M(\underline{\eta})(f)} \in \F(N).$$
It holds that $f_*\alpha_{\underline{\eta}}=\beta_{\underline{\eta}}$. 

The supports of 
$\alpha_{\underline{\eta}}$ and $\beta_{\underline{\eta}}$  
are the  $\underline{\eta}$-singular locus $\overline{\underline{\eta}(f)} \subset M$ 
and its image $f(\overline{\underline{\eta}(f)})\subset N$, respectively: 
those functions take the values $|{\rm Aut}(\eta_2, \cdots, \eta_r) |$ 
and $|{\rm Aut}(\eta_1, \cdots, \eta_r) |$  
on the open parts of their supports, 
but  may take several different values on the boundary strata. 
The image constant function $\jeden_{f(M)}$  (resp. $\imageconst$) 
is written 
by a linear combination with rational coefficients 
of $\beta_{\underline{\xi}}$  (resp. $\alpha_{\underline{\xi}}$)
among multi-singularity types $\underline{\xi}$ adjacent to $\underline{\eta}$, 
for instance, 
$$\textstyle \jeden_{f(M)}=
\frac{1}{ |{\rm Aut}(\underline{\eta}) |} \cdot \beta_{\underline{\eta}} + 
\sum_{\mbox{\tiny boudary}} b_{\underline{\xi}} \cdot \beta_{\underline{\xi}}$$
for some $b_{\underline{\xi}} \in \Q$. 

We conjecture the existence of 
source and target universal Segre-SM classes for multi-singularity constructible functions, 
that  generalizes simultaneously Theorem \ref{tp^SM} on $tp^{\SM}$ for mono-singularities 
and Theorem \ref{kaz_thm} on $tp$ of multi-singularities. 
In some particular cases of low dimension, 
Thereoms \ref{imageSSM} and \ref{DisSSM} support that the conjecture is true. 

\begin{conj} 
For any stable multi-singularity type $\underline{\eta}$ 
in relative codimension $\kappa$, 
there exist power series 
$tp^{\SM}(\alpha_{\underline{\eta}})$ and 
$tp^{\SM}(\beta_{\underline{\eta}})$ 
 in quotient Chern classes $c_i\; (=c_i(f))$ and 
the Landweber-Novikov classes $s_I$ such that 
for any stable maps $f: M \to N$ of relative codimension $\kappa$ 
it holds that 
$$tp^{\SM}(\alpha_{\underline{\eta}})=c(TM)^{-1} C_*(\alpha_{\underline{\eta}}), \;\; 
tp^{\SM}(\beta_{\underline{\eta}})
=c(TN)^{-1}C_*(\beta_{\underline{\eta}})$$
in $H^*(M)$ and $H^*(N)$ respectively. 
\end{conj}

There two universal multi-singularity universal SSM classes are related in the following form 
by the naturality of $C_*$: We define 
$$\rho: H^*(M) \to H^*(N), \;\; \rho(\omega)=f_*(c(f)^{-1}\cdot \omega)$$
and then it holds that 
$$\rho ( tp^{\SM}(\alpha_{\underline{\eta}}))= tp^{\SM}(\beta_{\underline{\eta}}).$$

The conjecture implies a remarkable property that 
these universal series admit a very particular form; That is parallel to 
the argument on $tp$ for multi-singularities in \cite[\S 3]{Kaz03} and \cite[\S 2.6]{Kaz06}. 
For each stable multi-singularity type $\underline{\eta}$, 
let $R_{\underline{\eta}}$ be the polynomial in quotient Chern classes $c_i=c_i(f)$ 
so that 
$$tp^{\SM}(\alpha_{\underline{\eta}})=R_{\underline{\eta}} + 
\mbox{terms containing $f^*s_I$}.$$
We call $R_{\underline{\eta}}$ the {\it residual polynomial of $\underline{\eta}$}. 
Recall that the SSM class has a natural property for transverse pullback (Proposition \ref{segre}). 
Then the same argument as in \cite[\S 3]{Kaz03} shows that 
there is a universal recursive relation 
$$
tp^{\SM}(\alpha_{\underline{\eta}})
=
R_{\underline{\eta}} + 
\sum_I
R_{\underline{\eta}_I} f^*\rho(tp^{\SM}(\alpha_{\underline{\eta}_J})), 
$$
where the sum is taken over all proper subset 
$I \subset \{1, 2, \cdots, r\}$ containing the element $1$ 
and $J=[r]-I \not=\emptyset$. 
For example, 
{\small 
\begin{eqnarray*}
tp^{\SM}(\alpha_{\eta_1})
&=& R_{\eta_1} =tp^{\SM}(\overline{\eta_1}) 
\qquad  \mbox{(This is Theorem \ref{tp^SM})}\\
tp^{\SM}(\beta_{\eta_1})
&=& \rho(R_{\eta_1}) \\
tp^{\SM}(\alpha_{\eta_1, \eta_2})
&=& 
R_{\eta_1, \eta_2}+R_{\eta_1}\cdot \rho(R_{\eta_2}),\\
tp^{\SM}(\beta_{\eta_1, \eta_2})
&=&
\rho(R_{\eta_1, \eta_2})+\rho(R_{\eta_1})\cdot \rho(R_{\eta_2}),\\
tp^{\SM}(\alpha_{\eta_1, \eta_2, \eta_3})
&=& 
R_{\eta_1, \eta_2, \eta_3}
+R_{\eta_1, \eta_2}\cdot \rho(R_{\eta_3})
+R_{\eta_1, \eta_3}\cdot \rho(R_{\eta_2})\\
&&+R_{\eta_1}\cdot \rho(R_{\eta_2, \eta_3})
+R_{\eta_1}\cdot \rho(R_{\eta_2})\cdot \rho(R_{\eta_3}), \\
tp^{\SM}(\beta_{\eta_1, \eta_2, \eta_3})
&=& 
\rho(R_{\eta_1, \eta_2, \eta_3})
+\rho(R_{\eta_1, \eta_2})\cdot \rho(R_{\eta_3})
+\rho(R_{\eta_1, \eta_3}) \cdot \rho(R_{\eta_2})\\
&&+\rho(R_{\eta_1}) \cdot \rho(R_{\eta_2, \eta_3})
+\rho(R_{\eta_1}) \cdot \rho(R_{\eta_2})\cdot \rho(R_{\eta_3}). 
\end{eqnarray*}
}
In particular, 
this recursive relation provides 
an exponential generating function formula for those universal SSM classes. 
For a mono-singuarity type $\eta$, we take a distinguished variable $t_{\eta}$. 
For a multi-singularity type $\underline{\eta}=(\eta_1, \cdots, \eta_r)$, 
put $t^{\underline{\eta}}=t_{\eta_1} \cdots t_{\eta_r}$ 
(If we denote by $\xi_1^{k_1}\cdots \xi_s^{k_s}$ the entries in $\underline{\eta}$ 
(i.e., forgetting the order), 
then $t^{\underline{\eta}}=t_{\xi_1}^{k_1} \cdots t_{\xi_s}^{k_s}$ and 
$|Aut(\underline{\eta})|=k_1! \cdots k_s!$).  
Define 
the generating function of target Segre-SM classes of all stable multi-singularity types 
$$\mathcal{T}^{\SM}:=1+\sum_{\underline{\eta}} tp^{\SM}(\beta_{\underline{\eta}}) \cdot 
\frac{t^{\underline{\eta}}}{|Aut(\underline{\eta})|},  
$$
then by the above recursive relation we have 
$$
\mathcal{T}^{\SM}
=\exp \left( 
\sum_{\underline{\eta}} \rho(R_{\underline{\eta}}) \cdot 
\frac{t^{\underline{\eta}}}{|Aut(\underline{\eta})|}
\right).
$$

\subsection{Computing the image and discriminant Milnor numbers}\label{ImageMilnor}

We have seen in \S 5  
an application of Thom polynomials $tp$ to the problem on counting 
stable (multi-)singularities in generic deformation. 
Now we shall go on the same direction, 
but apply our higher Thom polynomial $tp^{\SM}$. 

\

\noindent
{\bf Image Milnor number}: 
Consider an $\A$-finitely determined weighted homogeneous map-germ 
$f: \C^m, 0 \to \C^{m+1}, 0$
which is not equivalent to any trivial unfolding of map-germ of smaller dimensions. 
Take a stable unfolding $F$ of $f$: 
$$\begin{array}{ccc}
\C^m &\stackrel{f}{\longrightarrow}& \C^{m+1}\\
\mbox{\footnotesize $i_0$} \downarrow \;\;\;\; & & \downarrow \mbox{\footnotesize $\iota_0$}  \\
\C^{m+k} &\stackrel{F}{\longrightarrow}&   \C^{(m+1)+k} 
\end{array}
$$
The image hypersurfaces of $f$ and $F$ 
relate  as $\Image(f)=\iota_0^{-1}(\Image(F))$.  
Take a generic (non-equivariant) section $\iota_t$, which yields 
a stable perturbation $f_t$ of $f_0=f$. 
Our interest is to compute the {\it vanishing Euler characteristics} of the section. 

\begin{definition}
$\mu_I(f):=(-1)^m(\chi(\Image(f_t))-1).$ 
\end{definition}

It was shown by D. Mond \cite{MararMond, Mond87, MondLect} that  
the singular Milnor fiber $\Image(f_t)$ has the homotopy type of a wedge of $m$-spheres, 
so the vanishing Euler number $\mu_I(f)$ is equal to the middle Betti number 
of the singular Milnor fiber, called the {\it image Milnor number} of $f$. 
In case of $m=1,2$, it is proved that 
$$\mu_I(f) \ge \A_e\mbox{\rm -codim}(f)$$
and the equality holds if $f$ is weighted homogeneous. 
The {\it Mond conjecture} claims that 
the same is true for any $m$ for which the pair $(m,m+1)$ is in Mather's nice dimensions, 
that has been unproven yet. 

Not only the image $\Image(f_t)$ 
but also the $k$-th multiple point locus $f_t(\overline{A_0^k(f_t)})$ in target has 
the same property about the homotopy type: 
 The {\it $k$-th image Milnor number} $\mu_{I_k}$ of $f$ is defined 
 in Houston \cite{Houston} 
 (of course, $\mu_I=\mu_{I_1}$).

Our strategy is the same as in \S 5: Using the natural torus action, 
 we deal with a global setting of universal maps associated to the above diagram 
of map-germs: we have the diagram of  
universal maps over $BT=\Proj^N\; (N\gg0)$ 
where $T=\C^*$: 
$$\begin{array}{ccc}
E_0 &\stackrel{f}{\longrightarrow}& E_1 \\
\mbox{\footnotesize $i_0$} \downarrow \;\;\;\; & & \;\;  \downarrow \mbox{\footnotesize $\iota_0$} \\
E_0\oplus E' &\stackrel{F}{\longrightarrow}&  E_1\oplus E' 
\end{array}
$$
Put $M = E_0$, $N=E_1$ the total spaces of source and target of the universal map 
for the original germ. 
A perturbation $\iota_t$ of $\iota_0$ is transverse to the image variety of 
the universal stable map $F$, which produces a stable perturbation $f_t: M \to N$. 

By Proposition \ref{segre} 
(the property of our Segre-SM class  for transversal pullback) and $\iota_0^*=\iota_t^*$, 
$$tp^{\SM}(\jeden_{f_t(M)})=\iota_0^*tp^{\SM}(\Image(F))
$$
which is thought of as the specialization of $tp^{\SM}(\Image(F))$ via $\iota_0$. 
Note that 
$$c(F)=c(f)=c(E_1-E_0) \in H^*(BT)\; (=H^*(E_0)=H^*(E_1)).$$
Then Theorem \ref{imageSSM} (or Remark \ref{rem_imageSSM}) shows that 
by the naturality of CSM classes 
$$c(E_1) \cdot tp^{\SM}(\jeden_{f_t(M)})=f_*(c(E_0)\cdot tp^{\SM}(\imageconst)).$$
On the other hand, the general slice of 
the image variety $f_t(M)$ via a fiber of the projection $N=E_1\to BT$ 
is isomorphic to $\Image(f_t) \subset \C^n$ in the local setting and 
$$\chi(\Image(f_t))=\int_{\C^n} \jeden_{{\rm Im}(f_t)} = \int_{\C^m} \imageconst(f_t).$$
By a similar argument of the proof of (2) in Theorem \ref{segre2}, 
we see that 
the $n$-dimensional component (some multiple of $a^n$) 
$$[c(E_1)\cdot tp^{\SM}(\jeden_{f_t(M)})]_n \in H^{2n}(BT)$$ 
is equal to the top Chern class $c_n(E_1)$ multiplied by the Euler number  
$\chi(\Image(f_t))$. 
In fact, the above arguments can properly be  stated in the $T$-equivariant setting: 
we then appeal to the Verdier specialization via $\iota_0$ and 
the Atiyah-Bott localization to the fixed point $0$ of $T$-equivariant CSM classes 
$C_*^T(\imageconst(F))$ and $C_*^T(\jeden_{\mbox{\tiny Image}(F)})$. 

Consequently, we have 

\begin{thm}
The following formula holds: 
$$ 
\chi(\Image(f_t))
=  \frac{[c(E_1)\cdot tp^{\SM}(\jeden_{f_t(M)})]_n}{c_n(E_1)} 
=  \frac{[c(E_0)\cdot tp^{\SM}(\imageconst)]_m}{c_m(E_0)},  
$$
where the notation in numerators $[\omega]_n$ means  
the coefficient of $a^n$ in $\omega \in H^*(BT)=\Q[[a]]$, 
and the denominators mean the products of weights and degrees: 
$c_m(E_0)=w_1\cdots w_m a^m$ and $c_n(E_1)=d_1\cdots d_n a^n$. 
In particular,  this formula enables us to compute the image Milnor number $\mu_I(f_0)$. 
\end{thm}

Notice that our formula above is valid 
for weighted homogeneous $\A$-finite germs with {\it any corank}. 
Comparing the above theorem 
with Thereom \ref{counting_0_stable_inv}, 
their similarity is clear. 

In the following examples, we compute the image Milnor number $m$-to-$(m+1)$ map-germs. 
Recall that for stable maps in relative codimension one, 
there is  a unique universal  Segre-SM class 
$tp^{\SM}(\imageconst)$ for the image of maps (Theorem \ref{imageSSM}).

\begin{exam}{\rm 
$(m,n)=(2,3)$: 
For weighted homogeneous map-germs $\C^2,0 \to \C^3, 0$, 
$$ 
c(f_\eta)=\frac{(1+d_1a)(1+d_2a)(1+d_3a)}{(1+w_1a)(1+w_2a)}, $$
$$s_0=f_{\eta*}(1)=\frac{d_1d_2d_3}{w_1w_2}a, \;\; s_I=f_{\eta*}(c^I)=c^I s_0,$$ 
$$
C^T_*(\imageconst)
=(1+w_1a)(1+w_2a) \cdot tp^{\SM}(\imageconst)(f_0), $$
$$ c_{top}(E_0)=w_1w_2a^2.  $$
Our computation on $\mu_I$ is straightforward like Example \ref{counting23}. 
We have the following result, which 
completely coincides with D. Mond's computation \cite{Mond}, 
the methods are quite different, though.  

\begin{eqnarray*}
&&\mu_I=\textstyle
-1+\left[\frac{1}{w_1w_2}(1+w_1a)(1+w_2a) \cdot tp^{\SM}(\imageconst)(f_0)\right]_2\\
&&\textstyle
=\frac{1}{6w_1^3w_2^3}\left(d_1^2(d_2^2d_3^2-w_1^2w_2^2)-w_1^2w_2^2
\{d_2^2+d_3^2+5w_1^2\right.\\
&&  +9w_1w_2+5w_2^2-6d_3(w_1+w_2)+3d_2(d_3-2(w_1+w_2)\}\\
&& \left.-3d_1w_1w_2\{w_1w_2(d_3-2(w_1+w_2))+d_2(w_1w_2+d_3(w_1+w_2))\}\right).\\
\end{eqnarray*}
}
\end{exam}

\begin{exam}{\rm 
$(m,n)=(3,4)$: 
For weighted homogeneous map-germs $\C^3,0 \to \C^4, 0$, 
 the image Milnor numbers $\mu_I$ and $\mu_{I_2}$  are given 
in the following Tables \ref{table_image34} and \ref{table_image34_2}. 

\begin{table}[b]
{\small 
\begin{eqnarray*}
&& 
\mu_I=1-\left[\frac{(1+w_1a)(1+w_2a)(1+w_3a)}{w_1w_2w_3}\cdot tp^{\SM}(\imageconst)(f_0)\right]_3\\
&&\textstyle
 =\frac{1}{24 w_1^4 w_2^4 w_3^4}
(d_1^3 (d_2^3 d_3^3 d_4^3+2 d_2^2 d_3^2d_4^2 w_1 w_2 w_3-d_2 d_3 d_4 w_1^2 w_2^2 w_3^2-2 w_1^3 w_2^3 w_3^3)
\\&&
+2d_1^2 w_1 w_2 w_3 (d_2^3 d_3^2 d_4^2+2 (d_3+d_4) w_1^2 w_2^2 w_3^2+d_2w_1 w_2 w_3\\&&
 (-9 d_3^2 d_4+2 w_1 w_2 w_3
+9 d_3 d_4 (-d_4+w_1+w_2+w_3))\\&&
+d_2^2d_3 d_4 (d_3^2 d_4-9 w_1 w_2 w_3+d_3 d_4 (d_4-3 (w_1+w_2+w_3))))\\&&
+2w_1^3 w_2^3 w_3^3 (-d_2^3-d_3^3+2 d_3^2 d_4-d_4^3
+2 d_2^2 (d_3+d_4)+d_4w_1^2-9 d_4 w_1 w_2\\&&
+9 w_1^2 w_2+d_4 w_2^2+9 w_1 w_2^2-9 d_4 w_1 w_3
+9w_1^2 w_3-9 d_4 w_2 w_3\\&&
+27 w_1 w_2 w_3+9 w_2^2 w_3+d_4 w_3^2+9 w_1w_3^2+9 w_2 w_3^2\\&&
+d_3 (2 d_4^2+w_1^2+w_2^2-9 w_2w_3+w_3^2-9 w_1 (w_2+w_3)
-3d_4 (w_1+w_2+w_3))\\&&
+d_2 (2 d_3^2+2 d_4^2+w_1^2-9 w_1 w_2+w_2^2-9 w_1w_3-9 w_2 w_3+w_3^2\\&&
-3 d_4 (w_1+w_2+w_3)+d_3 (9 d_4-3 (w_1+w_2+w_3))))\\&&
-d_1w_1^2 w_2^2 w_3^2 (d_2^3 d_3 d_4
+2 d_2^2 (9 d_3^2 d_4+9 d_3 d_4 (d_4-w_1-w_2-w_3)\\&&
-2w_1 w_2 w_3)-2 w_1 w_2 w_3 (2 d_3^2+2 d_4^2+w_1^2-9 w_1 w_2+w_2^2-9w_1 w_3\\&&
-9 w_2 w_3+w_3^2-3 d_4 (w_1+w_2+w_3)+d_3 (9 d_4-3 (w_1+w_2+w_3)))\\&&
+d_2(d_3^3 d_4+18 d_3^2 d_4 (d_4-w_1-w_2-w_3)+6 w_1 w_2 w_3\\&&
 (-3 d_4+w_1+w_2+w_3)+d_3(d_4^3-18 w_1 w_2 w_3-18 d_4^2(w_1+w_2+w_3)\\&&
+d_4 (17 w_1^2+17 w_2^2+6w_2 w_3+17 w_3^2+6 w_1 (w_2+w_3))))))
\end{eqnarray*}
}
\caption{\small Image Milnor numbers for $\C^3,0 \to \C^4,0$. }
\label{table_image34}
\end{table}

\begin{table}
{\small 
\begin{eqnarray*}
&&
\mu_{I_2}=
1-\left[\frac{(1+w_1a)(1+w_2a)(1+w_3a)}{w_1w_2w_3}\cdot tp^{\SM}(\imageconst(2))(f_0)\right]_3\\
&&\textstyle
=\frac{1}{24 w_1^4 w_2^4 w_3^4}
(d_1^3 (3 d_2^3 d_3^3 d_4^3-2 d_2^2 d_3^2d_4^2 w_1 w_2 w_3-3 d_2 d_3 d_4 w_1^2 w_2^2 w_3^2\\&&
+2 w_1^3 w_2^3 w_3^3)+2w_1^3 w_2^3 w_3^3 (d_2^3+d_3^3+d_4^3-24 d_4^2 w_1+47 d_4 w_1^2-24 w_1^3\\&&
-24d_4^2 w_2+57 d_4 w_1 w_2-33 w_1^2 w_2+47 d_4 w_2^2-33 w_1 w_2^2-24 w_2^3\\&&
-24d_4^2 w_3+57 d_4 w_1 w_3-33 w_1^2 w_3+57 d_4 w_2 w_3-51 w_1 w_2 w_3-33w_2^2 w_3\\&&
+47 d_4 w_3^2-33 w_1 w_3^2-33 w_2 w_3^2-24 w_3^3\\&&
+d_3^2 (22 d_4-24 (w_1+w_2+w_3))+d_2^2(22 d_3+22 d_4-24 (w_1+w_2+w_3))\\&&
+d_3 (22 d_4^2+47 w_1^2+47 w_2^2+57 w_2 w_3+47w_3^2+57 w_1 (w_2+w_3)-69 d_4 (w_1\\&&
+w_2+w_3))+d_2 (22 d_3^2+22 d_4^2+47w_1^2+57 w_1 w_2+47 w_2^2+57 w_1 w_3\\&&
+57 w_2 w_3+47 w_3^2-69 d_4 (w_1+w_2+w_3)+d_3(75 d_4-69 (w_1+w_2+w_3))))\\&&
-2 d_1^2 w_1 w_2 w_3 (d_2^3 d_3^2 d_4^2+2w_1^2 w_2^2 w_3^2 (-11 d_3-11 d_4+12 (w_1+w_2+w_3))\\&&
-d_2 w_1 w_2 w_3 (-21d_3^2 d_4+22 w_1 w_2 w_3-3 d_3 d_4 (7 d_4-9 (w_1+w_2+w_3)))\\&&
+d_2^2 d_3d_4 (d_3^2 d_4+21 w_1 w_2 w_3+d_3 d_4 (d_4+3 (w_1+w_2+w_3))))+d_1w_1^2 w_2^2 w_3^2\\&&
 (-3 d_2^3 d_3 d_4+2 w_1 w_2 w_3 (22 d_3^2+22 d_4^2+47w_1^2+57 w_1 w_2+47 w_2^2+57 w_1 w_3\\&&
 +57 w_2 w_3+47 w_3^2-69 d_4 (w_1+w_2+w_3)+d_3(75 d_4-69 (w_1+w_2+w_3)))\\&&
 +d_2^2 (-42 d_3^2 d_4+44 w_1 w_2 w_3-6 d_3d_4 (7 d_4-9 (w_1+w_2+w_3)))\\&&
 -3 d_2 (d_3^3 d_4+2 d_3^2 d_4 (7 d_4-9 (w_1+w_2+w_3))+2w_1 w_2 w_3 (-25 d_4\\&&
 +23 (w_1+w_2 +w_3))+d_3 (d_4^3-50 w_1 w_2 w_3-18d_4^2 (w_1+w_2+w_3)\\&&
 +d_4 (17 w_1^2+17 w_2^2 +18 w_2 w_3+17 w_3^2+18 w_1 (w_2+w_3))))))
\end{eqnarray*}
}
\caption{\small Second image Milnor numbers for $\C^3,0 \to \C^4,0$. }
\label{table_image34_2}
\end{table}

For corank one map-germs $\C^3,0\to \C^4,0$, 
take weights $w_0, w_1, w_2$ and degrees 
$d_1, d_2, d_3=w_1, d_4=w_2$, 
then we obtain a new general formula for corank one germs: 

$$
\mu_I=\frac{(w_0-d_1) (w_0-d_2)}{24 w_0^4 w_1 w_2} \left(
\begin{array}{l}
d_1^2 \left(d_2^2+3 d_2 w_0+2 w_0^2\right)\\
+d_1 w_0 \left(3 d_2^2-d_2 (19  w_0+4 (w_1+w_2))\right. \\
\left.+2 w_0(w_0-2 (w_1+w_2))\right)\\
+2 w_0^2 \left(d_2^2+d_2 (w_0-2 (w_1+w_2))\right.\\
\left.+2(5 w_0 (w_1+w_2)+3 w_1 w_2)\right)\\
   \end{array}
   \right).
$$

The classification of $\A$-simple germs of corank one can be seen 
in  \cite{Houston}, and it is checked that 
for weighted homogeneous germs appearing in the list, 
our formulas above recover the same answers on image Milnor numbers 
as computed in \cite{Houston}. 
For instance, 
$$Q_k: (x, y, xz+yz^2, z^3+y^kz)$$
 has 
weights $(k, 2, k+2)$ and degrees $(k+2, 2, 2k+2, 3k)$, and 
the above formula gives $\mu_I=k$ and $\mu_{I_2}=0$. 

Some examples of corank $2$ germs of $\C^3,0 \to \C^4,0$ 
are recently considered in  \cite{Altintas} 
in a completely different approach. 
It would be nice to compare the computations: 
As a test, let us take 
\begin{eqnarray*}
&&\Hat{A}_k: (x, y^k+xz+x^{2k-2}y, yz, z^2+y^{2k-1})\\
&&\Hat{B}_{2k+1}: (x, y^2+xz, x^2+xy, y^{2k+1}+y^{2k-1}z^2+z^{2k+1}). 
\end{eqnarray*}
Those are  $\A$-finite germs, and 
 weights and degrees are  
 $(1,2,2k-1)$ and  $(1,2k, 2k+1, 2(2k-1))$  for $\Hat{A}_k$, and 
  $(1,1,1)$ and  $(1,2,2,2k+1)$ for $\Hat{B}_{2k+1}$. 
Our formula gives the  answer in Table \ref{imageMcorank2}. 

For another example, 
$$(x^2+z^\ell y, y^2-z^\ell x, x^3+x^2y+x y^2-y^3, z)$$
 we have $\mu_I=45\ell-12$ 
 which coincides with \cite[Prop.4.4, 4.6]{Altintas}. 

 \begin{table}
$$
\begin{array}{c | c | c} 
\mbox{\footnotesize type} & \mbox{$\mu_I$} & k=2,3,4,5,6, \cdots \\
\hline 
\mbox{\tiny $\Hat{A}_k$} & \frac{1}{3}k(-3+15k-20k^2+6k^3+2k^4) &
\mbox{\footnotesize $18, 186, 844, 2620, 6510,\cdots$}\\
\mbox{\tiny $\Hat{B}_{2k+1}$} & 3k^2(1+10k) &
\mbox{\footnotesize $252,837,1968,3825,6588,\cdots$}\\ 
\end{array}
$$
\caption{}\label{imageMcorank2}
\end{table}
}
\end{exam}

\noindent
{\bf Discriminant Milnor number}: 
Next, let us consider $f: \C^m, 0 \to \C^n, 0$ in case of $m \ge n$. 
Assume that $f$ is $\A$-finitely determined. 
In the same way as above, 
we set the vanishing Euler characteristics: 

\begin{definition}
 $\mu_\Delta(f_0):=(-1)^{n-1}(\chi(D(f_t))-1).$
 \end{definition}

It is shown by Damon-Mond \cite{DM} that 
the discriminant $D(f_t)$ of a stable perturbation 
has the homotopy type of a wedge of $(n-1)$-spheres, 
so the vanishing Euler number $\mu_\Delta(f)$ 
is equal to the middle Betti number of $D(f_t)$, 
called the {\it discriminant Milnor number}. 
It is proved in \cite{DM} that  if $(m,n)$ is in nice dimensions, 
$$\mu_\Delta(f) \ge \A_e\mbox{\rm -codim}(f)$$
and the equality holds if $f$ is weighted homogeneous.

For a finitely determined weighted homogeneous germ $f$,  
we compute $\mu_\Delta(f)$ by  
localizing our higher Thom polynomials: 

\begin{thm} It holds that 
$$
 \chi(D(f_t)) =  \frac{[c(E_1)\cdot tp^{\SM}(\jeden_{D(f_t)})]_n}{c_n(E_1)}
 = \frac{[c(E_0)\cdot tp^{\SM}(\disconst)]_m}{c_m(E_0)}. 
$$
Thus we can compute the discriminant Milnor number $\mu_\Delta(f_0)$ 
in terms of weights and degrees. 
\end{thm}

Recall the {\it discriminant} Segre-SM class 
$tp^{\SM}(\disconst)$ for $m$-to-$m$ maps is given in 
Theorem \ref{DisSSM}. We use the low degree terms of this power series 
for the study of vanishing topology of 
germs  $\C^m,0 \to \C^m,0$, $m=2,3$.

\begin{exam}{\rm 
$(m,n)=(2,2)$: 
For weighted homogeneous map-germs $\C^2,0 \to \C^2, 0$, 
we recover the computational result in Gaffney-Mond  \cite{GM1} 
in a completely different way. 

\begin{eqnarray*}
&\mu_\Delta&=\textstyle
1-\left[\frac{1}{w_1w_2}(1+w_1a)(1+w_2a) \cdot tp^{\SM}(\disconst)(f_0)\right]_2\\
&&=
\frac{1}{2w_1^2w_2^2}(d_1d_2-2w_1w_2) \\
&&\quad (d_1^2+d_2^2+w_1^2+2d_1(d_2-w_1-w_2)+w_2^2-2d_2(w_1+w_2))
\end{eqnarray*}

}
\end{exam}

\begin{exam}{\rm 
$(m,n)=(3,3)$: 
For discriminant Milnor number of finitely-determined 
weighted homogeneous finite germs $\C^3,0 \to \C^3, 0$, 
we have the following formula in Table \ref{table_dis33}. 
In particular, 
for corank one map-germs, 

$$
\mu_\Delta=\frac{d-2 w_0}{6 w_0^3 w_1 w_2}
 \left(
 \begin{array}{l}
 d^4-4 d^3 w_0+d^2 w_0 (8 w_0-3 (w_1+w_2))\\
 +2 d w_0^2(3 (w_1+w_2)-4 w_0)\\
 +3 w_0^2 \left(w_0^2-w_0 (w_1+w_2)+2 w_1w_2\right)
 \end{array}
   \right). 
 $$
 
 This general formula also seems to be new.  It can be checked that 
this agrees with known computational results 
for weighted homogeneous germs appearing in $\A$-classification, e.g. \cite{MT}. 
\begin{table}
{\small 
\begin{eqnarray*}
&& 
\mu_\Delta=-1+\left[\frac{1}{w_1w_2w_3}(1+w_1a)(1+w_2a)(1+w_3a)\cdot tp^{\SM}(\disconst)(f_0)\right]_3\\
&&\textstyle
=\frac{1}{6 w_1^3 w_2^3 w_3^3}
(d_1^5 d_2^2 d_3^2+3 d_1^4 d_2 d_3 (d_2^2d_3+d_2 d_3 (d_3-w_1-w_2-w_3)\\&&
-w_1 w_2 w_3)+w_1^2 w_2^2 w_3^2(d_2^3+d_3^3-6 w_1^3-7 w_1^2 w_2-7 w_1 w_2^2-6 w_2^3-7 w_1^2 w_3\\&&
-15 w_1 w_2w_3-7 w_2^2 w_3-7 w_1 w_3^2-7 w_2 w_3^2-6 w_3^3-8 d_3^2 (w_1+w_2+w_3)\\&&
+d_3(13 w_1^2+13 w_2^2+15 w_2 w_3+13 w_3^2+15 w_1 (w_2+w_3))\\&&
+2 d_2^2 (7 d_3-4(w_1+w_2+w_3))+d_2 (14 d_3^2+13 w_1^2+13 w_2^2+15 w_2 w_3\\&&
+13 w_3^2+15 w_1(w_2+w_3)-27 d_3 (w_1+w_2+w_3)))\\&&
+d_1^3 (3 d_2^4 d_3^2+6 d_2^3 d_3^2(d_3-w_1-w_2-w_3)+w_1^2 w_2^2 w_3^2-3 d_2 d_3 w_1 w_2 w_3\\&&
 (5 d_3-4(w_1+w_2+w_3))+3 d_2^2 d_3 (d_3^3-5 w_1 w_2 w_3-2 d_3^2 (w_1+w_2+w_3)\\&&
+d_3(w_1+w_2+w_3)^2))+d_1^2 (d_2^5 d_3^2+3 d_2^4 d_3^2 (d_3-w_1-w_2-w_3)\\&&
-2w_1^2 w_2^2 w_3^2 (-7 d_3+4 (w_1+w_2+w_3))\\&&
+3 d_2^3 d_3 (d_3^3-5 w_1 w_2w_3-2 d_3^2 (w_1+w_2+w_3)+d_3 (w_1+w_2+w_3)^2)\\&&
-d_2 w_1 w_2 w_3(15 d_3^3-14 w_1 w_2 w_3-30 d_3^2 (w_1+w_2+w_3)\\&&
+3 d_3 (5 w_1^2+5 w_2^2+8w_2 w_3+5 w_3^2+8 w_1 (w_2+w_3)))\\&&
+d_2^2 d_3 (d_3^4-3 d_3^3 (w_1+w_2+w_3)+30w_1 w_2 w_3 (w_1+w_2+w_3)\\&&
+3 d_3^2 (w_1+w_2+w_3)^2-d_3 (w_1^3+3 w_1^2(w_2+w_3)+(w_2+w_3)^3\\&&
+3 w_1 (w_2^2+14 w_2 w_3+w_3^2))))+d_1w_1 w_2 w_3 (-3 d_2^4 d_3-3 d_2^3 d_3\\&&
 (5 d_3-4 (w_1+w_2+w_3))+w_1w_2w_3 (14 d_3^2+13 w_1^2+13 w_2^2+15 w_2 w_3\\&&
 +13 w_3^2+15 w_1 (w_2+w_3)-27 d_3(w_1+w_2+w_3))-d_2^2 (15 d_3^3-14 w_1 w_2 w_3\\&&
 -30 d_3^2 (w_1+w_2+w_3)+3d_3 (5 w_1^2+5 w_2^2+8 w_2 w_3+5 w_3^2 +8 w_1 (w_2+w_3))) \\&&
 -3 d_2 (d_3^4-4d_3^3 (w_1+w_2+w_3)+9 w_1 w_2 w_3 (w_1+w_2+w_3)\\&&
 +d_3^2 (5 w_1^2+5 w_2^2+8w_2 w_3+5 w_3^2+8 w_1 (w_2+w_3))\\&&
 -d_3 (2 w_1^3+4 w_1^2 (w_2+w_3)
 +w_1(4 w_2^2+21 w_2 w_3+4 w_3^2)+2 (w_2^3+2 w_2^2 w_3\\&&
 +2 w_2 w_3^2+w_3^3)))))
\end{eqnarray*}
}
\caption{\small Discriminant Milnor number for germs $\C^3,0 \to \C^3,0$}
\label{table_dis33}
\end{table}

As examples of corank $2$ singularity types,  
for $(x^2+y^2+xz, xy, z)$, $\mu_\Delta=1$, 
and for $(x^9+y^2+xz, xy, z)$, $\mu_\Delta=183$.

}
\end{exam}


\end{document}